\documentclass[9pt,twoside]{article}
\usepackage{amssymb}
\usepackage{makeidx}
\usepackage{color}
\usepackage{latexsym}
\usepackage{epic}
\usepackage{graphicx}
\usepackage{ifpdf}
\makeatletter

\@addtoreset{equation}{section} \makeatother

\newtheorem{theorem}{Theorem}[section]
\newtheorem{remark}[theorem]{Remark}
\newtheorem{lemma}[theorem]{Lemma}
\newtheorem{example}[theorem]{Example}
\newtheorem{proposition}[theorem]{Proposition}
\newtheorem{corollary}[theorem]{Corollary}
\newtheorem{definition}[theorem]{Definition}
\newtheorem{convention}[theorem]{Convention}

\newcommand{\ras}{evenly pairing boundary}

\newcommand{\cvd}{\ \rule{0.5em}{0.5em}}

\newcommand{\be}{\begin{equation}}
\newcommand{\ee}{\end{equation}}

\newcommand{\ben}{\begin{enumerate}}
\newcommand{\een}{\end{enumerate}}
\newcommand{\bit}{\begin{itemize}}
\newcommand{\eit}{\end{itemize}}
\newcommand{\edoc}{\end{document}}

\newcommand{\bdefi}{\begin{definition}}
\newcommand{\btheo}{\begin{theorem}}
\newcommand{\bprop}{\begin{proposition}}
\newcommand{\brema}{\begin{remark}}
\newcommand{\bcoro}{\begin{corollary}}
\newcommand{\blemm}{\begin{lemma}}
\newcommand{\bexam}{\begin{example}}
\newcommand{\bconv}{\begin{convention}}

\newcommand{\edefi}{\end{definition}}
\newcommand{\etheo}{\end{theorem}}
\newcommand{\eprop}{\end{proposition}}
\newcommand{\erema}{\end{remark}}
\newcommand{\ecoro}{\end{corollary}}
\newcommand{\elemm}{\end{lemma}}
\newcommand{\eexam}{\end{example}}
\newcommand{\econv}{\end{convention}}

\newcommand{\N}{{\mathbb N}}

\newcommand{\R}{{\mathbb R}}

\newcommand{\dqm}{d_Q^+}
\newcommand{\F}{F^{{\rm rev}}}
\newcommand{\bcaum}{\partial_C^+M}
\newcommand{\bcaumr}{\partial_C^-M}
\newcommand{\bcaums}{\partial_C^sM}
\newcommand{\caum}{ M_C^+}
\newcommand{\caumr}{ M_C^-}
\newcommand{\caums}{ M_C^s}
\newcommand{\cauchyd}{{\rm Cau}(M)} 
\newcommand{\cauchydr}{{\rm Cau}^{-}(M)} 
\newcommand{\rd}{d^{{\rm rev}}}
\newcommand{\gromp}{M_G^+}
\newcommand{\gromm}{M_G^-}

\newcommand{\dpl}{d^+}
\newcommand{\dm}{d^-}
\newcommand{\ds}{d^s}
\newcommand{\lip}{{\cal L}_1(M,d^R)}

\newcommand{\lipdM}{{\cal L}^{\rm Max}_1(M,d)}
\newcommand{\lipd}{{\cal L}_1(M,d)}

\newcommand{\cur}{C(M)}

\newcommand{\bpbus}{\partial_{\cal B}M}
\newcommand{\bbus}{\partial_BM}
\newcommand{\bus}{M_B}

\newcommand{\dri}{d^R}

\newcommand{\bpgro}{\partial_{\cal G}M}

\newcommand{\bgro}{\partial_GM}
\newcommand{\cambios}{{}}

\def\br#1\er{\textcolor{red}{#1}} %

\begin{document}

\parindent=5mm
\date{}

\medskip


\title{{\bf\LARGE Gromov, Cauchy and causal boundaries for Riemannian, Finslerian and
Lorentzian manifolds}}

\author{{\bf\large J.L. Flores$^*$,
J. Herrera$^*$,
M. S\'anchez$^\dagger$}\\
{\it\small $^*$Departamento de \'Algebra, Geometr\'{\i}a y Topolog\'{\i}a,}\\
{\it \small Facultad de Ciencias, Universidad de M\'alaga,}\\
{\it\small Campus Teatinos, E-29071 M\'alaga, Spain}\\
{\it\small $^\dagger$Departamento de Geometr\'{\i}a y Topolog\'{\i}a,}\\
{\it\small Facultad de Ciencias, Universidad de Granada,}\\
{\it\small Campus Fuentenueva s/n, E-18071 Granada, Spain}}

\maketitle

\begin{abstract}

Recently, the old notion of {\em causal boundary} for a spacetime
$V$ has been redefined consistently. The computation of this
boundary $\partial V$ on any standard conformally stationary
spacetime $V=\R\times M$, suggests a natural compactification
$M_B$ associated to any Riemannian metric on $M$
  or, more generally, to any Finslerian one. The
corresponding boundary $\partial_BM$  is constructed in terms of
Busemann-type functions. Roughly, $\partial_BM$
 represents the set of all the directions  in $M$ including both,
 asymptotic and ``finite'' (or ``incomplete'')  directions.

This Busemann boundary $\partial_BM$ is related to two classical
boundaries: the Cauchy boundary $\partial_{C}M$ and the Gromov
boundary $\partial_GM$. In a natural way $\partial_CM\subset
\partial_BM\subset \partial_GM$, but the topology in $\partial_BM$
is coarser than the others. Strict coarseness reveals some
remarkable possibilities
---in the Riemannian case, either $\partial_CM$ is not locally
compact or $\partial_GM$ contains points which cannot be reached
as limits of ray-like curves in $M$.

In the non-reversible Finslerian case, there exists always a
second boundary associated to the reverse metric, and many
additional subtleties appear. The spacetime viewpoint interprets
the asymmetries between the two Busemann boundaries,
$\partial^+_BM (\equiv
\partial_BM)$, $\partial^-_BM$, and this yields natural relations between
some of their points.

Our aims are: (1) to study the subtleties of both, the Cauchy
boundary for any {\em generalized} (possibly non-symmetric) {\em
distance} and the Gromov compactification for any  (possibly
incomplete) Finsler manifold, (2) to introduce the new Busemann
compactification $M_B$, relating it with the previous two
completions, and (3) to give  a full description of the causal
boundary $\partial V$ of any standard conformally stationary
spacetime.
\end{abstract}

\tableofcontents

\newpage

\section{Introduction} 


In Differential Geometry there are quite a few boundaries which
can be attached to a space, depending on the problem one would
like to study. For a Riemannian manifold $(M,g)$, when $g$ is
incomplete the Cauchy completion yields a simple boundary. When
$g$ is complete, Gromov introduced a general compactification by
using quotients of Lipschitz functions \cite{Gr81}. Such a
compactification coincides with Eberlein and O'Neill's one for a
Hadamard manifold, which defines the boundary points as  Busemann
functions associated to rays, up to additive constants,  and uses
the {\em cone topology} \cite{EO}. This construction can be
extended to more general spaces (as the CAT(0) ones), but it was
not conceived for an arbitrary Riemannian manifold. Among the
different boundaries in Lorentzian Geometry (Schmidt's {\em bundle
boundary}, Geroch's {\em geodesic boundary}, Penrose's {\em
conformal boundary}...), the so-called {\em causal boundary} (or
{\em c-boundary} for short) becomes especially interesting. The
c-boundary was introduced by Geroch, Kronheimer and Penrose
\cite{GKP} by using a conformally invariant construction, which is
explicitly intrinsic and systematic  (in advantage with the
conformal boundary, widely used in Mathematical Relativity). But
some problems about its consistency originated a long sequence of
redefinitions of this boundary (see \cite{Orl-M} for a critical
review). Recently, with the additional stimulus of finding a
general boundary for the AdS-CFT correspondence,  the notions of
c-boundary and c-completion have been widely developed
\cite{Ha1,H2,MR2,F,Orl-M}, and the recent  detailed study in
\cite{FHS0} justifies that a satisfactory definition is available
now. Moreover, several new ideas have been introduced for the
computation of the c-boundary (or several non-problematic elements
of it) in relevant cases \cite{AF,FH,FS,H}.

In this article we carry out a systematic study of the c-boundary
of a natural class of spacetimes, the (conformally) {\em
stationary} ones. However, our motivation is not only to compute
this boundary for a remarkable class of spacetimes, but also the
revision of other classical boundaries in Differential Geometry
---which turn out to be related with the causal one. More
precisely, in order to  describe the c-boundary of a standard
stationary spacetime $V=\R\times M$, the space of, say,
``diverging directions of curves on $M$'' appears in a natural
way. These directions must be computed with a Riemannian metric
$g$ in the particular case that the spacetime is {\em static}, and
with a Finsler metric $F$ of Randers type (and its reverse metric
$\F$) in the general stationary case. Such directions are computed
from Busemann-type functions constructed for arbitrary curves of
bounded velocity, and its topology is naturally defined from the
 {\em chronological topology} for any c-completion. Of
course, in the particular case that $(M,g)$ is a Hadamard
manifold, the boundary agrees with Eberlein and O'Neill's one.
But, in general, one obtains a new compactification of both, any
Riemannian manifold and any Randers manifold, the latter trivially
extensible to any Finslerian manifold.
Then, it is natural to compare this new {\em Busemann
compactification} of any Finslerian manifold with the classical
Gromov and Cauchy completions. As far as we know, neither the
Gromov compactification has been studied systematically for an
arbitrary Finsler manifold, nor the Cauchy completion has been
developed for a {\em generalized distance}, i.e. a non-necessarily
symmetric distance (as the ones induced by Finsler metrics). As
they present some non-evident properties (which will turn out to
be essential for the c-boundary)
our first aim will be to study both completions in some detail.
Then, we will   introduce the Busemann boundary of a Finsler
manifold (with no reference to the causal boundary or Lorentzian
Geometry) and compare it with the other two completions. Finally,
with all these elements at hand, a detailed picture of the c-boundary of a standard stationary spacetime will be achieved.

The contents of the paper are organized as follows. The first
subsection of Section \ref{s2} provides   some generalities on the
c-boundary construction for any strongly causal spacetime $V$. We
emphasize some subtleties on this boundary $\partial V$: (a) at a
point set level, the c-boundary is constructed by making pairs of
elements from the future and past causal (pre)boundaries
$\hat\partial V, \check\partial V$ so that a point, say
$P\in\hat\partial V$, may be paired with more than one point
$F,F'\in\check\partial V$ yielding more than one pair $(P,F),
(P,F')$ in $\partial V$, and (b) even when the points in
$\hat\partial V, \check\partial V$ yield only a single pair, the
topology in $\partial V$ cannot be regarded as a plain topological
quotient of these two boundaries. In fact, such involved cases
originated a long list of redefinitions of the c-boundary ---and,
remarkably, they will occur even is stationary spacetimes. So,  we
introduce the notion of {\em simple} c-boundary for the cases when
these subtle possibilities do not occur (Definition
\ref{simpletop}). One of the goals of the present article will be
to prove that, under very mild hypotheses, the c-boundary of a
stationary spacetime is simple. The other two subsections recall
some properties of Finsler manifolds and conformally stationary
spacetimes. The latter summarizes the approach in  \cite{CJS}, and
it allows to understand that the conformally invariant properties
of a stationary spacetime are codified in a Finsler metric. After
these preliminaries, the following topics are developed.
\smallskip


{\em Section \ref{s3}}. The Cauchy completion of a non-necessarily
symmetric distance is studied. Such distances   are well-known
since the old work by Zaustinsky \cite{Zau}, and our definitions
agree with this author. We distinguish between {\em generalized
distances} (which are close to the symmetric case, and appear in
Finsler manifolds) and {\em quasi-distances} (which will appear in
the extensions of the generalized distances to the completions,
and retain less properties), see Definition  \ref{gms}, Remark
\ref{ra4}. We study the Cauchy completion of a generalized
distance in some detail, because they contain some non-trivial
subtleties, as the following:

(a) There are three types of natural Cauchy completions, the
forward, the backward and the symmetrized ones (Definition
\ref{gg}, \ref{dbs}, Proposition \ref{prop11}). Of course, the
three completions coincide for any (symmetric) distance, but even
the symmetrized boundary of a Finsler manifold may behave in a
highly non-Riemannian way (Remark \ref{r11}, Example \ref{illus}).

(b) There are two natural definitions of (forward) Cauchy
sequence. The usual one (Definition. \ref{ff}, which agrees with
\cite{BCS}, for example) is a straightforward generalization of
the symmetric case. The alternative one (Definition  \ref{f}, Proposition
\ref{palt}) seems more natural for the Cauchy completion
(Proposition  \ref{paux}). These notions {\em are not equivalent}
(Example \ref{r2}), but they do yield {\em equivalent Cauchy
completions} (Theorem \ref{co}, Remark \ref{rco}), supporting  the
remainder of our approach.

(c) The generalized distance $d$ can be extended to the (say,
forward) Cauchy completion $\caum$ but, then, {\em the extension
$d_Q$ is only a quasi-distance} (Proposition \ref{ppp1'}). This
lies in the core of many subsequent difficulties; in fact, when
$d_Q$ is a generalized distance on $\caum$ then $\caum = \caumr
=\caums$ (Proposition \ref{genersym}). But for the quasi-distance
$d_Q$, the topology associated to the open forward balls may be
different to the one generated by the backward balls (and none  of
them is $T_1$). Moreover, only the latter topology ensures that
$\caum$ is truly a completion, i.e. any forward Cauchy sequence
will be convergent in $\caum$ with the topology generated by the
backward balls (Remark \ref{r324}, Convention \ref{conv}, Theorem
\ref{tcompltcompl}).

(d) It is even possible to extend further $d_Q$ (Proposition
\ref{Srelat}), yielding
 some relations between the forward and backward completions.
Such relations  will become quite natural under the viewpoint of
the c-boundary (in fact, the asymmetries between $M_C^+$ and
$M_C^-$ will reflect the asymmetries between the future and past
directions in the spacetime). They support the notion of a
generalized metric space with {\em evenly pairing boundaries}
(Definition \ref{dtotsep}; see Figure \ref{fig2}), which will be
useful for the computation of the
c-boundary 
of a stationary spacetime $V$.
\smallskip

\noindent We emphasize that all the explained subtleties already
appear in the simplest case of the generalized metric space
associated to a Finsler metric of Randers type (the relevant case
for the c-boundary of a stationary spacetime). So, the examples
which exhibit these subtleties are chosen  of this type.
\smallskip

{\em Sections \ref{s4}, \ref{s5}}. Here, Gromov and Busemann
completions and boundaries are extended to arbitrary Finsler
manifolds. For the convenience of the reader, some properties are
studied first  in Section \ref{s4} for the particular case of
Riemannian manifolds (or somewhat more general spaces as
reversible Finsler manifolds, see Convention
\ref{cginebra}).
 As far as we know, the Finslerian
versions in Section \ref{s5} had not been considered before. Even
though our study is independent of Lorentzian Geometry, our
definitions are guided by the detailed correspondence between
global Finsler properties and global conformally invariant
Lorentzian elements (causal objects) developed in \cite{CJS}. Some
connections between Finsler and Lorentzian elements had been also
introduced in \cite{Gibbons} under a  local viewpoint, and in
\cite{CJM} in the framework of variational methods.

For the Riemannian Gromov completion $M_G$, we distinguish between
the proper Gromov boundary $\partial_{\cal G}M$, which contains
``points at infinity'', and the Cauchy-Gromov one
$\partial_{CG}M$. The latter includes the Cauchy boundary
$\partial_CM$ and, roughly, it can be understood as a local
compactification of $\partial_CM$  obtained by means of both, a
coarser topology and the introduction of new boundary points which
constitute the residual Gromov boundary $\partial_{RG}M$. In
the Cauchy-Gromov completion $M_{CG}=M\cup\partial_{CG}M$, a
Heine-Borel property holds;
 so, the natural condition to ensure that the Cauchy completion
$M_C$ is topologically embedded in $M_{CG}$, is that $\partial_CM$
is locally compact (Corollary \ref{ccg}). All these properties are
proved in the Riemannian case in Subsection \ref{s41}, and then,
only the specific difficulties of the Finslerian case are analyzed
in Section \ref{secgro}. Essentially, these Finslerian
difficulties involve three facts. First, the appropriate choices
for notions such a Lipschitz function (Definition \ref{dflips},
Remark \ref{rrr1}) or oriented-equicontinuity (Definition
\ref{doreq}). 
Second, the existence of forward and backward completions $M_G^+,
M_G^-$, which yield more possibilities to deal with (for
example, the continuity of the inclusion of $M_C^+$ in $M_G^+$
holds only under an additional hypothesis, see Lemma \ref{lcf} and
Remark \ref{rcf}). And, third, some specific technicalities must
be analyzed now, as the  necessity of a version of Arzela's
theorem for oriented-equicontinuous curves (Theorem \ref{taf}),

About the Busemann completion, it is inspired on two elements. On
one hand (as a point set), in the commented Eberlein and O'Neill
completion of Hadamard manifolds and, on the other
(topologically), in the interpretation of the graphs of
Busemann-type functions (on Riemannian and Randers-type Finsler
manifolds) as points in the (future) causal completion of a
stationary spacetime. Some properties of the Busemann completion
for a Riemannian manifold had been already studied in \cite{H,
FH}. Here, we present a systematic study, which includes all the
Finslerian manifolds. For the convenience of the reader, we
explain first the Busemann completion $M_B$ as a point set in the
Riemannian case (Section \ref{s42}). Such a $M_B$ (Definition
\ref{db}, Remark \ref{remm1}) contains the points in Gromov's one
$M_G$ which can be obtained from Busemann-type functions.
Intuitively, this means that these points in $M_G$ can be
``reached by curves asymptotically like a ray'' (see Remarks
\ref{rgrape}, \ref{rgrape2}). {\cambios So, in spite of their very
different definitions, Gromov and Busemann completions agree in
most cases as point sets. However, the cases when they do not
coincide become also interesting, and the existence of such
differences will be equivalent to the existence of differences at
the topological level}. For reasons of space, the study of the
Busemann topology (or {\em chr. topology}) is postponed to the
general Finlerian case. But the reader interested only in
Riemannian Geometry will be able to
 reconstruct easily this part just by making some simplifications.

The Busemann completions $M_B^\pm$ in the
general Finslerian case (or the somewhat more general spaces under
Convention \ref{cginebra2}) are analyzed in Section
\ref{sectbusf}.
 Its natural
topology (the {\em chronological topology})  is adapted to the
structure of the Busemann functions. It becomes somewhat
technical, but  it will admit  natural interpretations from the
spacetime viewpoint. The original space $M$ can be embedded
naturally in $M_B^\pm$, and these completions yield two
compactifications of $M$. The  Cauchy completion $M_C^+$ is
included in $M_B^+$ so that $\partial_B^+M$ becomes the disjoint
union of $\partial_C^+ M$ (``incomplete finite directions'') and
the properly Busemann boundary $\partial_{\cal B}^+ M$
(``asymptotic directions''). The inclusion at the topological
level is also discussed, as in the Gromov case. The comparison
between Gromov's and Busemann's completions is carried out in
Section \ref{s53}.
  As a point set, $M_B^+$ is included in $M_G^+$ and this inclusion may be strict
(in particular, no point in $\partial^+_B M$  belongs to the
residual Gromov boundary $\partial^+_{RG}M$). The chronological
topology on $M_B^+$ is coarser than the induced from
$M^+_G$, but also {\cambios yields a compactification of $M$}. This is
done at the price of being non-Hausdorff, and has a natural
explanation: while Busemann boundary compactifies only among
``asymptotically ray-like directions'',  Gromov one introduces
additional points, not always reachable as the limit of some curve
in $M$. {\cambios Gromov completion has the nice Hausdorff property,
but Busemann completion does not introduce points which may appear
as spurious
 (Remark \ref{rgrape2}). Remarkably, both topologies agree when the
Busemann completion $M_B^+$ is Hausdorff, and this is equivalent
to the equality of $M_B^+$ and $M_G^+$ as point sets (Theorem
\ref{below}).}

In conclusion, our study gives a precise description of the three
boundaries for a Finslerian manifold, which becomes especially
simple when the conditions which ensure both, $M^+_B=M^+_G$ and
the inclusion of $M_C^+$ in $M_B^+$ is an embedding, hold.
The results in Sections \ref{s4}, \ref{s5} are summarized in the
following theorem:
\begin{theorem}\label{tresumen} Let $(M,F)$ be {\cambios any (connected)} Finsler manifold
with associated generalized distance $d$.

\begin{enumerate}
\item The (forward) Gromov completion $M_G^+$ of $(M,d)$ endowed
with the quotient of the pointwise convergence topology (see
Definition \ref{def1}) of $(M,d)$ is a compact metrizable
topological space.


\item The (forward) Busemann completion $M_B^+$ of $(M,d)$ endowed
with the chronological topology (see Definitions \ref{ddd0},
\ref{ddd1}) is a possibly non-Hausdorff $T_1$
{\cambios sequentially} compact topological space.


The (possibly non-continuous) inclusion $i: M^+_B\rightarrow
M^+_G$ satisfies that $i^{-1}: i(M^+_B) \rightarrow M^+_B$ is
always continuous (i.e., the chr. topology is coarser than the
Gromov one).

The following statements are equivalent:\begin{itemize}\item[a)]
{\cambios No sequence in $M_B^+$ converges to more than one point in
$\partial_{B}^{+}M$ with the chronological topology.} \item[b)]
$M_B^+$ is {\cambios Hausdorff} with the chronological topology.
\item[c)] $i: M^+_B\rightarrow M^+_G$ is continuous (and thus, an
embedding). \item[d)]  $M_B^+= M_G^+$ both, as  point sets and
topologically. \item[e)] {\cambios $M_B^+$ is equal to $M_G^+$ as point
sets}.

\end{itemize}

\item The (forward) Cauchy completion $M_C^+$ of $(M,d)$ is a set
endowed with a quasi-distance $d_Q$ (see Definition \ref{gms})
which generates a natural $T_0$-topology (see Convention
\ref{conv}). {\cambios Its relations with  $M_G^+$ and $M_B^+$
are} the following:
\begin{itemize}
\item[A)] The natural map $j^+:M_C^+\rightarrow M_G^+$ is
injective, and it is continuous if and only if the backward
topology {\cambios on $M_C^+$} is finer than the forward one (i.e.
condition {\em (a4')} in Lemma \ref{lcf} is satisfied). In
particular, $j^+$ is continuous if the quasi-distance $d_Q$ in
$M_C^+$ is a generalized distance (this property always happens in
the Riemannian case).

The Cauchy-Gromov completion $M_{CG}^+$ (Section \ref{secgro},
below Theorem \ref{tgromovf})
 includes naturally the Cauchy
one $M_C^+$ as a point set, and it satisfies the Heine-Borel
property (i.e., the closure {\cambios in $M_{CG}^{+}$ of any bounded
set in $M$ }is compact); in particular, it is  locally compact.

Moreover, if $d_Q$ is a generalized distance and $M_C^+$ is
locally compact: a) $j^+$ is an embedding, b) the Cauchy
completion agrees with the bounded part of the Gromov one both, as
a point set and topologically (there exists no {\em residual}
Cauchy-Gromov boundary, i.e. $\partial_{RG}^+ M :=\partial_{CG}^+
M\setminus \partial_C^+ M =\emptyset$).

\item[B)] Let $j^+_B: M_C^+ \rightarrow M_B^+$ obtained by
restricting $j^+$ (as $j^+(M_C^+)\subset M_B^+$) and consider the
natural topologies (chronological on $M_B^+$ and the one in
Convention \ref{conv} on $M_C^+$). Then, $j^+_B$ is continuous if
the backward topology generated by $d_Q$ on $M_C^+$ is finer than
the forward one. In particular, $j^+_B$ is always continuous in
the Riemannian case.

\end{itemize}

\end{enumerate}
At any case, $M$ is naturally an open dense subset of $M_G^+,
M_B^+$ and $M_C^+$, each one endowed with its natural topology.
So, the boundaries $\partial_G^+M$ and $\partial_B^+M$
are sequentially compact 
and these two boundaries coincide if and only if {\cambios no pair of
points in $\partial_B^+M$ is the limit of a single sequence in
$M$}.
\end{theorem}

{\em Section \ref{s6}}. By using previous machinery, the
c-boundary of a (standard) stationary spacetime $V$ is computed
now. The results are summarized in Theorem \ref{th} below (Figures
\ref{figdiagram}, \ref{figfinal} and  \ref{figfinal2} may be
clarifying). But, previously, the following remarks are in order:

(a) First, we have to compute  the future causal boundary
$\hat{\partial}V$, {\cambios composed by indecomposable past sets, IPs
(plus the past causal boundary $\check{\partial}V$, composed by
indecomposable future sets, IFs). This boundary}  will be endowed
with Harris' chronological topology defined in \cite{H2}. At least
when $M_B^+$ is Hausdorff, $\hat{\partial}V$ has a {\cambios plain
structure: it is just a cone on the forward Busemann completion
$\partial_B^+M$ (Proposition \ref{pepe}, Theorem
\ref{pconocau})}. Moreover, the points of the symmetrized
Cauchy boundary $\partial^s_C M$ yield timelike lines of the cone,
and the other points of $\partial_B^+M$ yield lightlike lines.
These results refine those in \cite{FH, H} for the static case,
and extend them to the stationary one.

(b) As we have already commented, some problems in the consistency
of the notion of  c-boundary $\partial V$ for a spacetime has
persisted along almost three decades. {\cambios They affected to: (i)
the {\em identification problem} between points of
$\hat{\partial}V$ and $\check{\partial}V$ (transformed in a
pairing prescription
by \cite{MR2}), (ii) the minimal character of the c-boundary \cite{F,
Orl-M} and (iii) its appropriate topology \cite{F}. The careful
study by the authors in \cite{FHS0} gives an answer to these
questions, and we follow this approach. Even though some of these
subtleties had appeared in somewhat pathological spacetimes, we
can find examples of all of them in standard stationary spacetimes
(see, for example, Remark \ref{r621}, with Examples \ref{ex1'},
\ref{ex2'}, about the questions (i) and (ii)). In particular,
the S-relation, which relates pairs IP-IF and defines the
c-completion as a point set (see (\ref{eSz}), (\ref{eSz2}) below)
was trivial in the static case, but it may be much more
complicated here (Remark \ref{r621}). The subtleties relative to
the topology become especially interesting, even in the static
case, as} the correct topological (and chronological) description
of $\partial V$ requires the notion of Busemann completion $M_B$
for a Riemannian manifold, plus the identification of $M_C$ with a
subset of $M_B$ (endowed with the topology induced by $M_B$). This
is studied in detail in Sections \ref{s6.3}--\ref{topcaubo}.
{\cambios The  refinement} of the definition of c-boundary left as an
open possibility in our study  \cite{FHS0} (the distinction
between causal boundary and {\em properly causal} boundary) does
not apply here, as the c-boundary  in the standard stationary case
is always properly causal, Remark \ref{remprocau}.

(c) We also emphasize the remarkable role of  the Cauchy
boundaries $\partial^\pm_C M$ in the computation of $\partial V$.
If these boundaries where empty, no point in $\partial V$ would be
timelike and $\partial V$ would be automatically globally
hyperbolic. So, no pairings between points of $\hat\partial V$ and
$\check
\partial V$ would occur (i.e. all the pairs in $\partial V$ would be type $(P, \emptyset)$, $(\emptyset, F)$
for some $P\in \hat\partial V$ or $F\in \check \partial V$).
Therefore, a mathematically elegant assumption of completeness for
the generalized distance $d$ (as it is assumed in the study of the
classical Gromov boundary) would spoil many beautiful
possibilities for spacetimes. Moreover, we emphasize that the
Busemann boundaries $\partial^\pm_BM$ split in two parts, the
proper {\cambios Busemann} boundary $\partial^\pm_{\mathcal{B}}M$ and
the Cauchy one $\partial^\pm_CM$ (the latter with a  topology
which comes from the Busemann viewpoint, not from the viewpoint of
metric spaces). Roughly, the former $\partial^+_{\mathcal{B}}M$
corresponds with asymptotic directions in $M$ and yields (most of)
the lightlike part of $\hat{\partial}V$. The latter
$\partial^+_CM$ corresponds with diverging directions at bounded
distance in $M$, it yields the timelike points of $\partial V$ and
some very particular lightlike ones.

\smallskip \noindent Summing up: 

\begin{theorem}\label{th} Let $(V,g)$ be a standard (conformally) stationary spacetime
(see (\ref{ggg}) and normalize $\Lambda\equiv 1$ with no loss of
generality). Let $F^\pm$ be its {\cambios associated Finsler
metrics (or {\em Fermat metrics}, see (\ref{h}))}, and denote
as $d^\pm$ the generalized distances associated to $F^\pm$. Then,
the c-boundary $\partial V$ has the following structure:
\begin{itemize}
\item[1.] {\em Point set}: \begin{itemize} \item[(1A)] The future
(resp. past) c-boundary $\hat{\partial} V$ (resp. $\check{\partial} V$) is naturally a point set cone with base $\partial_B^+ M$ (resp. $\partial_B^-
M$) and apex $i^+$ (resp. $i^-$). \item[(1B)] A pair $(P,F)\in
\partial V$ with $P\neq \emptyset$ satisfies that $P=P(b_c^+)$ for
some $c\in C^+(M)$ (see Proposition  \ref{caractIP}) and:
\begin{itemize}
\item[(a)] If $b_c^+\equiv \infty$ then $P=V, F=\emptyset$.
\item[(b)] If $b_{c}^+\in {\cal B}^+(M)(\equiv
\R\times\partial_{\cal B}^+M, \hbox{ see (\ref{ezzz})})$ then $F=\emptyset$. \item[(c)] If
$b_c^+\in B^+(M)\setminus {\cal B}^+(M)$, then $b_c^+=d_p^+$ with $p=(\Omega^+,x^+)\in \R\times
\partial_C^+M$ (see (\ref{equuu2})), $P=P(d^+_p)$ and $F\subset
F(d_{p}^-)$. In this case, there are two exclusive possibilities:
\begin{itemize}
\item[(c1)] either $F=\emptyset$, \item[(c2)] or $F=F(d_{p'}^-)$
with $p'=(\Omega^-,x^-)\in \R\times \partial_C^-M$ and satisfying
$\Omega^- -\Omega^+=d^+(x^+,x^-)$ (in this case, $p'$ is not necessarily
unique).
\end{itemize}
Moreover, if $x\in \partial_C^s M$, then $p'=p$, $\uparrow
P=F(d_{p}^{-})$ and $P$ is univocally S-related with $F=F(d_p^+)$.
\end{itemize}
A dual result holds for pairs $(P,F)$ with $F\neq \emptyset$. So,
the total c-boundary is the disjoint union of lines $L(P,F)$ (see
Definition \ref{deflin}).

When $\partial V$ is simple as a point set (see Definition
\ref{simpletop}; in particular this happens when $(M,d^+)$ has a
forward or backward {\ras}, Definition \ref{dtotsep}), then $\partial {V}$ it is the quotient set $\hat{\partial} V\cup_d \check{\partial} V/\sim_S $ of the partial boundaries $\hat{\partial} V, \check{\partial} V$  under the S-relation.

\end{itemize}


\item[2.] {\em Causality}: For each pair $(P,F)\in \partial V$,
the line $L(P,F)$ is timelike if $P=P(d_p^+)$ and $F=F(d_p^-)$ for
some $p\in \R\times \partial_C^s M$, horismotic if either $P$ or
$F$ are empty and locally horismotic otherwise (recall Definition
\ref{caurellin}).


\item[3.] {\em Topology}:
\begin{itemize}\item[(3A)] If $M_B^+ $ (resp.
$M_B^- $) is Hausdorff, the future (resp. past) causal boundary
has the structure of a (topological) cone (Definition
\ref{defcono}) with base $\partial_B^+M$ (resp. $\partial_B^-M$)
and apex $i^+$ (resp. $i^-$). \item[(3B)] If $M_{C}^{s}$ is
locally compact and $d_Q^+$ is a generalized distance, then
$\overline{V}$ is simple {\cambios (see Definition
\ref{simpletop}) and so, it coincides with the quotient
topological space $\hat{V}\cup_d \check{V}/\sim_S$ of the partial
completions $\hat{V}$ and $\check{V}$ under the S-relation.}
\end{itemize}
Summarizing, if $M_C^s $ is locally compact, $d^+_Q$ is a
generalized distance and $M_B^\pm $ is Hausdorff, $\partial V$
coincides with the quotient topological space
$(\hat{\partial}V\cup_d \check{\partial}V)/\sim_{S}$, where
$\hat{\partial}V$ and $\check{\partial}V$ have the structure of
cones with bases $\partial_B^+ M, \partial_B^- M$ and apexes
$i^+,i^-$, resp.
\end{itemize}
\end{theorem}

\begin{remark}
{\em {\cambios Previous theorem can be extended to the corresponding
causal completions in the obvious natural way (replacing the
boundary of each type by the corresponding completion), at all the
levels, i.e., as a pointset, causally and topologically.}}
\end{remark}

\section{Preliminaries}\label{s2}
\subsection{Spacetimes and c-boundaries}\label{fff'}

We will use typical background and terminology in Lorentzian
Geometry as in \cite{BEE, MS, O}. A spacetime will be a
time-oriented connected smooth Lorentzian manifold $(V,g)$ (the
time-orientation, and so, the choice of a future cone at each
tangent space, will be assumed implicitly) of signature $(-,+,
\dots , +)$. A tangent vector $v\in T_{p}V$, $p\in V$ is called
{\em timelike} (resp. {\em lightlike}; {\em causal}) if $g(v,v)<0$
(resp. $g(v,v)=0$, $v\neq 0$; $v$ is either timelike or
lightlike); null vectors include both, the lightlike ones and the
0 vector. A causal vector is called {\em future} or {\em
past-directed} if it belongs to the future or past cone.
Accordingly, a smooth curve $\gamma:I\rightarrow V$ ($I$ real
interval) is called timelike, lightlike, causal and future or
past-directed if so is $\dot{\gamma}(s)$ 
for all $s\in I$.

 Two events $p,q\in V$ are {\em
chronologically related} $p\ll q$ (resp. {\em causally related}
$p\leq q$) if there exists some future-directed timelike (resp.
either future-directed causal or constant) curve from $p$ to $q$.
If $p\neq q$ and $p\leq q$ but $p\not\ll q$, then $p$ is said {\em
horis\-mo\-ti\-cal\-ly related} to $q$. The {\em chronological
past} (resp. {\em future}) of $p$, $I^{-}(p)$ (resp. $I^{+}(p)$)
is defined as:
\[
I^{-}(p)=\{q\in V: q\ll p\}\qquad(\hbox{resp.}\; I^{+}(p)=\{q\in
V: p\ll q\}).
\]
In what follows, we will be especially interested in the
chronological past $I^-[\gamma]=\cup_{s\in I}I^-(\gamma(s))$
(resp. future $I^+[\gamma]=\cup_{s\in I}I^+(\gamma(s))$) of any
future-directed (resp. past-directed) timelike curve
$\gamma:I\rightarrow V$.

\begin{remark} \label{r1piecewise}
{\em Usually, in the literature a causal curve $\gamma$ with a
compact domain $I=[a,b]$ is allowed to be {\em piecewise smooth},
that is, it admits a partition $a=t_0<t_1<\dots <t_n=b$ such that:
(i) the restriction of $\gamma$ to each subinterval
$[t_{j},t_{j+1}]$ is smooth and causal, and (ii) the left and
right velocities of $\gamma$ at each break $t_j$ lie in the same
cone. This is technically  convenient, even though typically it
does not represent a true increment of generality (for example, in
order to define the relations $\ll, \leq$), because any such curve
admits a variation through smooth causal ones with the same
endpoints.

Here, we will use typically, say, future-directed timelike curves
defined on a half-open interval $I=[a,b)$, which cannot be
continuously extended to $b$. One can assume also that such a
curve $\rho$ is piecewise smooth in the sense that there exists a
strictly increasing sequence $\{t_j\}\nearrow b$ included in $I$
such that conditions (i) and (ii) above hold. In fact, on one
hand, all the techniques to be used here will work trivially in
this alternative case. On the other, for any such a piecewise
smooth curve $\rho$ there exists a smooth future-directed causal
curve $\tilde \rho$ arbitrarily close to $\rho$  such that
$I^-[\tilde\rho]=I^-[\rho]$ (see for example \cite[Sect.
3.3]{FS})\footnote{If $\rho$ were allowed to be continuously
extensible to $b$ then the extension may be non-differentiable at
$b$. However, it would be locally Lipschitzian  and $H^1$ (see
\cite[Appendix A]{CFS}), and any future-directed (smooth) causal
curve $\tilde{\rho}$ with endpoint $\rho(b)$ would satisfy
$I^-[\tilde \rho]=I^-(\rho(b))=I^-[\rho]$.}. }\end{remark}
The {\em c-completion} of spacetimes is constructed by adding {\em
ideal points} to the spacetime in such a way that any timelike
curve acquires some endpoint in the new space \cite{GKP}. To
formalize this construction, which will be conformally invariant
and applicable to {\em strongly causal} spacetimes, previously we
need to introduce some basic notions.

A non-empty subset $P\subset V$ is called a {\em past set} if it
coincides with its past; i.e. $P=I^{-}[P]:=\{p\in V: p\ll
q\;\hbox{for some}\; q\in P\}$. The {\em common past} of $S\subset
V$ is defined by $\downarrow S:=I^{-}[\{p\in V:\;\; p\ll
q\;\;\forall q\in S\}]$. In particular, the past and common past
sets must be open. A past set that cannot be written as the union
of two proper subsets, both of which are also past sets, is called {\em
indecomposable past} set, IP. An IP which does coincide with the
past of some point of the spacetime $P=I^{-}(p)$, $p\in V$ is
called {\em proper indecomposable past set}, PIP. Otherwise,
$P=I^{-}[\gamma]$ for some inextendible future-directed timelike
curve $\gamma$, and it is called {\em terminal indecomposable past
set}, TIP. The dual notions of {\em future set}, {\em common
future}, IF, TIF and PIF, are obtained just by interchanging the
roles of past and future in previous definitions.

To construct a {\em future} and a {\em past causal completion},
first we have to identify each {\em event} $p\in V$ with its PIP,
$I^{-}(p)$, and PIF, $I^{+}(p)$. This is possible in any {\em
distinguishing} spacetime, that is, a spacetime which satisfies
that two distinct events $p, q$ have distinct chronological
futures and pasts ($p\neq q \Rightarrow I^\pm (p) \neq I^\pm
(q)$). In order to obtain consistent topologies in the
completions, we will define the c-completion for a somewhat more
restrictive class of spacetimes, the {\em strongly causal ones}.
These are characterized by the fact that the PIPs and PIFs
constitute a sub-basis for the topology of the manifold $V$.

%

Now, every event $p\in V$ can be identified with its PIP,
$I^-(p)$. So, the {\em future causal boundary} $\hat{\partial}V$
of $V$ is defined as the set of all the TIPs in $V$, and  {\em the
future causal completion} $\hat{V}$ becomes the set of all the
IPs:
\[
V\equiv \hbox{PIPs},\qquad \hat{\partial}V\equiv
\hbox{TIPs},\qquad\hat{V}\equiv \hbox{IPs}.
\]
Analogously, each $p\in V$ can be identified with its PIF,
$I^+(p)$. The {\em past causal boundary} $\check{\partial}V$ of
$V$ is defined as the set of all the TIFs in $V$, and  {\em the
past causal completion} $\check{V}$ is the set of all the IFs:
\[
V\equiv \hbox{PIFs},\qquad \check{\partial}V\equiv
\hbox{TIFs},\qquad\check{V}\equiv \hbox{IFs}.
\]

For the (total) causal boundary, the so-called S-relation comes
into play \cite{S1,S2}.
Denote $\hat{V}_{\emptyset}=\hat{V}\cup \{\emptyset\}$ (resp.
$\check{V}_{\emptyset}=\check{V}\cup \{\emptyset\}$). The
S-relation $\sim_S$ is defined in $\hat{V}_{\emptyset}\times
\check{V}_{\emptyset}$ as follows. If $(P,F)\in \hat{V}\times
\check{V}$, \be \label{eSz}  P\sim_S F \Longleftrightarrow \left\{
\begin{array}{l}
F \quad
\hbox{is included and is a maximal IF in} \quad \uparrow P \\
P \quad \hbox{is included and is a maximal IP in} \quad \downarrow
F.
\end{array} \right.
\ee Recall that, as proved by Szabados \cite{S1}, $I^-(p) \sim_S
I^+(p)$ for all $p\in V$, and these are the unique S-relations
(according to our definition (\ref{eSz})) which involve proper
indecomposable sets. For $(P,F)\in \hat{V}_{\emptyset}\times
\check{V}_{\emptyset}$, with $(P,F)\neq (\emptyset,\emptyset)$, we
also put \be \label{eSz2} P\sim_S \emptyset, \quad \quad
(\hbox{resp.} \; \emptyset \sim_S F )\ee if $P$ (resp. $F$) is a
(non-empty, necessarily terminal) indecomposable set and it is not
S-related by (\ref{eSz}) to any other indecomposable set. Now, we
can introduce the notion of c-completion, according to
\cite{FHS0}:
\begin{definition}\label{d1} As a point set, the {\em c-completion} $\overline{V}$ of a strongly causal spacetime $V$ is formed by all
the pairs $(P,F)\in
\hat{V}_{\emptyset}\times\check{V}_{\emptyset}$ such that
$P\sim_{S} F$. The {\em c-boundary} $\partial V$ is defined as
$\partial V:=\overline{V}\setminus V$, where $V\equiv
\{(I^{-}(p),I^{+}(p)): p\in V\}$.
\end{definition}

%


The chronological relation of the spacetime is extended to the
completion in the following way. We say that $(P,F), (P',F')\in
\overline{V}$ are {\em chronologically related}, $(P,F)\ll
(P',F')$, if $F\cap P'\neq\emptyset$.

Finally, the topology of the spacetime is also extended to the
completion by means of the so-called {\em chronological topology }
({\em chr. topology}, for short).
The motivation for this topology comes from the following
observation \cite{H2}. The natural topology of $V$ can be
characterized by a limit operator $\hat{L}$
which determines the possible limits for any sequence $\sigma=\{x_n\}$ in $V$. More
precisely, we write $x\in\hat{L}(\sigma)$ if and only if
\begin{itemize}
\item[(a)] for all $y\ll x$, for all sufficiently large $n$, $y\ll
x_{n}$, and \item[(b)] for any $x'$ with $I^{-}(x)\subset
I^{-}(x')$, if for all $y\ll x'$, there is infinitely many $n$'s
with $y\ll x_{n}$, then $x=x'$.
\end{itemize}
 The operator $\hat L$ can be rewritten in terms of IP's and,
 then, it can be applied to any sequence in $\hat V$. Moreover, it can be
 used to define the natural {\em (future) chr. topology} on $\hat V$ by
 imposing that a subset $C\subset \hat V$ is closed if and
 only if  $\hat L(\sigma)\subset C$ for any sequence
$\sigma\subset C$ (recall that, now, the limit operator $\hat L$
determines the convergence of sequences in an indirect way
\cite[p. 562]{H2}). A dual construction can be carried out with a
limit operator $\check L$ which defines the {\em (past) chr.
topology} on $\check V$. For the chr. topology on the (total)
causal boundary $\overline{V}$, a similar construction is done by
using a single operator $L$ constructed from the previous two ones
$\hat L, \check L$.

Rigorously,   the   limit operator $L$ for sequences in
$\overline{V}$ is defined as follows: given a sequence
$\sigma=\{(P_{n},F_{n})\}\subset\overline{V}$ and $(P,F)\in
\overline{V}$,
$$(P,F)\in L(\sigma)\iff \hbox{ if } P\neq \emptyset, P\in \hat{L}(P_n) \hbox{ and if }F\neq \emptyset, F\in \check{L}(\{F_n\}).$$ where\footnote{By LI and LS we
mean  the usual inferior and superior limits of sets: i.e.
LI$(\{A_{n}\})\equiv
\liminf(A_{n}):=\cup_{n=1}^{\infty}\cap_{k=n}^{\infty}A_{k}$ and
LS$(\{A_{n}\})\equiv
\limsup(A_{n}):=\cap_{n=1}^{\infty}\cup_{k=n}^{\infty}A_{k}$.}
\[
\begin{array}{c}
\hat{L}(\{P_{n}\}):=\{P'\in\hat{V}: P'\subset {\rm
LI}(\{P_{n}\})\;\;\hbox{and}\;\; P'\;\;\hbox{is a maximal IP into}\;\; {\rm LS}(\{P_{n}\})\} \\
\check{L}(\{F_{n}\}):=\{F'\in\check{V}: F'\subset {\rm
LI}(\{F_{n}\})\;\;\hbox{and}\;\; F'\;\;\hbox{is a maximal IF into}\;\;
{\rm LS}(\{F_{n}\})\}.
\end{array}
\]
 Then, the operator $L$ defines the {\em closed sets} for the
chronological topology on $\overline{V}$ as those subsets
$C\subset \overline{V}$ such that $L(\sigma)\subset C$ for any
sequence $\sigma\subset C$. The chr. topologies on the partial
completions $\hat V, \check V$ are expressible analogously in
terms of $\hat L, \check L$, respectively.

\begin{remark}\label{propsimplepunt}{\em We remark the following assertions about the chronological
topology:

(1) Clearly, if $(P,F)\in L(\{(P_{n},F_{n})\})$ then
$\{(P_{n},F_{n})\}$ converges to $(P,F)$. When the  converse
happens, $L$ is called {\em of first order} (see the discussion
\cite[Section 3.6]{FHS0}).

(2) Given a pair $(P,F)\in \partial V$, any timelike curve
defining $P$ (or $F$) converges to $(P,F)$ with the chronological
topology (see \cite[Th. 3.27]{FHS0}).

}
\end{remark}
\smallskip

These definitions for the c-boundary construction involve some
particular 
subtleties, which are
essentially associated to the following two facts: an IP (or IF)
does not determine a unique pair in the c-boundary; and the
topology does not always agree with the S-relation, in the sense
that $$P\in \hat{L}(P_n)\not\Leftrightarrow F\in \check{L}(F_n).$$
This makes natural to distinguish the following special cases:
\begin{definition}\label{simpletop}
A spacetime $V$ has a {\cambios c-completion} $\overline{V}$ which is
{\em simple as a point set} if each TIP (resp. each TIF)
determines a unique pair in $\partial V$.

Moreover, the {\cambios c-completion}
is {\em simple} if it is simple as a point set  {\cambios and also {\em topologically}, i.e.
$(P,F)\in L(P_{n},F_{n})$ holds when either $P\in
\hat{L}(\{P_{n}\})$ or $F\in \check{L}(\{F_{n}\})$. }
\end{definition}

Observe that the condition of being simple as a point set is,
indeed, a condition only on $\partial V$ i.e., the notion of being
simple as a point set makes sense for $\partial V$ {\cambios
---this obeys to the original idea by Geroch, Kronheimer and Penrose
\cite{GKP} of regarding the causal boundary as a quotient of
$\hat\partial V$ and $\check \partial V$. Nevertheless, the
simplicity at the topological level involves also the spacetime
$V$. Moreover, from the definition of $L$, when $P\in
\hat{L}(\{P_{n}\})$ and $F\in \check{L}(\{F_{n}\})$ one has
$(P,F)\in L(P_{n},F_{n})$ always, and if, say,  $P$  is the empty
set, then only the second
 of the previous conditions is required (the first one would not make sense, as the
 empty set is never an
 element of $\hat L(P_n)$). So, the topological simplicity must be
 checked only when $P\neq \emptyset \neq F$.}

{\cambios When $\overline V$ is simple as a point set then one can
consider the disjoint union $\hat V \cup_d \check V$ and the
quotient $\hat V \cup_d \check V/\sim_S$ obtained by identifying
each TIP  $P$ (resp. IF $F$) with itself and, eventually, with the
unique $F$ such that $P\sim_S F$. Clearly, $\hat V \cup_d \check
V$ inherits a natural topology such that  $\hat V$ and $\check V$
are its connected parts and, then, $\hat V \cup_d \check V/\sim_S$
is also regarded as a quotient topological space.
}
\begin{proposition}\label{propsimplepunt1}
Let $(V,g)$ be a spacetime:

(1) If the chronological topology
associated to the c-completion $\overline{V}$ is Hausdorff, then
the {\cambios c-completion} $\overline{V}$ is simple as a point set.

(2) {\cambios If $\overline V$ is simple, then the } the
c-completion {\cambios is naturally homeomorphic to  $\hat V \cup_d
\check V/\sim_S$ through the mapping:}
\begin{equation}\label{epss}
 p:\overline{V}\rightarrow
\hat{V}\cup_d \check{V}/_{\sim S}, \quad \quad
p(P,F):=\left\{\begin{array}{cc} \hbox{[P]} & \hbox{ if $P\neq
\emptyset$}\\ \hbox{[F]} & \hbox{ if $F\neq \emptyset$}
\end{array}\right. .\end{equation}

\end{proposition}
{\cambios {\it Proof. } (1)  From Remark \ref{propsimplepunt} (2), any
timelike curve defining a TIP $P\neq \emptyset$ converges to any
pair of the form $(P,F)\in \overline{V}$ with the chronological
topology. Then, as $\overline{V}$ is Hausdorff, $F$ is unique.

(2) First, let us  characterize the closed sets in $\hat{V}\cup_d
\check{V}/_{\sim S}$. As the natural projection
$\pi:\hat{V}\cup_d \check{V}\rightarrow
\hat{V}\cup_d \check{V}/_{\sim S}$ is
continuous, and for any set $C\in \hat{V}\cup_d
\check{V}/_{\sim S}$
\[\pi^{-1}(C)=\hat{C}\cup_d \check{C}\subset
\hat{V}\cup_d \check{V}\] for some $\hat
C\subset \hat V, \check C\subset \check V$, we have that $C$ is
closed if and only if so are $\hat C$ and $\check C$.

Now, recall that $p$ in (\ref{epss}) is well defined and bijective
because of the simplicity as a point set of $\overline V$. To
check that $p$ is continuous, let $C$ be a closed subset of
$\overline{V}\rightarrow \hat{V}\cup_d \check{V}/_{\sim S}$ and
assume by contradiction that $p^{-1}(C)$ is not closed in
$\overline{V}$. Then, there exist a sequence $\{(P_n,F_n)\}\subset
p^{-1}(C)$ and some $(P,F)\notin p^{-1}(C)$ such that $(P,F)\in
L(P_n,F_n)$. From the definition of the limit operator $L$, if
$P\neq \emptyset$  (otherwise, reason with $F\neq \emptyset$) then
$P\in \hat{L}(P_n)$. Nevertheless, this is a contradiction with
the closedness of $C$ (and then $\hat C$), because  by hypothesis
$\{P_n\}\subset \hat{C}$ but $P\notin \hat{C}$.

 For the continuity
of $p^{-1}$, let $C$ be now a closed set in $\overline{V}$ and
assume by contradiction that $p(C)$ is not  closed in
$\hat{V}\cup_d \check{V}/_{\sim_S}$. That
is, $\pi^{-1}(p(C))=\hat{p(C)}\cup_d \check{p(C)}$ is not closed
in $\hat{V}\cup_d \check{V}$ and, with no loss of generality,
assume that $\hat{p(C)}$ is not closed in $\hat{V}$.
Then, there exist a sequence $\{P_{n}\}\in \hat{p(C)}$ and some
$P\notin \hat{p(C)}$ such that $P\in \hat{L}(P_n)$. As
$\overline{V}$ is (topologically) simple, $(P,F)\in L(P_n,F_n)$,
with $(P_n,F_n)\in C$ and $(P,F)\notin C$, in contradiction with
the closed character of $C$. \cvd }

\subsection{Finsler manifolds}
In this section, some basic elements of Finsler manifolds are
outlined (for a general background, see for example \cite{BCS}).

\begin{definition} A {\em Finsler metric} on a smooth manifold $M$ is a non-negative
function $F$ defined {\cambios on the tangent bundle} $TM$ which
is continuous on $TM$, $C^{\infty}$ on $TM\setminus \{0\}$, it
vanishes only on the zero section, it is fiberwise positively
homogeneous of degree one (i.e. $F(x,\lambda v)=\lambda F(x,v)$
$\forall x\in M,\; v\in T_{x}M$ and $\lambda>0$), and $F^2$ is
fiberwise strongly convex (i.e. the matrix
$[\partial^{2}F^{2}/\partial v^{i}\partial v^{j}(x,v)]$ is
positive definite for any $(x,v)\in TM\setminus \{0\}$, {\cambios
where the $v^i$'s denote linear coordinates in $T_xM$}).

A {\em Finsler manifold} is a pair $(M,F)$ composed by a smooth
manifold $M$ and a Finsler metric $F$.
\end{definition}
 Notice that, in general, a
Finsler metric $F$ is not reversible, i.e. it may happen that
$F(x,v)\neq F(x,-v)$
for some $(x,v)$. 
The {\em reversed Finsler metric} $\F$ of $F$ is
then
$$\F(x,v):=F(x,-v)\qquad\hbox{for all $(x,v)$.}$$
In what follows, the classical notation above will be simplified
by putting $F(v)$ ($v\in TM$) instead of $F(x,v)$ ($x\in M, v\in
T_xM$). The {\em length} of a piecewise smooth curve
$c:[s_{0},s_{1}]\rightarrow M$ with respect to a Finsler metric
$F$ on $M$ is given by:
\[
{\rm length}_{F}(c):=\int_{s_{0}}^{s_{1}}F(\dot{c}(s))ds.
\]
Because of non-reversibility, this  length is not invariant under
reparametrizations which change the orientation. However, one
still has:
\begin{proposition}\label{ppp1} For any Finsler manifold $(M,F)$ the map $d_{F}:M\times
M\rightarrow \R$
$$d_{F}(x,y):=\inf_{c\in
C(x,y)}{\rm length}_{F}(c)\geq 0,$$ where $C(x,y)$ is the set of
piecewise smooth curves starting  at $x$ and ending at $y$,
satisfies all the axioms of a distance on $M$ but symmetry (i.e.,
$d_F(x,y)\neq d_F(y,x)$, in general).
\end{proposition}
The conclusion can be restated by saying that $(M,d_F)$ is a
generalized metric space according to Defn. \ref{gms} below (see
Remark \ref{rguev}).

A {\em Randers metric} on $M$ is a Finsler metric $F$ which can be
written as $F(v)= \sqrt{g_0(v,v)+\omega(v)^2}+\omega(v)$ for some
Riemannian metric $g_0$ and 1-form $\omega$ on $M$. Randers
metrics constitute a typical class of  non-reversible Finsler
metrics, and will appear naturally in standard stationary
spacetimes.
\subsection{Conformally stationary spacetimes}

A {\em conformally stationary} (resp. {\em stationary}) spacetime
  is  a spacetime $(V,g)$ which admits a timelike conformal-Killing
(resp. Killing) vector field $K$; in this case, the vector field
$K$ is also called {\em conformally stationary (resp.
stationary)}. $(V,g)$ is {\em standard conformally stationary} if
it can be written as: \be\label{ggg} V=\R\times
M\qquad\hbox{and}\qquad g=\Lambda(-dt^2+\pi^*\omega\otimes
dt+dt\otimes \pi^*\omega +\pi^*g_{0}), \ee where $(M,g_{0})$ is a
Riemannian manifold, $\omega$ a 1-form on $M$,  $\Lambda$ a
positive function on $V$ and $t: \R\times M \rightarrow \R$,
$\pi:\R\times M\rightarrow M$ the natural
projections.
Locally, any conformally stationary spacetime looks like a
standard one with $K=\partial_t$. {\cambios If $\Lambda$ in
(\ref{ggg}) is independent of the $t\in \R$ } coordinate then
$(V,g)$ is called {\em standard stationary}.

Along this work, we will assume two harmless simplifications.
First, the conformally stationary vector field $K$ will be always
unitary ($g(K,K)\equiv -1$) and, thus, stationary. The reason is
that we will be interested only in conformally invariant elements
(as so is the c-boundary) and, thus, we can change the original
metric $g$ by $g^*=-g/g(K,K)$ with no loss of generality.

The second one is to consider always the standard stationary case
i.e., globally the expression (\ref{ggg}) {\cambios (with $\Lambda
\equiv 1$)}. This is a natural simplification because a recent
result in \cite{JS} shows that a stationary spacetime $V$ is
standard if and only if: (i) it is {\em distinguishing}, and (ii)
admits a {\em complete} stationary vector field $K$. As we have
seen, condition (i) is necessary to deal with the c-boundary (in
fact, it is less restrictive than strong causality). Condition
(ii) is not fulfilled in general: typically, one can put some
holes in a standard stationary spacetime and both, the standard
character and the completeness of
$K$ are dropped. 
However, in such an example, the spacetime would admit an open
conformal embedding in the original standard stationary spacetime,
(and, for example, the results in \cite{FHS0} which relate the
causal and conformal boundaries may be applicable). So, the
restriction to the standard case allows us to focus on the
specific problems of the c-boundary for a  stationary spacetime.

 As we will see along this
paper, the structure of the c-boundary for these spacetimes is
closely related to properties of two Finsler metrics $F^+, F^-$ of
Randers type on the spatial part $M$. These metrics (or {\em
Fermat metrics}) are defined in terms of $g_{0}$ and $\omega$ (for
$\Lambda\equiv 1$); concretely,
\begin{equation}\label{h}
F^\pm(v)= \sqrt{g_0(v,v)+\omega(v)^2 } \pm \omega(v)\qquad \forall v\in TM .
\end{equation}
We will put simply $F=F^+$, and then, $\F=F^{-}$. It is worth
pointing out that the same standard stationary spacetime can be
decomposed as standard stationary in more than one way (with
different $g_0$ and $\omega$). Then, all the Fermat metrics
associated to some decomposition will share those properties which
are associated to the geometry of the spacetime $V$. These
properties include the possible global hyperbolicity or causal
simplicity  of the spacetime (see \cite{CJS}).

\section{Cauchy completion of a generalized metric space}
\label{s3}

\subsection{Basic definitions}
We start by outlining some basic elements of generalized metric
spaces (for a general background see \cite{Zau}).
\begin{definition}\label{gms}
A {\em generalized metric space} is a pair $(M,d)$, where $M$ is a
set of points and $d:M\times M\rightarrow \R$ is a {\em
generalized distance}
(or {\rm generalized metric}), i.e. $d$ satisfies the following
axioms:
\begin{itemize}
\item[(a1)] $d(x,y)\geq 0$ for all $x,y\in M$. \item[(a2)]
$d(x,y)=d(y,x)=0$ if and only if $x=y$. \item[(a3)] $d(x,z)\leq
d(x,y)+d(y,z)$ for all $x,y,z\in M$. \item[(a4)] Given a sequence
$\{x_n\}\subset M$ and $x\in M$, then $\lim_{n\rightarrow
\infty}d(x_n,x)=0$ if and only if $\lim_{n\rightarrow
\infty}d(x,x_n)=0$.
\end{itemize}
A {\em quasi-distance} on $M$ is a map $d_Q: M\times M\rightarrow
\R \cup \{\infty\}$ which satisfies the axioms (a1), (a2) and (a3)
above. 
\end{definition}
The {\em reversed generalized distance} $\rd$ of a generalized
metric $d$ is defined as $\rd (x,y)=d(y,x)$, and the {\em
symmetrized distance} $\ds$ as
\begin{equation}\label{symdist}
\ds(x,y)=\frac{1}{2}(d(x,y)+\rd (x,y)).
\end{equation}
It is direct to check that $\ds$ is a true (symmetric) distance.

One can define for $d$ the notions of {\em forward and backward
balls} of center $x\in M$ and radius $r\geq 0$; say, for open
balls,
$$B^+(x,r)=\{y\in M : d(x,y)<r\}, \quad B^-(x,r)=\{y\in M :
d(y,x)<r\},\quad \hbox{resp.}$$ 
These balls coincide with the corresponding backward and forward
balls for $\rd$, resp. In what follows, we will use ``forward
elements'',  being the backward ones analogous (and corresponding
to $\rd$). The balls for the symmetrized distance will be denoted
 $B^s(x,r)$.

{\cambios From the triangle inequality {\em (a3)}, if $y\in
B^{+}(x,r)$ one has $B^{+}(y,r')\subset B^{+}(x,r)$ whenever
$0<r'<r-d(x,y)$, and the forward open balls constitute a basis for
a (first countable) topology. Moreover, with this topology a
sequence $\{x_n\}\subset M$ converges to some $x\in M$ if and only
if $d(x,x_n)\rightarrow 0$. By {\em (a4)} this convergence is
equivalent to the convergence with the topology generated by the
open backward balls (which will be equal to the former forward
topology).}

\begin{remark} \label{ra4}{\em  Previous definitions are extended naturally
to any quasi-distance $d_Q$. Recall that for such a $d_Q$ the
hypothesis {\it (a4)} in Defn. \ref{gms} holds if and only if the
topology generated by the forward balls coincide with the topology
generated by the backward ones --and, then, also by the open
symmetrized balls. {\cambios So, any generalized metric or
quasi-distance satisfying {\it (a4)} will be regarded as a
topological space with the topology generated by the open forward
and backward  balls, and the limits in {\it (a4)} characterize the
convergence of $\{x_n\}$ to $x$ in the unique natural sense.}

Moreover, this condition {\it (a4)} implies the following natural
property: If $d_Q(x,y)=0$ then $d_Q(y,x)=0$, and thus, by {\it
(a2)}, $x=y$ {\cambios (in particular, the topology is Hausdorff)}. So,
the strengthening of {\it (a2)}, $d(x,y)=0$ iff $x=y$, holds in
any generalized metric space $(M,d)$, but it is easy to check that
this property does not imply {\it (a4)} (see Example \ref{illus}).
}\end{remark}
\begin{remark}\label{rguev}
{\em We will be interested in the study of Finsler manifolds
$(M,F)$, which yield the generalized metric spaces $(M,d_{F})$,
$(M,d_{F^{\rm rev}})$: in fact, these spaces satisfy condition
{\it (a4)} in Defn. \ref{gms} because the topologies associated to
$d_F$ and $d_F^{\rm rev}=d_{F^{\rm rev}}$ coincides with the
topology of the manifold (see for example \cite[Sect. 6.2
C]{BCS}).}
\end{remark}

\subsection{Forward and backward Cauchy boundaries as  point sets}

\subsubsection{Cauchy sequences and completions.}
\begin{definition}\label{ff} A sequence $\sigma=\{x_{n}\}$ in a generalized metric space $(M,d)$ is a
{\em (forward) Cauchy sequence} if for all $\epsilon>0$ there
exists $n_{0}$ such that $d(x_{n},x_{m})<\epsilon$ whenever
$n_{0}\leq n\leq m$. The space  of all the Cauchy sequences will be
denoted ${\rm Cau}^{+}(M)$ or, simply, $\cauchyd$. A {\em
backward Cauchy sequence} is a (forward) Cauchy sequence for
$\rd$. The space of all the backward Cauchy sequences will be denoted
by $\cauchydr$.
\end{definition}
Obviously, these definitions agree with the usual ones for
(symmetric) distances, and are also extendible to quasi-distances.
{\cambios It is easy to check that any convergent sequence (with
the natural topology, recall Remark \ref{ra4}) is Cauchy for a
generalized metric (see Proposition \ref{propsim}(3) below) but
this property fails for a quasi-distance (Remark
\ref{r324}(4)). }

Note the following straightforward assertion on double limits,
which will be claimed several times:
\begin{lemma} \label{claim} Two sequences $\sigma=\{x_{n}\}, \sigma'=\{x'_{n}\}\subset M$
satisfy $\lim_{n}(\lim_{m}d(x_{n},x'_{m}))$ $= \ell$, $\ell\in \R$
(resp. $\ell=\infty$) if and only if for all $\epsilon>0$ there
exists $n_0>0$ such that for all $n\geq n_{0}$ there is some
$m_{0}(n)$ satisfying: if $m\geq m_{0}(n)$ then
$|d(x_{n},x'_{m})-\ell | <\epsilon$ (resp.
$d(x_{n},x'_{m})>\epsilon$).
\end{lemma}
As a first technical result we have:
\begin{lemma}\label{lll1} For any $x\in M$ and $\{x_{n}\}$, $\{y_{n}\}\in
\cauchyd$, the sequence $\{d(x,y_{m})\}$  converges in $\R$, the
sequence $\{d(y_{m},x)\}$ converges in $[0,\infty]$, and the
double limit
\begin{equation}\label{defdist}
\ell= \lim_{n}(\lim_{m}d(x_{n},y_{m}))\end{equation} exists in
$[0,\infty]$. Moreover, there exist subsequences $\{x_{n_k}\}_k,
\{y_{m_k}\}_k$ such that $\ell = \lim_k d(x_{n_k}, y_{m_k})$.
\end{lemma}
{\it Proof.} For the first assertion (the second one is
analogous), assume first that the sequence does not converge in
$[0,\infty]$. Then, there exist two
 subsequences converging to different limits. Arranging
 appropriate subsequences $\{y_{m_{k}}\}$, $\{y_{m'_{k}}\}$
such that $m_{k}<m'_{k}$ and
$$0<\epsilon_0  \leq d(x,y_{m'_{k}})-d(x,y_{m_{k}})\leq
d(y_{m_{k}},y_{m'_{k}}) \qquad \forall k\in \N,
$$ 
a contradiction with the fact that $\{y_{m}\}$ is a Cauchy
sequence is obtained. For the finiteness of the limit, observe
also that $\{d(x,y_{m})\}_m$ is clearly bounded, and thus, it must
converge in $\R$.

For the limit (\ref{defdist}), denote $L_n=
\lim_{m}d(x_{n},y_{m})$ and assume by contradiction that there
exist subsequences $L_{n_k}, L_{n'_k}$ such that
$$0< 2 \epsilon_1 < L_{n'_k} - L_{n_k}.$$
Without loss of generality, we can also assume $n'_{k}<n_{k}$.
Then, there exist two sequences $\{y_{m_{k}}\}$ and
$\{y_{m'_{k}}\}$ such that $m_{k}<m'_{k}$ and
\begin{equation}\label{eqeq2} 0<\epsilon_{1}\leq
d(x_{n'_{k}},y_{m'_{k}})-d(x_{n_{k}},y_{m_{k}}) \leq
d(x_{n'_{k}},x_{n_{k}})+d(y_{m_{k}},y_{m'_{k}}) \qquad 
\forall k\in \N 
\end{equation}
in contradiction to the fact that $\{x_{n}\}$, $\{y_{m}\}$ are
both Cauchy sequences.

Finally, the last assertion follows easily from Lemma \ref{claim}.
\cvd

\begin{remark}\label{rll1} {\em As a technical question to be used later,
notice that, if $\{x_{n}\}\in \cauchyd$ and $\{y_{n}\}\in
\cauchydr$, we still have that $\lim_{n}(\lim_{m}d(x_{n},y_{m}))$
exists in $[0,\infty]$. To check it, just interchange the roles of
$y_{m_k}$ and $y_{m'_k}$ in formula (\ref{eqeq2}). Of course, by
applying Lemma \ref{lll1} to $\rd$, we also have that
$\{d(x,y_{m})\} (=\{\rd(y_{m},x)\})$ converges in $[0,\infty]$ and
$\{d(y_{m},x)\} (=\{\rd(x,y_{m})\})$ converges in $\R$.}
\end{remark}

\begin{definition}\label{rcs} Two Cauchy sequences $\sigma=\{x_{n}\}, \sigma'=\{x'_{n}\}
\in\cauchyd$ are {\em related}, $\sigma\sim\sigma'$, if
\be\label{sim}\lim_{n}(\lim_{m}d(x_{n},x'_{m}))=\lim_{n}(\lim_{m}d(x'_{n},x_{m}))=0.\ee
\end{definition}

\begin{proposition} \label{propsim} (1) The binary relation $\sim$ is a relation of  equivalence on $\cauchyd$.

(2) Let $x\in M$, $\sigma=\{x_{n}\}, \sigma'=\{x'_{n}\}
\in\cauchyd$. If $\sigma\sim\sigma'$ then $\lim_nd(x_n,x)=\lim_nd(x'_n,x) (\in [0,\infty])$.

(3) If a representative of a class of equivalence converges, then
any other representative of that class converges to the same
limit. If two sequences $\{x_{n}\}, \{x'_{n}\}$ converge to the
same point $x$, then both sequences are Cauchy and $\{x_{n}\}\sim
\{x'_{n}\}$.
\end{proposition}
{\it Proof.} (1) Symmetry and reflexivity are direct from the
definitions. 
To prove
transitivity,
consider Cauchy sequences  $\{x_{n}\} \sim \{x'_{n}\} \sim
\{x''_{n}\}$. Let $\epsilon>0$, by applying Lemma \ref{claim} to
the first two
sequences on the one hand, and to the last two sequences on the other hand, one obtains the corresponding natural numbers 
$n_0=:n_{12}, m_0(n)=:m_{12}(n)$ and $n_{23}$, $m_{23}(n)$, resp.
To check that $\{x_{n}\}  \sim \{x''_{n}\}$, define $l(n):=$
max$\{m_{12}(n),n_{23}\}$ for each $n\geq n_{12}$. Then, use the
claim with $n_{13}:=n_{12}$, $m_{13}(n):=m_{23}(l(n))$, and take
into account that, for $n\geq n_{13}, m\geq m_{13}(n)$, it is
\[
d(x_{n},x''_{m})\leq d(x_{n},x'_{l(n)})+d(x'_{l(n)},x''_{m})
<2 \epsilon.
\]
So, $\lim_{n}(\lim_{m}d(x_{n},x''_{m}))=0$, and the limit
$\lim_{n}(\lim_{m}d(x''_{n},x_{m}))=0$ is deduced analogously.

(2) From Lemma \ref{lll1} both limits exist in $[0,\infty]$, and
$$\lim_n d(x_n,x)\leq \lim_n( \lim_m (d(x_n,x'_m) + d(x'_m,x)))=\lim_m
d(x'_m,x).$$ So,  $\lim_n d(x_n,x)\leq \lim_m d(x'_m,x)$
 and the converse is
analogous.

 (3) The first assertion is obvious from item (2) {\cambios as the convergence is characterized by
 the limits with the distance (recall Remark \ref{ra4}, and its
 paragraph above)}

{\cambios For the last assertion, recall that $\lim_n d(x_n,x)=0=\lim_n d(x'_n,x)=\lim_n d(x,x'_n)$ (where the last equality follows from condition (a4) in Definition \ref{gms}) and $d(x_n,x'_m)\leq d(x_n,x)+d(x,x'_m)$. Then, for any $\epsilon>0$ there exists $n_0$ such that, for all $m,n\geq n_0$, $d(x_n,x'_m)<\epsilon$. This last property ensures that $\{x_n\}\sim \{x'_n\}$. For the Cauchy character of (each) sequence, put $x_n=x'_n$ in previous proof (obviously, the sequence is then both, forward and backward Cauchy).}
\cvd


\begin{remark} \label{rquasi} {\em The first two assertions in Proposition \ref{propsim} can be extended
directly to  Cauchy sequences of any set with a quasi-distance
(for the third one, a discussion on the topologies associated to
the quasi-distance is required, see Remark \ref{r324} below). This
would allow to define also the Cauchy completion of any set
endowed with a quasi-distance, as we will do next for generalized
metric spaces (Defn. \ref{gg}). Even though we will not be
interested in this possibility, the extension of Proposition
\ref{propsim}{\cambios (2)}
to quasi-distances will
be used below in order to ensure that the Cauchy completion of a
generalized metric space
becomes truly complete (see Theorem \ref{tcompltcompl}
(iii)).}\end{remark}

Now, we are in conditions to define the (generalized) Cauchy
completions.
\begin{definition}\label{gg} The {\em (forward) Cauchy completion} $\caum$ of a generalized metric space $(M,d)$ is the quotient space $\caum
:=\cauchyd/\sim$. Then, $M$ is regarded naturally as the subset of
$\caum$ consisting of the classes of equivalence of all the
sequences converging in $M$ ($\subset \caum$). The {\em (forward)
Cauchy boundary} of $(M,d)$ is $\bcaum :=\caum \setminus M$.

The {\em backward Cauchy completion} $\caumr$ and {\em boundary}
$\bcaumr$ are defined analogously by using the reversed generalized distance
$\rd$.
\end{definition}

\brema \label{r1}{\rm (a) In the Finslerian case, $d_F$ and $d_{\F}$
are generalized distances associated to a {\em length space} (in
the sense of \cite{Gr, Gr81}); this means that such distances are
computable  from the infimum of the (previously defined) lengths of
suitable curves connecting two points. In particular, this implies
that an appropriate version of Hopf-Rinow Theorem holds (see
\cite[Appendix]{CJS} for further details).

(b) Thus, the corresponding Cauchy boundaries can be constructed
by using curves of finite length instead of Cauchy sequences.
Namely, given a piecewise smooth curve  (as in the case of
spacetimes discussed in Remark \ref{r1piecewise}, one can consider
alternatively smooth or piecewise smoooth curves)
$\gamma:[0,b)\rightarrow M, b<\infty$, parametrized with the
length, any sequence $\{s_n\}\nearrow b$ yields a Cauchy sequence
$\{\gamma(s_n)\}$ --and all the so-obtained Cauchy sequences are
equivalent for $\sim$. Conversely, given a Cauchy sequence, a
curve $\gamma$ as above is obtained by connecting appropriate
terms of the sequence by means of arc-parametrized segments with
lengths approaching fast to the distances between their endpoints.
This viewpoint will be used later, when the Cauchy completion is
regarded as a part of the Busemann one. }\end{remark}

\subsubsection{Alternative Cauchy sequences yield the same
completions}

The equivalence relation defined by the double limits in
(\ref{sim}) suggests the following alternative way to introduce
the notion of Cauchy sequence.
\begin{definition}\label{f}  A sequence $\sigma=\{x_{n}\}$ in a generalized metric space $(M,d)$ is a {\em (forward) alternative Cauchy sequence} if
$\lim_{n}(\lim_{m}d(x_{n},x_{m}))=0$. The space of all the
alternative Cauchy sequences will be denoted ${\rm Cau}^{+}_{{\rm
alt}}(M)$ or, simply, ${\rm Cau}_{{\rm alt}}(M)$. A {\em backward
alternative Cauchy} sequence is a (forward) alternative Cauchy
sequence for $\rd$. The space of all the backward alternative
Cauchy sequences will be denoted ${\rm Cau}^{-}_{{\rm alt}}(M)$.
\end{definition}
Let us summarize the elementary properties of alternative Cauchy
sequences.
\begin{proposition}\label{palt}
Let $\sigma=\{x_n\}$ be a sequence in a generalized metric space
$(M,d)$. Then:
\begin{itemize}
\item[(1)] $\sigma$ is an alternative Cauchy sequence if and only
if for all $\epsilon>0$ there exists $n_0$ such that for all
$n\geq n_0$ there is $m_0(n)$ satisfying: if $m\geq m_0(n)$ then
$d(x_n,x_m)<\epsilon$. \item[(2)] Any Cauchy sequence is also an
alternative Cauchy sequence.

\item[(3)] If $d$ is symmetric then the alternative Cauchy sequences
coincide with the Cauchy ones.
\end{itemize}
\end{proposition}
{\it Proof.} (1) Apply Lemma \ref{claim} with $\sigma=\sigma'$.

(2) From the definition of alternative Cauchy sequence, the
property in (1) holds by putting $m_{0}(n)=n$.

(3) It reduces to prove the converse of (2). So, let
$\sigma=\{x_n\}$ be an alternative Cauchy sequence, and
$\epsilon>0$. Take $n_0$ and $m_0(n)$ as in (1), and define
$m_{12}=$max$\{m_0(n_1),m_0(n_2)\}$ for each $n_1,n_2\geq n_0$.
Then:
$$d(x_{n_1},x_{n_2})\leq d(x_{n_1},x_{m_{12}})+d(x_{m_{12}},x_{n_2})=d(x_{n_1},x_{m_{12}})+d(x_{n_2},x_{m_{12}})<
2\epsilon ,$$ 
as required. \cvd
\begin{example}\label{r2} {\em Let us show that, in general, there are alternative Cauchy sequences
which are not Cauchy sequences. Consider the Finsler manifold
$(M,F)$ given by
\[
M=(0,\infty),\;\;\;F=\sqrt{g_{0}+\omega^{2}}+\omega,\qquad\hbox{where}\;\;
g_{0}=\frac{dx^{2}}{1+x^{2}},\;\;\; \omega=-dx,
\]
and the sequence $\sigma=\{x_n\}\subset M$, $x_{n}=n+(-1)^{n+1}$.
Observe that $d_{F}(x_{2n-1},x_{2n})>2$, and so, $\sigma$ is not a
Cauchy sequence according to Defn. \ref{ff}. However, $\sigma$ is
an alternative Cauchy sequence according to Defn. \ref{f}, since
$$\lim_{n}(\lim_{m}d_{F}(x_{n},x_{m}))=\lim_{n}(\lim_{m}{\rm
length}_{F}(\alpha\mid_{(x_{n},x_{m})}))=0,\quad\hbox{where}\;\;
\alpha(t)=t.$$}
\end{example}
The relevance of the notion of alternative Cauchy sequence for the
Cauchy completion is stressed in the first assertion of the
following proposition.
\begin{proposition}\label{paux}
(1) Let $\sigma$, $\sigma'$ be two arbitrary sequences in a
generalized metric space $(M,d)$. If the double limits (\ref{sim})
hold then $\sigma$ and  $\sigma'$ are alternative Cauchy
sequences.

(2) The extension of the binary relation $\sim$ (Defn. \ref{rcs})
to the set ${\rm Cau}_{{\rm alt}}(M)$ of all the alternative
Cauchy sequences is also a relation of equivalence.

(3) If a representative of a class of equivalence in ${\rm Cau}_{{\rm alt}}(M)/\sim$ converges, then
all the representatives of that class converge to the same limit.
\end{proposition}
{\it Proof.} We will prove that $\sigma'$ is an alternative Cauchy
sequence by using Prop. \ref{palt} (1), being the proof of the
remainder assertions similar to the case of Cauchy sequences.
Given $\epsilon>0$, use the values of $n_1$ and $m_1(n)$ (for
$n\geq n_1$) in Lemma \ref{claim} so that $d(x'_n,x_m)<\epsilon$
whenever $m\geq m_1(n)$. Analogously, let $n_2>0$ and $m_2(n)$
(for $n\geq n_2$) so that $d(x_n,x'_m)<\epsilon$ whenever $m\geq
m_2(n)$. Put $n_{0}:=n_{1}$ and, for $n\geq n_{0}$, $k(n):={\rm
max}\{m_{1}(n), n_{2}\}$ and $m_{0}(n):=m_{2}(k(n))$. Then, if
$n\geq n_{0}$ and $m\geq m_{0}(n)$, it is
$$d(x'_{n},x'_{m})\leq
d(x'_{n},x_{k(n)})+d(x_{k(n)},x'_{m})< 
2\epsilon, $$
as required. \cvd

\begin{theorem}\label{co}
Any alternative Cauchy sequence $\sigma$ contains some subsequence
$\sigma'$ which is a Cauchy sequence, and $\sigma\sim \sigma'$
(i.e. (\ref{sim}) holds). Moreover, $\sigma$ converges if and only
if so does $\sigma'$, and the natural projection
$$
({\rm Cau}_{{\rm alt}}(M)/\sim ) 
\rightarrow \left({\rm Cau}(M)/\sim\right) ,
 \quad [\sigma] \mapsto [\sigma']$$
is well-defined and  bijective.
\end{theorem}
{\it Proof.} Let $\sigma=\{x_n\}$, and define inductively
 the subsequence $\sigma'$  as follows. For $k=0$, put $n_0=m(n_0)=1$ and, for each integer $k\geq 1$
take $\epsilon=1/k$. From Lemma \ref{claim}  there exists $n_k\geq
m(n_{k-1})$ and $m(n_k)\geq n_k$ such that $d(x_{n_k},x_{m})<1/k$
for all $m\geq m(n_k)$. Then, the subsequence
$\{x_{n_{k}}\}_{k}\subset \{x_{n}\}_{n}$ is a Cauchy sequence
according to Defn. \ref{ff}. In fact, for every $k$, if $k_2\geq
k_1\geq k$ then $n_{k_1}\geq n_k$ and $n_{k_2}\geq
m(n_{k_2-1})\geq m(n_{k_1})$, and so
$d(x_{n_{k_1}},x_{n_{k_2}})<1/k$.  The remainder is
straightforward. \cvd

\smallskip

\smallskip

\begin{remark}\label{rco} {\rm Theorem \ref{co} ensures that, if we
define the (forward) Cauchy completion by using alternative Cauchy
sequences instead of Cauchy sequences, the resulting space remains
the same. This shows the consistency of our approach against
possible alternatives.

Notice that the consistency also includes the convergence of
sequences. In particular, the Finslerian version of Hopf-Rinow
Theorem also holds if, instead of considering  Cauchy sequences,
one uses alternative Cauchy sequences.

In what follows, we will use Cauchy sequences instead of
alternative ones, as the first set is smaller (${\rm
Cau}(M)\subset {\rm Cau}^{\rm alt}(M)$) and more classical.
}\end{remark} The next result  summarizes the natural consistency
between the notion of Cauchy sequence and the relation of
equivalence $\sim$. It is straightforward from Prop. \ref{paux}
(1) and Theorem \ref{co}.
\begin{corollary} Let $\sigma \in \cauchyd$ and $\sigma'$ any sequence in $M$.
If the double limits (\ref{sim}) hold, then $\sigma'$ contains
some subsequence $\sigma''\in\cauchyd$ satisfying $\sigma \sim
\sigma''$. In this case, $\sigma$ is convergent if and only if so
is $\sigma''$, and the limits of both sequences coincide.
\end{corollary}

\subsubsection{The symmetrized Cauchy boundary}

Even though the forward and backward boundaries have been defined
in terms of different objects,  the following technical result
allows to establish a precise relation between them. We always
consider $\sigma=\{x_{n}\}$, $\sigma'=\{x'_{n}\}$.

\begin{proposition}\label{p2-} If $\sigma, \sigma'\in \cauchyd$ are
equivalent for $d$ and $\sigma\in \cauchydr$, then $\sigma'\in
\cauchydr$ and $\sigma$, $\sigma'$ are also equivalent for $\rd$.
\end{proposition}
{\it Proof.} To prove  $\sigma'\in \cauchydr$, take $\epsilon>0$
and let $n_1$, $m_1(n)$, $n_3$, $m_3(n)$ be the natural numbers
obtained by applying Lemma \ref{claim} to the equalities
\[
\lim_{n}(\lim_{m}d(x_{n},x'_{m}))=0, \quad
\lim_{n}(\lim_{m}d(x'_{n}, x_{m}))=0\quad\hbox{(recall that
$\sigma\sim \sigma'$ for $d$).}
\]
Since $\sigma\in \cauchydr$, there exists $n_{2}$ such that if
$m\geq n\geq n_{2}$ then
\[
d(x_{m},x_{n})<\epsilon. 
\]
Let $n_{0}:=max\{n_{1},n_{2}\}$,
$N_{0}:=max\{n_{3},m_{1}(n_{0})\}$. For any $m\geq n\geq N_{0}$,
take some $k\geq max\{m_{3}(m),n_{0}\}$. Then,
\[
d(x'_{m},x'_{n})\leq
d(x'_{m},x_{k})+d(x_{k},x_{n_{0}})+d(x_{n_{0}},x'_{n})<
3\epsilon,
\]
and $\sigma'\in \cauchydr$. To prove that $\sigma$, $\sigma'$ are
related for $\rd$, consider the natural numbers $n_1, m_1(n),
n_2$, $n_{0}=max\{n_{1},n_{2}\}$ obtained as above, and let
$N'_{0}:=max\{n_{0},m_1(n_{0})\}$. If $m\geq n\geq N'_{0}$ then
\[
d(x_{m},x'_{n})\leq
d(x_{m},x_{n_{0}})+d(x_{n_{0}},x'_{n})<2\epsilon.
\]
Therefore, $\lim_{n}(\lim_{m}\rd (x'_{n},x_{m}))=0$. The other
limit follows by interchanging the roles of both sequences.
\cvd

\smallskip

\noindent With this result in mind, one can identify points in
$\bcaum$ with points in $\bcaumr$ by requiring that they are
represented by the same class of Cauchy sequences. In particular,
the following natural definition makes sense:
\begin{definition} \label{dbs} Let $(M,d)$ be a generalized metric space. The {\em
symmetrized Cauchy  boundary} of $M$ is the intersection of its
two Cauchy boundaries: $$\bcaums:=\bcaum\cap \bcaumr\quad
(\hbox{and then,}\;\; M^{s}_{C}:=M\cup \bcaums).$$
\end{definition}
Easily, $\bcaums$ can be also computed in terms of the
symmetrized distance  in (\ref{symdist}):
\begin{proposition} \label{prop11}
A  sequence $\sigma$ is Cauchy for $\ds$ if and only if it belongs
to $\bcaums$. Thus, the symmetrized Cauchy boundary $\bcaums$ is
equal to the Cauchy boundary for the symmetrized distance $\ds$.
\end{proposition}

\brema \label{r11} {\rm For a Riemannian metric, one has a unique
 Cauchy boundary, namely $\bcaums$, which is associated to a
length space (see Remark \ref{r1}). Nevertheless,  in the
Finslerian case, the symmetrized distance $d^s$  is not associated
to a length space (see \cite[Appendix]{CJS} and \cite[Example
2.3]{CJS} for details). So,  the construction of
$\partial^{s}_{C}M$ cannot be carried out  by means of lengths of
curves in the natural way. This prevents us to consider the
symmetrized Cauchy boundary $\partial^{s}_{C}M$ as a sort of
Riemannian boundary in general. }\erema The next simple example
justifies some of the cautions adopted with double limits along
this section, and motivates the next one.



\bexam\label{illus} {\em In this  Finslerian example we show two
Cauchy sequences in $\partial^{s}_{C}M$ satisfying that only one
of the double limit equalities in (\ref{sim}) holds.
Let $D_*\subset \R^2$ be the open disk of radius 1 where the
segment $[-1,0]$ of the horizontal axis has been removed. Define
$F= \sqrt{dr^2+r^2d\theta^2+(1-r)^2d\theta^2}-
(1-r)d\theta$  in polar coordinates $r, \theta$. 
Consider  the sequences
$\{x_n=(r=1/n,\theta=0)\}_n,
\{x'_n=(r=1/n,\theta=\frac{\pi}{2})\}_n$, whose classes of
equivalence clearly belong to $\bcaums$. Note that:
$$\lim_{n}(\lim_{m}d_{F}(x_{n},x'_{m}))=0,\quad \quad
\lim_{n}(\lim_{m}d_{F}(x'_{n},x_{m}))\geq\pi/4.$$ In fact, for the
first limit consider arc curves (of constant radius) from $x_n$ to
$x'_n$, and then, radial curves from $x'_{n}$ to $x'_{m}$. For the
second one, consider any curve $\gamma:[a,b]\rightarrow D_*,
\gamma(s)=(r(s),\theta(s))$ joining $x'_n$ and $x_m$, and let
$J:=\{s\in [a,b]: \dot{\theta}(s)\leq 0\}$. Then,
$[0,\frac{\pi}{2}]\subset \theta(J)$, and so,
$${\rm length}_F(\gamma)\geq \int_J F(\dot{\gamma})=
\int_J \left(\sqrt{\dot{r}^2+r^ 2\dot{\theta}^2+(1-r)^2\dot{\theta}^2}-(1-r)\dot{\theta}\right)>\frac{1}{2}\int_J |\dot{\theta}| \geq \frac{\pi}{4}.$$  }\eexam

\subsection{Quasi-distance and topology on the Cauchy completions}
Example \ref{illus} also shows that, if a generalized distance $d$
is extended naturally to $\caum$, then such an extension may not
be a generalized distance. However, it will be still a
quasi-distance (Defn. \ref{gms}):
\begin{proposition}\label{ppp1'} The map $d_Q:\caum \times \caum
\rightarrow [0,\infty]$ given by
\begin{equation}\label{defdistbar}
d_Q([\{x_n\}], [\{y_n\}]):=\lim_{n}(\lim_{m}d(x_{n},y_{m}))
\end{equation}
is well-defined and becomes a quasi-distance on $\caum$.

Analogously, $\rd$ extends to a quasi-distance $\rd_Q$ on
$\caumr$, and $d^s$ extends to a distance $\overline d^s: \caums
\times \caums \rightarrow \R$, which satisfies $\overline
d^s=(d_{Q}+\rd_Q)/2$ on $\caums$.



\end{proposition}
{\it Proof.} For the assertions on $d_{Q}$ (the assertions on
$d_{Q}^{{\rm rev}}$ are proved analogously), note that, from Lemma
\ref{lll1}, the limits type $\lim_{n}(\lim_{m}d(x_{n},y_{m}))$
always exist. In particular, they still exist if the sequences
$\{x_{n}\}$, $\{y_{n}\}$ are replaced by subsequences
$\{x_{n_k}\}_{k}$, $\{y_{n_k}\}_{k}$, and
$\lim_{n}(\lim_{m}d(x_{n},y_{m}))=\lim_k d(x_{n_k},y_{m_k})$.

If $\{x_{n}\}\sim \{x'_{n}\}$ and $\{y_{m}\}\sim \{y'_{m}\}$, take
subsequences $\{x_{n_{k}}\}_k$, $\{x'_{n'_{k}}\}_k$ and
$\{y_{m_{l}}\}_l$, $\{y'_{m'_{l}}\}_l$ such that
$\lim_{k}d(x_{n_{k}},x'_{n'_{k}})=0, 
\lim_{l}d(y'_{m'_{l}},y_{m_{l}})=0$. Then,
\[
\begin{array}{rl}
\lim_{n}(\lim_{m}d(x_{n},y_{m})) &
=\lim_{k}(\lim_{l}d(x_{n_{k}},y_{m_{l}}))
\\ & \leq
\lim_{k}\left(\lim_{l}(d(x_{n_{k}},x'_{n'_{k}})+d(x'_{n'_{k}},y'_{m'_{l}})+d(y'_{m'_{l}},y_{m_{l}}))\right)
\\ & =\lim_{k}(\lim_{l}d(x'_{n'_{k}},y'_{m'_{l}})) \\ &
=\lim_{n}(\lim_{m}d(x'_{n},y'_{m})),
\end{array}
\]
and analogously for the reversed inequality. Therefore, $d_Q$ is
well-defined.

In order to prove that $d_Q$ is a quasi-distance on $M^{+}_{C}$,
first observe that the equivalence $d_Q(x,y)=d_Q(y,x)=0$ iff $x=y$
becomes obvious; in fact, the former holds iff $\{x_{n}\}\sim
\{y_{n}\}$, for $x\equiv [\{x_{n}\}]$, $y\equiv [\{y_{n}\}]$, i.e.
iff $x=y$. For the triangle inequality, consider $x\equiv
[\{x_{n}\}]$, $y\equiv [\{y_{m}\}]$, $z\equiv [\{z_{l}\}]$ and
take subsequences $\{y_{m_{n}}\}_n$, $\{z_{l_{n}}\}_n$ such that
\[
d_Q(x,y)=\lim_{n}d(x_{n},y_{m_{n}}),\;\;
d_Q(x,z)=\lim_{n}d(x_{n},z_{l_{n}}),\;\;
d_Q(y,z)=\lim_{n}d(y_{m_{n}},z_{l_{n}})
\]
(for the last one, take a subsequence of the previous $l_n$ and
rename it). Then,
\[
d_Q(x,z)=\lim_{n}d(x_{n},z_{l_{n}})\leq
\lim_{n}(d(x_{n},y_{m_{n}})+d(y_{m_{n}},z_{l_{n}}))=d_Q(x,y)+d_Q(y,z).\;
\]
Finally, the assertions on $d^s$ follow easily from Prop.
\ref{prop11}, and the usual properties of distances. \cvd

\begin{figure}
\centering
\ifpdf
  \setlength{\unitlength}{1bp}%
  \begin{picture}(243.45, 202.54)(0,0)
  \put(0,0){\includegraphics{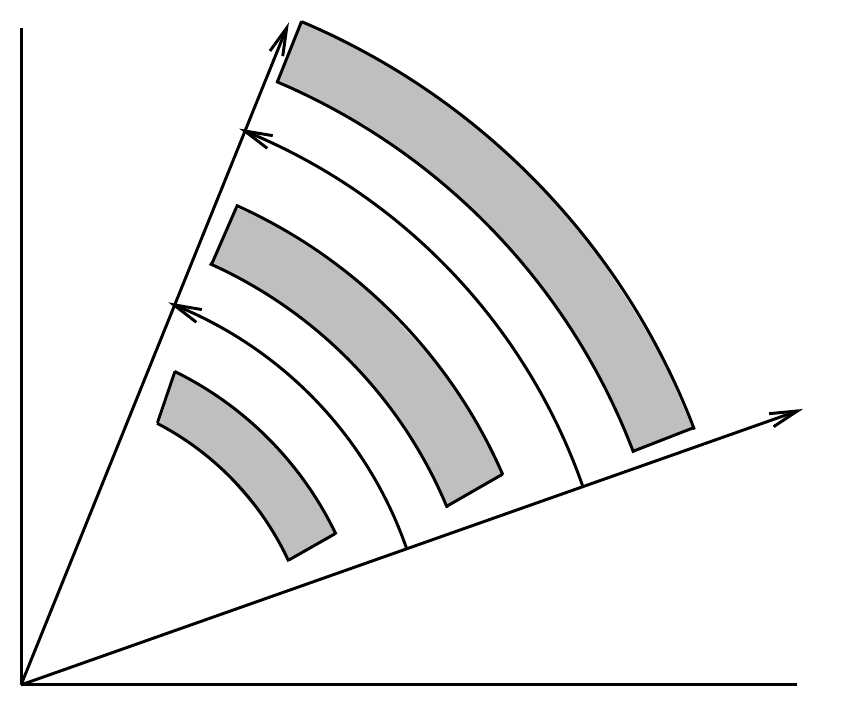}}
  \put(210.07,68.89){\fontsize{11.38}{13.66}\selectfont $c_1$}
  \put(60.44,185.02){\fontsize{11.38}{13.66}\selectfont $c_2$}
  \put(124.58,110.01){\fontsize{8.54}{10.24}\selectfont $\gamma_n$}
  \end{picture}%
\else
  \setlength{\unitlength}{1bp}%
  \begin{picture}(243.45, 202.54)(0,0)
  \put(0,0){\includegraphics{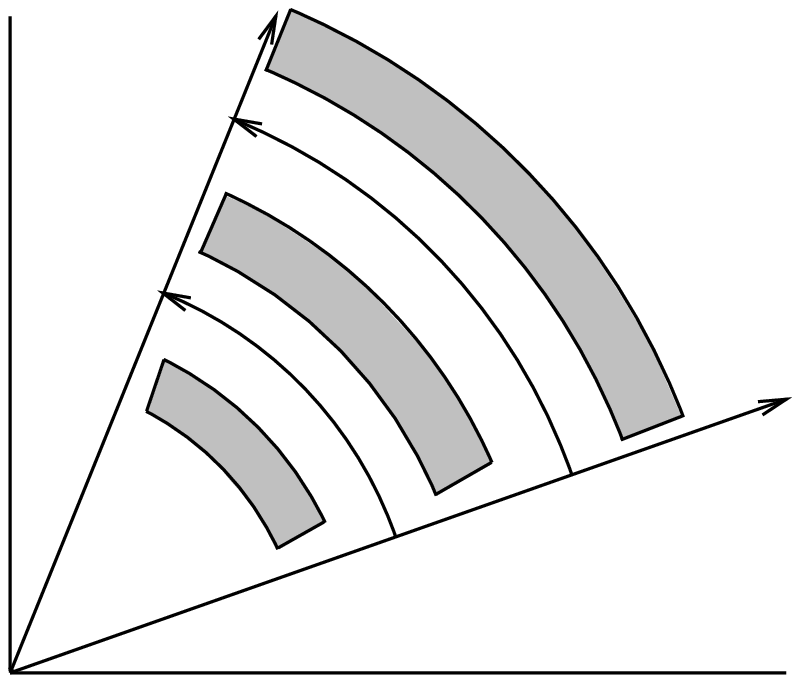}}
  \put(210.07,68.89){\fontsize{11.38}{13.66}\selectfont $c_1$}
  \put(60.44,185.02){\fontsize{11.38}{13.66}\selectfont $c_2$}
  \put(124.58,110.01){\fontsize{8.54}{10.24}\selectfont $\gamma_n$}
  \end{picture}%
\fi  \caption{\label{fig1} Limit situation for the quasi-distance
$d_Q$
(see Remark \ref{r324}). \newline An open subset $M$ of $\R^2$
(the region in grey is removed) with a Randers metric
$F=\sqrt{dx^2+dy^2+\omega^2}+\omega$. The radial dependence
assumed for $\omega$ makes the rays $c_1, c_2$ to have finite
length in the direction of the arrows; so, they  determine some
$z_1,z_2\in \bcaum$. The length is infinity in the opposite
direction. The angular dependence assumed for $\omega$ makes the
lengths of curves $\gamma_n$ to approach to $0$ in the direction
of the arrows and to infinite in the opposite direction. So,
$d_Q(z_1,z_2)=0$, $d_Q(z_2,z_1)=\infty$.}
\end{figure}

\begin{remark}\label{r324} {\em


(1) As $d_Q$ is only a quasi-distance, there are two natural topologies on $\caum$: the one generated by
the forward open balls of $d_Q$, and the one generated by the
backward open balls.
If $\{z_n\}\subset \caum$ then: 
\begin{equation}\label{ebt}
\begin{array}{lccl}
z_n\rightarrow z\in\caum & \hbox{{\rm (forward balls topology)}} &
\Leftrightarrow & d_Q(z,z_n)\rightarrow
0 \\
z_n\rightarrow z\in\caum & \hbox{\rm (backward balls topology)} &
\Leftrightarrow &  d_Q(z_n,z)\rightarrow 0.
\end{array}
\end{equation}
These two topologies are different in general. In fact, Example
\ref{illus} shows two points $z=[\{x_n\}], z'=[\{x'_n\}] \in
\partial_C^s M \subseteq \bcaum$ such that $d_Q(z,z')=0$ and
$d_Q(z',z)>0$. So, any forward ball containing $z$ also contains
$z'$ (and any backward ball containing $z'$ also contains $z$) but
the converse does not hold. In particular,  the topology in
$\caum$ generated by the backward balls
is not $T_1$ ---even though it is necessarily $T_0$. Note also
that $\overline d^s$ becomes a true distance (Prop. \ref{ppp1'}),
in spite of the fact that
$z,z'\in \bcaums$. 

(2) The character of  quasi-distance for $d_Q$ is
sharpened in Fig. \ref{fig1}.
In this figure, the Randers manifold $(M,F)$ admits two boundary
points $z_1,z_2 \in  \bcaum $ such that $d_Q(z_1,z_2)=0$ and
$d_Q(z_2,z_1)=\infty$. So, again the forward and backward
topologies on $\caum$ are different and non-$T_1$. Note also that
$\overline d^s$ remains finite, as $\caums=M$.

 (3) Even though the names suggest  that the topology generated by the
forward open balls is the natural one in $\caum$, this is not the
case, because of the following completeness property: {\em with
the  topology generated by the backward open balls, every
(forward) Cauchy sequence $\{x_{n}\} \subset M$ converges to the
point $z=[\{x_{n}\}]$ which defines in $ M_{C}^{+}$}. In fact,
$$\lim_{n}d(x_{n},z)=\lim_{n}(\lim_{m}d(x_{n},x_{m}))=0,$$
in agreement with the second equivalence in (\ref{ebt}). In
general, this property does not hold for the topology generated by
the forward balls (see Examples \ref{r2}, \ref{illus}).

{\cambios (4) As a relevant difference with the case of generalized
distances, a convergent sequence for  a quasi-distance (with the
above natural topology) may not be Cauchy. In fact, consider again
in Figure \ref{fig1} the sequence $\{y_{n}\}$ given by
$y_{2n+1}=c_{1}(n)$ (the starting point of $\gamma_n$) and
$y_{2n}=c_{2}(n)$ (the endpoint). As $c_1$ defines $z_{1}$,
$c_{2}$ defines $z_{2}$ and $d(z_{1},z_{2})=0$, the sequence
$\{y_{n}\}$ satisfies $\lim_{n}d(y_{n},z_2)=0$, even though
$\{y_n\}$ is not Cauchy, as
$d(y_{2n},y_{2n+1})$ does not converge to 0.}}
\end{remark}

\bconv \label{conv} {\em The (forward) Cauchy completion $\caum$
will be regarded as a topological space with the topology
generated by the backward balls for the quasi-distance $d_Q$ in
(\ref{defdistbar}). So, the convergence of a sequence in $\caum$
is characterized by the second relation in (\ref{ebt}), and each
$\{x_n\}\in $ Cau$(M)$ converges to $[\{x_n\}]\in \caum$. } \econv

At any case, if $d_Q$ is a generalized distance, everything is simplified:
\begin{proposition}\label{genersym}
If $d_Q$ is a generalized distance, then $\partial_C^\pm M=\partial_C^s M$.
\end{proposition}
{\it Proof.} Consider $z=[\{x_n\}]\in \partial_C^+M$ and observe
that $\lim_n d_Q(x_n,z)=0$ (recall Remark \ref{r324} (3)). From the
generalized character of $d_Q$, $\lim_n d_Q(z,x_n)=0$ and, from
definition of $d_Q$, $\lim_n \lim_m d(x_m,x_n)=0$, i.e. the
sequence $\{x_n\}$ is an alternative backward Cauchy sequence.
Finally, from Theorem \ref{co}, there exists a subsequence
$\{x_{n_k}\}$ which is a backward Cauchy sequence, i.e.
$z=[\{x_{n_k}\}]\in
\partial_C^- M\cap \partial_C^+ M=\partial_C^s M$. \cvd

\btheo \label{tcompltcompl} The (forward) Cauchy completion $\caum$ of a generalized
 metric space $(M,d)$ endowed with its natural (backward) topology in Convention \ref{conv} satisfies:

(i) The topology is $T_0$, and
if two distinct points of $\caum$ are not Hausdorff related then
they both lie at the boundary $\bcaum$.

 (ii) The  inclusion $i:
M\rightarrow \caum$ is a homeomorphism onto its image $i(M)$, and
$i(M)$ is dense in $\caum$.

(iii) The quasi-distance $d_Q$ is {\em forward complete}, i.e. any
forward Cauchy sequence $\{z_n\}\subset \caum$  is convergent in
$\caum$.

The analogous statements hold for $\rd_Q$  on $\caumr$, while
$(\caums, \overline d^s)$ is a (necessarily $T_2$) metric space.
\etheo {\it Proof.} (i) If $x, y\in \caum$ are not $T_0$-separated, then
$y\in B^-(x,\frac{1}{n})$ and $x\in B^-(y,\frac{1}{n})$ for all
$n$. Thus, $d_Q(x,y)=d_Q(y,x)=0$ and $x=y$, as $d_Q$ is a
quasi-distance.

Now, let $x,y\in \caum$ be non-Hausdorff related, and assume $x\in
M$. Taking into account that $B^-(x,\frac{1}{n})\cap
B^-(y,\frac{1}{n})\neq \emptyset$ for all $n\in \N$, consider a
sequence $\{x_n\}\subset M$ with $x_n\in B^-(x,\frac{1}{n})\cap
B^-(y,\frac{1}{n})$. Then,  $\lim_n d_Q(x_n, x)=0$ and, as $d_Q$
agrees with the generalized distance $d$ on $M$,
 $\lim_n d_Q(x,x_n)=0$. Therefore
 $$0\leq d_Q(x,y)\leq d_Q(x,x_n)+d_Q(x_n,y) \rightarrow 0\quad\hbox{as $n\rightarrow\infty$}.$$
Therefore, if $y=[\{y_n\}]$, necessarily
$0=d_Q(x,y)=\lim_n d_Q(x,y_n)$, i.e. $\{y_n\}$ converges to $x$
with the topology generated by the forward balls. As this topology
coincides on $M$ with the one generated by the backward balls,
one also have $\lim_n d_Q(y_n,x)=0$, i.e. $x$ is identifiable to
$[\{y_n\}]=y$.

(ii) Straightforward.

(iii) Let $\{z_n\}_n\subset \caum$, $z_n=[\{z_n^i\}_i]$ be a
Cauchy sequence. There is no loss of generality if we suppose that
$d_Q(z_n,z_{n+1})<\frac{1}{n^2}$ for all $n$. In fact, otherwise
some subsequence $\{z_{n_k}\}_{k}\sim\{z_n\}_n$ would satisfy this
property, and Proposition \ref{propsim}(2) would be applicable (recall Remark \ref{rquasi}). Consider the
sequence $\{x_n\}$ defined as follows. For each $n$, Lemma
\ref{lll1}
implies that, up to subsequences, there exist $i_n,j_n$ such that,
for $i\geq i_n, j\geq j_n, d(z_n^i,z_{n+1}^j)<1/n^2$. Then, choose
$x_1=z_1^{i_1}$ and $x_n=z_n^{k(n)}$ where $k(n)=$ Max
$\{i_{n},j_{n-1}\}$.

One has to prove that $\{x_n\}$ is a Cauchy sequence, and the
point $x=[\{x_n\}]\in \caum$ satisfies $\lim_n d_Q(z_n,x)=0$. For
the first property, just notice that
$d(x_n,x_{n+1})=d(z_n^{k(n)},z_{n+1}^{k(n+1)})<\frac{1}{n^2}$,
and, so, $d(x_n,x_{n+h_0})< \sum_{h=0}^\infty (n+h)^{-2}$. The
second one follows from the following estimation:
\smallskip

$\begin{array}{rl} d_Q(z_n,x)= & \lim_i (\lim_m d(z_n^i,x_m)) \\
\leq & \lim_i (\lim_m
(d(z_n^i,x_{n+1})+\sum_{j=n+1}^{m-1}d(x_j,x_{j+1})))\\= & \lim_i(
d(z_n^i,z_{n+1}^{k(n+1)}))+\lim_m(\sum_{j=n+1}^{m-1}d(x_j,x_{j+1}))\\
< & \frac{1}{n^2}+\lim_m \sum_{j=n+1}^{m-1}\frac{1}{j^2}
\\
= &
\sum_{j=n}^{\infty}\frac{1}{j^2}\end{array}.$

Finally, for the last assertion recall Prop. \ref{ppp1'}. \cvd



\smallskip

For Riemannian manifolds the Cauchy boundary is a closed subset of
the Cauchy completion. Of course, this property cannot be expected
for an arbitrary metric space $(X,d)$ (as its Cauchy completion is
also the Cauchy completion of any dense, non necessarily open,
subset $A\subset X$). But such a property can be extended to the
Finsler case (as well as more general spaces as those under
Conventions
 \ref{cginebra}, \ref{cginebra2}
below, as local compactness will be the essential property needed
for the proof).

\bprop \label{p329} For any Finsler manifold $(M,F)$, the boundary
$\bcaum$ is a closed subset of $\caum$. \eprop {\it Proof.} Notice
that, for any $x\in M$, any closed ball centered at $x$ with small
radius $r>0$ is a compact subset of $M$ and, so, cannot contain a
non-converging Cauchy sequence (nor a point of $\bcaum$). \cvd

\subsection{Relating $\caum$ and $\caumr$ through the extended quasi-distance}

Our aim here is to show that the extended quasi-distance $d_Q$ can
be extended further to a domain larger than $\caum$.

\begin{convention}\label{conv2}{\rm From now on, and whenever there is no possibility of confusion, we will drop the notational distinction between $d,
d_Q$ (resp. $\rd, \rd_Q$) and the further extensions of the
generalized distance to be defined in the next proposition. So,
$d$ and $\rd$ will be applied in each case to the biggest possible
domain. }\end{convention}

\begin{proposition}\label{Srelat} Let $(M,d)$ be a generalized metric space. Then:
\begin{itemize} \item[(i)] The maps $d:M^{+}_{C}\times
(M^{+}_{C}\cup M^{-}_{C})\rightarrow [0,\infty]$ and $\rd:
M^{-}_{C}\times (M^{+}_{C}\cup M^{-}_{C})\rightarrow [0,\infty]$,
given formally by the double limit (\ref{defdistbar}), are
well-defined. \item[(ii)] $d(x,y)=\rd (y,x)$ for all $x\in
M^{+}_{C}$, $y\in M^{-}_{C}$.
\end{itemize}
\end{proposition}
{\it Proof.} Let $x=[\{x_{n}\}]=[\{x'_{n}\}]\in M^{+}_{C}$,
$y=[\{y_{n}\}]=[\{y'_{n}\}]\in M^{-}_{C}$.

(i) From Remark \ref{rll1}, the double limits
$\lim_{n}(\lim_{m}d(x_{n},y_{m}))$ and
$\lim_{n}(\lim_{m}d(x'_{n},y'_{m}))$ exist in $[0,\infty)$. Then,
the same proof in Prop. \ref{ppp1'} shows that they coincide.



(ii) 
Clearly,
$$\begin{array}{rl}
\rd (y,x)=& \lim_{n}(\lim_{m}\rd (y_{n},x_{m})) =
\lim_{n}(\lim_{m}d(x_{m},y_{n})) \\ \leq &
\lim_{n}(\lim_{m}(\lim_{k}(d(x_{m},y_{k})+d(y_{k},y_{n})))) \\
= & \lim_{m}(\lim_{k}d(x_{m},y_{k}))+\lim_{n}(\lim_{k}\rd
(y_{n},y_{k})) \\ = & d(x,y)+\rd (y,y)=d(x,y), \end{array}$$ and
the reversed inequality holds by interchanging the roles of $x$
and $y$.
\cvd

\begin{remark}\label{motrem}{\em  {\cambios In the previous proposition, the value $\infty$ is allowed for $d$. Only when
 $x\in M$ and $x^+\in \caum$, we can ensure $d(x,x^+)<\infty$
(Lemma \ref{lll1}), and analogously for $\rd$.}

In the case of length spaces, one can connect any $x\in M$ with
$x^+ = [\{x_n\}]\in \bcaum$ by means of a curve of finite length:
in fact, for certain subsequence $\{n_i\}$, one can join $x$ with
$x_{n_{1}}$ by a curve of finite length, and each $x_{n_{i}}$ with
$x_{n_{i+1}}$ by a curve of length smaller than $i^{-2}$.
If the reversed curve of $\gamma$ had also finite length, then one
would also have $d(x^+,x)<\infty$ and $x^+\in \bcaums$. However,
even if we know that $d(x^+,x)<\infty$, we cannot ensure that
$x^+$ belongs to $\bcaums$. The reason is that $x^+$ may be
non-connectable to any $x\in M$ by means of a curve of finite
length, but all the elements of $\{x_n\}$ may be connectable with
$x$ by means of a sequence of curves of bounded length. Figure 2
shows an example of this situation.}
\end{remark}

\begin{figure}
\centering
\ifpdf
  \setlength{\unitlength}{1bp}%
  \begin{picture}(339.59, 147.66)(0,0)
  \put(0,0){\includegraphics{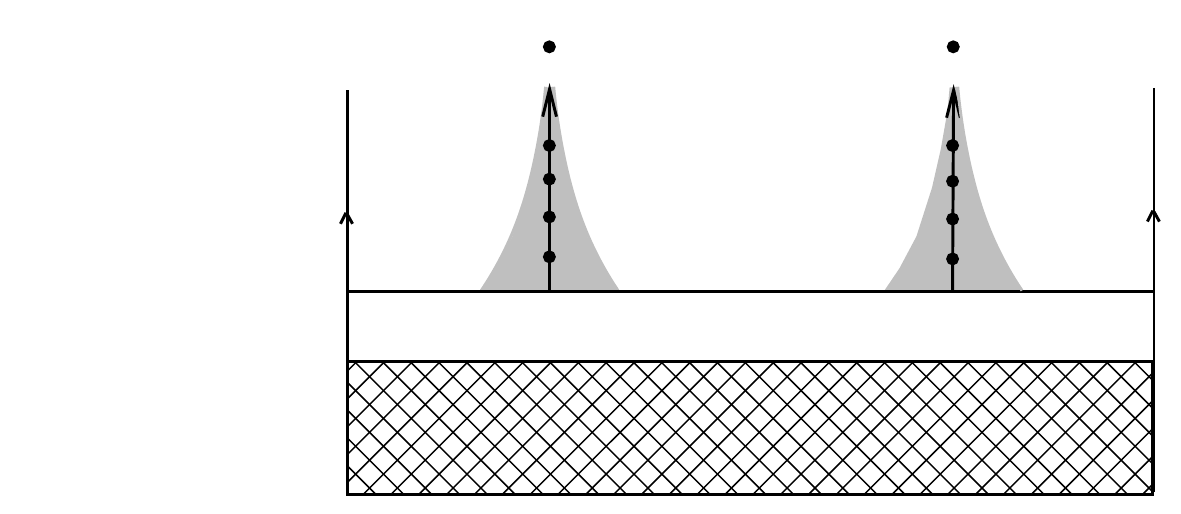}}
  \put(154.50,54.86){\fontsize{9.72}{11.66}\selectfont -3}
  \put(272.75,55.18){\fontsize{9.72}{11.66}\selectfont 3}
  \put(5.67,110.67){\fontsize{9.72}{11.66}\selectfont $\partial^{+}_{C}M\equiv \{z^{+}\}$}
  \put(5.67,81.98){\fontsize{9.72}{11.66}\selectfont $\partial^{-}_{C}M\equiv \{z^{-}\}$}
  \put(5.67,54.86){\fontsize{9.72}{11.66}\selectfont $\partial_{C}^{s} M=\emptyset$}
  \put(168.47,112.33){\fontsize{9.72}{11.66}\selectfont $c_{1}$}
  \put(279.57,112.33){\fontsize{9.72}{11.66}\selectfont $c_{2}$}
  \put(101.03,55.41){\fontsize{9.72}{11.66}\selectfont -6}
  \put(323.74,54.86){\fontsize{9.72}{11.66}\selectfont 6}
  \put(160.91,135.32){\fontsize{8.54}{10.24}\selectfont $z^+$}
  \put(278.43,134.11){\fontsize{8.54}{10.24}\selectfont $z^-$}
  \put(145.35,82.70){\fontsize{8.54}{10.24}\selectfont $x_n$}
  \put(261.68,82.70){\fontsize{8.54}{10.24}\selectfont $x'_n$}
  \end{picture}%
\else
  \setlength{\unitlength}{1bp}%
  \begin{picture}(339.59, 147.66)(0,0)
  \put(0,0){\includegraphics{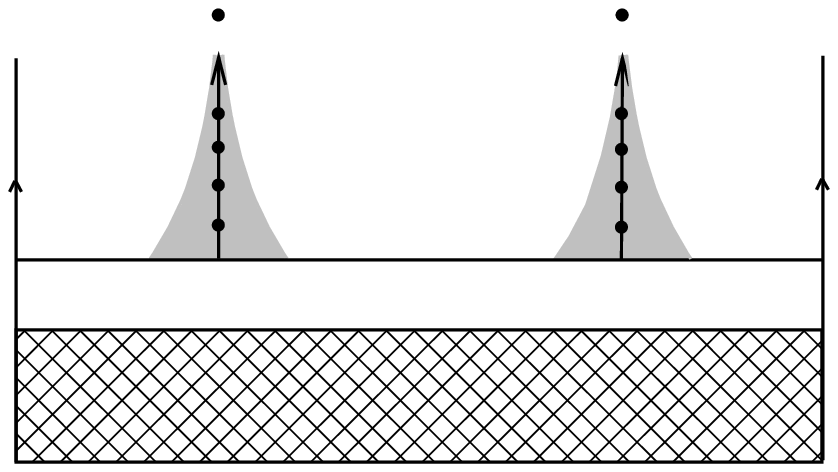}}
  \put(154.50,54.86){\fontsize{9.72}{11.66}\selectfont -3}
  \put(272.75,55.18){\fontsize{9.72}{11.66}\selectfont 3}
  \put(5.67,110.67){\fontsize{9.72}{11.66}\selectfont $\partial^{+}_{C}M\equiv \{z^{+}\}$}
  \put(5.67,81.98){\fontsize{9.72}{11.66}\selectfont $\partial^{-}_{C}M\equiv \{z^{-}\}$}
  \put(5.67,54.86){\fontsize{9.72}{11.66}\selectfont $\partial_{C}^{s} M=\emptyset$}
  \put(168.47,112.33){\fontsize{9.72}{11.66}\selectfont $c_{1}$}
  \put(279.57,112.33){\fontsize{9.72}{11.66}\selectfont $c_{2}$}
  \put(101.03,55.41){\fontsize{9.72}{11.66}\selectfont -6}
  \put(323.74,54.86){\fontsize{9.72}{11.66}\selectfont 6}
  \put(160.91,135.32){\fontsize{8.54}{10.24}\selectfont $z^+$}
  \put(278.43,134.11){\fontsize{8.54}{10.24}\selectfont $z^-$}
  \put(143.35,82.70){\fontsize{8.54}{10.24}\selectfont $x_n$}
  \put(261.68,82.70){\fontsize{8.54}{10.24}\selectfont $x'_n$}
  \end{picture}%
\fi  \caption{\label{fig2} Non-evenly pairing (forward and
backward) boundaries,
see Remark \ref{motrem}, Definition \ref{dtotsep}.
\newline Starting at a strip of the plane, the two vertical lines
at $x=\pm 6$ are identified (the dashed region will be
irrelevant). Then, a Randers metric
$F=\sqrt{dx^2+dy^2+\omega^{2}}+\omega$ is chosen such that $c_1$
(resp. $c_2$) determines the unique forward (resp. backward)
boundary point $z^+\in \bcaum$ (resp. $z^-\in \bcaumr$): namely,
the support of $\omega$ is included in the grey region and
$\omega(\partial_y)$ decreases  (resp. increases) strictly to $-1$
(resp. 1) on $c_1$ (resp. $c_2$) so that the $F$-length of $c_1$
and $-c_2$ (the reversed parametrization of $c_2$) becomes finite
but the $F$-length of $-c_1$ and $c_2$ is infinite. The distance
$d(z^+,z^-) $
is finite ---essentially, one could use horizontal lines to
connect the points in the Cauchy sequence $\{x_n\}$ which defines
$z^+$ with the points in the Cauchy sequence $\{x'_n\}$ which
defines $z^-$. Moreover, $d(z^+,x)$, $d(x,z^-)$ are finite too for
any $x\in M$ ---say, for large $n$, go from $x_n$ to $c_2$ by
using a horizontal line,
 go down by using
$c_2$ until the height of $x$, and use again an horizontal segment
to reach $x$. However, $z^+,z^-$ do not belong to $\bcaums$ (in
fact, $\bcaums =\emptyset$).}
\end{figure}



In spite of examples as the illustrated one in Figure \ref{fig2}, the
length spaces satisfying that any two points $x^+$, $x$ such that
$d(x^+,x)<\infty$ admit a curve of finite length from $x^+$ to $x$
(and thus, $x^+\in \bcaums$), constitute a natural
{\cambios subclass: say, as any
point $x^+$  of the boundary and any point $x$ in $M$ can be
connected by means of some curve in $M$ (as explained in the
previous remark), if $d(x^+,x)<\infty$, one may expect that such a
curve can be chosen of finite length.
}
  The following
definition includes, in particular, this class.
\begin{definition}\label{dtotsep}
A generalized metric space $(M,d)$ is said to have {\em forward
(resp. backward)} {\em \ras} if the following condition holds:
$$\begin{array}{rlcll}
& d(x^+,x)= \infty &  & \forall x^+\in \bcaum
\setminus \bcaums & \hbox{for some, and then any,} \; x\in \caumr \\ 
(\hbox{resp.} & d(x,x^-) = \infty &   & \forall x^-\in \bcaumr
\setminus \bcaums & \hbox{for some, and then any,} \; x\in \caum).
\end{array}$$
\end{definition}
We have emphasized in this definition that, when the first (resp.
second) condition is satisfied for some $x\in M$, then it holds
for any $x\in \caumr$ (resp. $x\in M_{C}^{+}$). So, in the
practice, it suffices to check them for some point of $M$. {\cambios
For example, if   Figure \ref{fig2} were modified so that the
support of $\omega$ became included  in only one of the two grey
regions, the boundaries would be clearly evenly pairing.}
Trivially, from Proposition \ref{genersym}:

\begin{corollary}\label{genersym2}
If $d_Q$ is a generalized distance, then $(M,d)$ has {\cambios forward
and backward evenly pairing boundaries}.
\end{corollary}

\section{Riemannian Gromov and Busemann completions}\label{s4}

Eberlein and O'Neill \cite{EO} developed a well-known
compactification for any  Hadamard manifold (or, in general, any
CAT(0)-space), based in classes of equivalence of rays. Gromov
\cite{Gr81} introduced a universal compactification for any
complete Riemannian manifold, valid also for some more general
complete metric spaces \cite[p. 184f]{Gr81}.
The latter compactification becomes equivalent to the former in
Hadamard manifolds; in order to prove this, Busemann functions
play an essential role (see, for example, the books \cite{BGS,
BH}). In this paper, we are interested in analogous boundaries for
arbitrary Riemannian and Finslerian manifolds.

We start by reviewing briefly Gromov's construction for any
Riemannian manifold, including the non-complete case, with a
threefold aim.
First, to make clear in which sense the Cauchy boundary can be
regarded as a part of Gromov's one. Second, to show the appearance
of the Busemann boundary as an intermediate boundary between the
Cauchy and Gromov's ones. And third, to introduce the precise
framework to be used in the Finslerian case and the remainder of
the paper.

Our objective is to study a metric space $(M,\dri)$, being $M$ a
(connected) smooth manifold and $\dri$ the distance associated to
a Riemannian metric $g_R$ on $M$. However, the results will be
stated in a more general setting, concretely:

\begin{convention}\label{cginebra} {\em Along the present Section \ref{s4},
$(M,\dri)$ will denote any metric space which satisfies: $M$ is
locally compact with a countable dense set (or equivalently, which
is second countable), and $\dri$ is derived from a length space.
Recall that the metric space associated to any {\em reversible}
Finsler manifold is included. Typical notations in manifolds are
extended to this general setting with no further mention. In
particular, if $c$ is a smooth curve, the role of Riemannian norm
$|\dot c (t)| =\sqrt{g_R(\dot{c}(t),\dot{c}(t))}$ in a general
length space will be played by the {\em dilatation}. In general,
we will adopt the terminology in \cite[Chapter 1]{Gr}.
}\end{convention}
\begin{remark}\label{remark8s}{\em As above,
$\partial_CM$ denotes the  Cauchy boundary, and $M_C$ the Cauchy
completion, the latter endowed with the continuous extension of
$\dri$ (denoted again by $\dri$). Notice that $(M_C, \dri)$ is
clearly a length space, but it may not lie under Convention
\ref{cginebra} as $M_C$ may not be locally compact, even if
$(M,d^{R})$ comes from a Riemannian manifold.}\end{remark}

\subsection{Gromov completion}\label{s41} For each  $x\in M_C$,
let $\dri_x$ be the Lipschitz function   on $M$ given by $y\mapsto
\dri(y,x)$. Consider the space $\lip$ of all the Lipschitz functions
 (always assumed with  Lipschitz
constant equal to 1) on $M$, endowed with the topology of
pointwise convergence. Notice that this
space is equivalent to the space of Lipschitz functions on $M_C$;
in particular, functions $f\in \lip$ will be extended to $M_C$
with no further mention.

\begin{remark}\label{rmetrizable}{\em
The topology of pointwise convergence on $\lip$ coincides with the
topology of uniform convergences on compact sets and the
compact-open topology; moreover, it is metrizable. In fact,
consider the space ${\cal C}(M)$ of all the continuous functions
from $M$ to $\R$. In this space the topology of uniform
convergence on compact sets coincides with the compact-open
topology (see \cite[p. 230]{KE}) and, when restricted to $\lip$,
with the pointwise convergence topology (see \cite[p. 232]{KE}).
It is not hard to prove that ${\cal C}(M)$ is Hausdorff, regular
and second countable (as so are $M$ and $\R$).
{\cambios So, by Urysohn
theorem, ${\cal C}(M)$, and then $\lip$, are metrizable. Even
more, the following metric $d^1$ on $\lip$ generates the topology.
Choose $x_0\in M$ and an auxiliary complete metric $\tilde{d}^R$
such that $\tilde{d}^R\geq d^R$, and define the metric}
\begin{equation}\label{emetrizable}{\cambios
d^1(f_1,f_2):={\rm sup}_{x\in M}
\frac{|f_1(x)-f_2(x)|}{1+\tilde{d}^R (x,x_0)^2}\quad \quad
\forall\;f_1,f_2\in \lip. }
\end{equation}
 Let $f_n,f\in \lip$, $n\in \N$. Clearly, if $\{d^1(f_n,f)\}$
converges to $0$, then the sequence $\{f_n\}$ converges pointwise
to $f$. For the converse,   {\cambios the Lipschitz condition
yields the following bound for the expression corresponding to the
fraction in (\ref{emetrizable}):
\[\frac{|f(x_0)|+|f_n(x_0)|+2\;d^R
(x,x_0)}{1+\tilde{d}^R(x,x_0)^2} \leq
2\frac{1+|f(x_0)|+\tilde{d}^R (x,x_0)}{1+\tilde{d}^R(x,x_0)^2},
\]} (the last inequality for large $n$). So, the uniform convergence of $\{f_n\}$ on $\tilde{d}^R$-bounded sets
completes the result.
}\end{remark}

\begin{lemma}\label{l41} Let $\{f_n\}$ be a sequence in $\lip$.
Assume that there exists $x_0\in M_C$ such that $\{f_n(x_0)\}$
converges to some $k_0\in[-\infty,\infty]$. If $k_0\in \R$ then
$\{f_n\}$ admits a subsequence converging in $\lip$. If $k_0=\pm
\infty$ then $\{f_n\}$ converges uniformly on compact sets to $\pm
\infty$.
\end{lemma}
{\em Proof}. Only the case $k_0\in\R$ will be considered, as
$k_0=\pm \infty$ is similar.

Let $\{x_{n}\}_{n=0}^{\infty} \subset
M$ be a countable dense subset of $M$. Remove a finite subset of
$\{f_n\}$, to obtain a subsequence  $\{f_n^0\}$ with
 $|f^0_{n}(x_{0})-k_{0}|<1$
for all $n$. From the  Lipschitz condition, one has
$$-\dri(x_{0},x_{m})-1\leq f^0_{n}(x_{m})-k_{0}\leq \dri(x_{0},x_{m})+1\quad\forall n,m.$$
Now, we can construct inductively a subsequence
$\{f_{n}^{m}\}_n$ of $\{f_{n}^{m-1}\}_n$ such that
$\{f_{n}^{m}(x_{i})\}_n$ is convergent for $i=1,\ldots,m$. Then,
the diagonal subsequence $\{f_{n}^{n}\}_n$ is convergent in
$\lip$, as it converges on each $x_m$. \cvd

\begin{proposition} \label{p41}
The map\footnote{For convenience, we change the sign of the
natural definition here. Later, this will affect only to the sign
of Busemann functions, which will be the opposite of the one in
the usual convention.} $$\hat j: M_C\rightarrow \lip, \quad \quad
x\mapsto -\dri_x$$ is a topological embedding of $M_C$ in $\lip$.

If $K\subset M_C$ is any bounded subset
then the closure of $\hat j(K)$ in $\lip$ is compact.
\end{proposition}
{\it Proof.} The first assertion is standard {\cambios (simplify Lemma \ref{lcf} below).
For the last one, we only
have to show that the closure of $\hat j(K)$ is sequentially
compact, since $\lip$ is metrizable by Remark
\ref{rmetrizable}
By Lemma \ref{l41}, it suffices to check that
any sequence $\{f_{n}\}$ in the closure
of $\hat j(K)$ admits some $x\in M_C$ such that
$\{|f_{n}(x)|\}$ is bounded. Notice that $f_n =\lim_i \hat
j(y^i_n)$ for some sequence $\{y_n^i\}_{i=0}^\infty \subset K$.
Hence, for any $x\in M_C$:
$$
\begin{array}{rl}
|f_n(x)-(\hat j(y_n^0))(x)| = & |\lim_i (\hat j(y_n^i))(x)-(\hat
j(y_n^0))(x)| \\ \leq & \lim_i |\dri(x,y_n^i)-\dri(x,y_n^0)| \leq \lim_i
\dri(y_n^0, y_n^i).\end{array}$$ But the last term is bounded by the
diameter of $K$, as required. \cvd

\smallskip

$\lip_*=\lip/\R$ will be the topological quotient of $\lip$ by the
1-dimensional subspace of
the (real) constant functions. The natural quotient topology will be called {\em Gromov} or {\em (quotient) pointwise topology}: 
\begin{proposition}\label{pemb}
The composition map 
$$j: M_C\rightarrow \lip_*, \quad
\quad x\mapsto [-\dri_x]$$ is continuous and injective.

Moreover, the restriction $j|_M: M\rightarrow \lip_*$ is a
topological embedding.
\end{proposition}
{\it Proof.} 
The continuity of $j$ follows from the continuity of $\hat j$, and
the injectivity from the following property: if $x\neq y$ then
$\dri_x-\dri_y$ changes its sign at $x, y$, and so, cannot be constant.
So, one only has to prove the continuity of $j^{-1}$ on $j(M)$.

To this aim, let $x_n,x\in M$ such that $\{j(x_n)\}\rightarrow
j(x)$. Then, there exists a sequence $\{t_{n}\}\subset\R$ such
that $\{t_{n}-\dri_{x_{n}}\}_{n} \rightarrow -\dri_x$ with the pointwise
topology. Assume by contradiction that $\{x_{n}\}\not\rightarrow
x$. Up to a subsequence, there exists $\epsilon>0$ small enough
such that $x_{n}\not\in \overline{B(x,2\epsilon)}$ for all $n$. We
can also assume that the spheres $S(x,r):=\{y\in M:\,\dri(y,x)=r\}$
are compact for $r=\epsilon, 2\epsilon$. Notice that, for any
$y\in S(x,2\epsilon)$ and $z\in S(x,\epsilon)$,\footnote{The first
equality is justified  by the existence of a finite limit stated
in the second one. We will use such a posteriori justifications
with no further mention.}
\begin{equation}\label{o}
\begin{array}{rl}\lim_{n}(\dri(z,x_{n})-\dri(y,x_{n}))= & \lim_{n}(t_{n}-\dri(y,x_{n}))-\lim_{n}(t_{n}-\dri(z,x_{n}))\\
=& -\dri(y,x)+\dri(z,x)=-\epsilon <0 .\end{array}
\end{equation}
On the other hand, for each $n$ we can choose  $y_{n}\in
S(x,2\epsilon)$ and $z_{n}\in S(x,\epsilon)$ such that
\begin{equation}\label{o2}
\dri(z_{n},x_{n})-\dri(y_{n},x_n)\geq 0.
\end{equation} In fact, consider a sequence
of curves $\gamma_n^j$ from $x$ to $x_n$ with $\{{\rm
length}(\gamma_n^j)\}_{j}\rightarrow \dri(x,x_n)$. Each $\gamma_n^j$
will intersect $S(x,2\epsilon)$, $S(x,\epsilon)$ in some points
$y_n^j,z_n^j$, resp., and, with no loss of generality, one can
assume the convergence  $\{y_n^j\}_{j}\rightarrow y_n$ and
$\{z_n^j\}_{j}\rightarrow z_n$. Then, by the continuity of
$\dri(\cdot,x)$, such points $\{y_n\}$ and $\{z_n\}$ satisfy
$\dri(z_n,x_n)=\dri(z_n,y_n)+\dri(y_n,x_n)\geq \dri(y_n,x_n)$.
Finally,  the compactness of the spheres allows to assume that
$\{y_n\}, \{z_n\}$ converge, obtaining a contradiction with
(\ref{o2}), (\ref{o}). \cvd

\begin{remark}\label{mclc}{\em Notice that when $M_C$ is locally compact,
{\cambios it also lies under Convention \ref{cginebra} (Remark
\ref{remark8s}), and so, the proof of Proposition \ref{pemb} also
ensures that $j$ is an embedding on all $M_C$.}
}\end{remark}

The following technical result will be useful for future referencing.
\begin{lemma} \label{l4.5}
Consider $j:M_C\rightarrow \lip_*$ and $x_n,x\in M_C$. If
$\{[-\dri_{x_n}]\}_{n}\rightarrow [-\dri_x]$ then $$\lim_n
\left(\dri(\cdot,x)+\dri(x,x_n)-\dri(\cdot,x_n)\right)=0.$$
\end{lemma}
{\it Proof.} The limit $\{[-\dri_{x_n}]\}_{n}\rightarrow [-\dri_x]$
provides the existence of $\{t_n\}$ such that
$\{t_n-\dri(\cdot,x_n)\}_{n}$ converges pointwise to $-\dri(\cdot,x)$.
Then, if we evaluate this expression at $x$, we deduce
$t_n-\dri(x,x_n)\rightarrow 0$. Therefore,
\[
\begin{array}{c}\lim_n(\dri(\cdot,x)+\dri(x,x_n)-\dri(\cdot,x_n)) \qquad\qquad\qquad\qquad\qquad\qquad\qquad\qquad\qquad \\ \quad\qquad\qquad\qquad =\dri(\cdot,x)-\lim_n(t_n-\dri(x,x_n))+\lim_n(t_n-\dri(\cdot,x_n)) \\ =  \dri(\cdot,x)-\dri(\cdot,x)=0. \qquad\qquad\qquad \end{array}
\]\cvd
\begin{example}\label{ex43}
{\em Here we are going to use bounded (i.e. without ``directions
at infinity'') examples to illustrate
the following four properties: (a) {\cambios $j$ may be not an
embedding on all $M_C$ and $\hat j (M_C)$ may} be non-closed
in $\lip$; (b) even if $j$ is a topological embedding, $j(M_C)$
may be non-closed in $\lip_*$; (c) even if $j(M_C)$ is closed, the
map $j:M_C\rightarrow \lip_*$ may not be a topological embedding;
and (d) even if $j$ is a topological embedding, $M_C$ may not be
locally compact. These items will be illuminating to understand
the way in which the Cauchy boundary lies in the Gromov one, and
they also show the optimal character of Propositions \ref{p41} and
\ref{pemb}.

Consider the metric space $(M,\dri)$ with
\begin{equation}\label{mmprima} M:=\cup_{n=1}^\infty \{(n^{-1},y):\; y\in [0,1)\} \cup
\{(x,0):\; x\in(0,1]\}\end{equation} and $\dri(x,y)$ obtained as
the infimum of the usual lengths of piecewise smooth curves in $M$
connecting $x$ and $y$ (see Fig. \ref{peine2})\footnote{Notice
that $(M,\dri)$ fulfills the essential properties in Convention
\ref{cginebra}. Even though, for simplicity, this space is not a
manifold, the example can be easily transformed into a Riemannian
2-manifold by enlarging each vertical line in a strip (thinner as
$n$ grows) and including also the forth quadrant, $x>0, y<0$.
Analogously, $M'$ in (\ref{mprima}) below can be also transformed
into a manifold.}. The sequences type
$\{(n^{-1},y_0)\}_n$ 
with $0<y_0<1$ do not contain any convergent subsequence in
$M_{C}$; in fact, $\partial_CM$ is naturally identifiable to
$(0,0)$ plus the non-Cauchy sequence
$\{(n^{-1},1)\}_{n=1}^{\infty}$. However, for any $0\leq y_0\leq
1$, each sequence $\{-\dri_{(n^{-1},y_0)}\}_{n}$ converges to a
function in $\lip$, which will be denoted $-\dri_{(0,y_0)}$
---notice that it corresponds with the natural (minus) distance to
$(0,y_0)$ in the closure of $M$ as a subset of $\R^2$. In
particular, this shows (a).

For (c), note that the sequence of functions
$\{y_0-\dri_{(n^{-1},y_0)}\}_{n}$, $y_0\in [0,1]$, converges
pointwise to $-\dri_{(0,0)}$. {\cambios So, the sequence of classes
$\{[-\dri_{(n^{-1},y_0)}]\}$ converges to $[-\dri_{(0,0)}]$ in
$\lip_*$ (intuitively, $[-\dri_{(0,0)}]=[-\dri_{(0,y_0)}]$), and
$j(M_C)$ is closed.}

\begin{figure}
\centering
\ifpdf
  \setlength{\unitlength}{1bp}%
  \begin{picture}(290.24, 152.91)(0,0)
  \put(0,0){\includegraphics{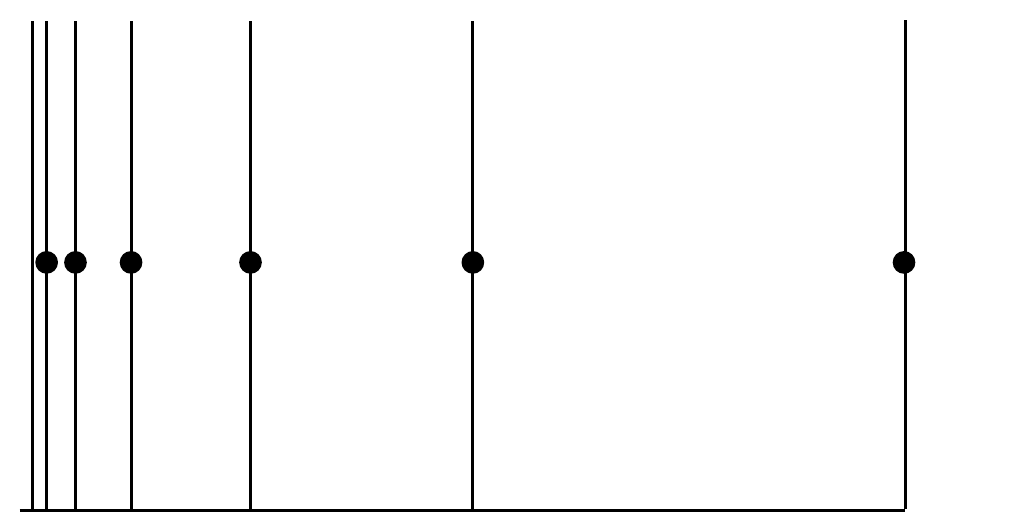}}
  \put(263.31,69.26){\fontsize{8.54}{10.24}\selectfont $x_1$}
  \put(139.14,69.63){\fontsize{8.54}{10.24}\selectfont $x_2$}
  \put(23.89,69.91){\fontsize{8.54}{10.24}\selectfont $x_n$}
 \put(-5.89,-3.01){\fontsize{8.54}{10.24}\selectfont $(0,0)$}
  \end{picture}%
\else
  \setlength{\unitlength}{1bp}%
  \begin{picture}(290.24, 152.91)(0,0)
  \put(0,0){\includegraphics{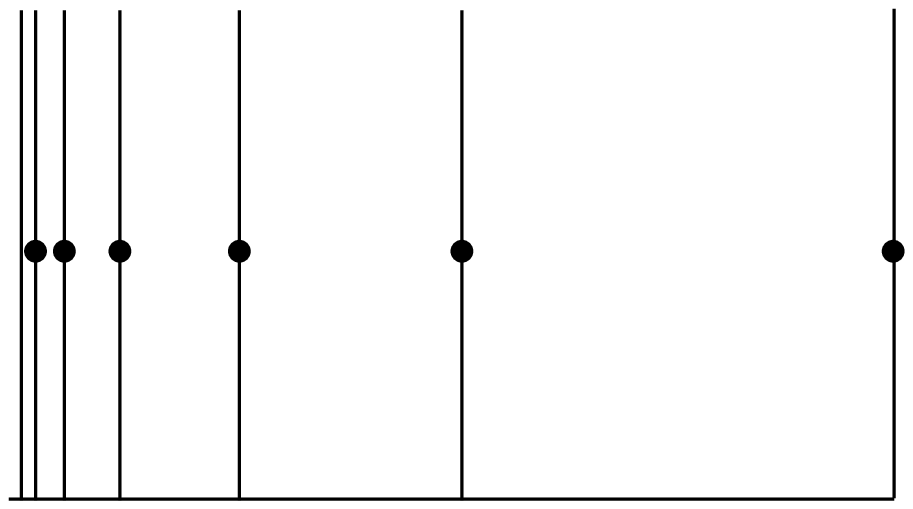}}
  \put(263.31,69.26){\fontsize{8.54}{10.24}\selectfont $x_1$}
  \put(139.14,69.63){\fontsize{8.54}{10.24}\selectfont $x_2$}
  \put(23.89,69.91){\fontsize{8.54}{10.24}\selectfont $x_n$}
\put(-5.89,-3.01){\fontsize{8.54}{10.24}\selectfont $(0,0)$}
  \end{picture}%
\fi \caption{\label{peine2} $M_G$ and $j(M_C)$ coincide as point
sets but not topologically.
\newline
$\hbox{$\quad$}$ Both, Gromov and Cauchy completions include
naturally $(0,0)$ as a boundary point, and coincide as point sets.
However, $\{x_n\}$ converges to $(0,0)$ only for Gromov topology.
Moreover, $M_C$ is not locally compact, and $j:M_C\rightarrow M_G$
is continuous but not a topological embedding.}
\end{figure}

For (b) and (d), modify slightly previous example by enlarging it
and redefining
$M$ as:
\begin{equation}\label{mprima} M:= M' \cup \{(x,1): \;
x\in(0,1]\}\end{equation}
where $M'$ is now given by the expression in (\ref{mmprima}) (see
Fig. \ref{peine} (A)). Observe that $j$ is a topological
embedding. In fact, 
from Prop. \ref{pemb} we only have to check the continuity of
$j^{-1}\mid_{j(M_C)}$. Observe that if $\{[-\dri_{x_{n}}]\}_{n}$
converges to $[-\dri_{x}]$ in $\lip_*$ for $x_n,x\in M_C$, Lemma
\ref{l4.5} ensures that $\lim_n
(\dri(\cdot,x)+\dri(x,x_n)-\dri(\cdot,x_n))=0$. If $\lim_n
\dri(x,x_n)\neq 0$, one can obtain the following contradiction:
there exist $y\in M, \epsilon>0$ and a subsequence $\{x_{n_k}\}$
such that $\dri(y,x)+\dri(x,x_{n_k})>\dri(y,x_{n_k})+\epsilon$. In
fact, considering $\{x_n\}$ as a sequence in $\R^2$, it will
converge, up to a subsequence, to some $x_0\in \R^2$ (with
$x_0\neq x$)
such that $x_0\in M_C$. Depending on the relative position of $x_0$ and
$x$, the claimed $y$ is chosen such that the strict triangle
inequality $d^R(y,x)+d^R(x,x_0)>d^R(y,x_0)$ holds.

Observe that, for $x_n=(n^{-1},1/2)$, the sequence
$\{[-\dri_{x_n}]\}$ converges in $\lip_*$ to the class
of the function 
${\rm max}\{-\dri_{(0,0)},-\dri_{(0,1)}\}$ (see Fig. \ref{peine} (B)),
which does not correspond with any element of $M_{C}$. This proves
(b) and, for (d), recall that $M_C$ is not locally compact, as any
neighborhood of $(0,0)$ contains a sequence of the form
$\{(n^{-1},\epsilon)\}$ for $\epsilon>0$ small enough, and this
sequence does not converge in $M_C$.
\begin{figure}
\centering
\ifpdf
  \setlength{\unitlength}{1bp}%
  \begin{picture}(298.64, 384.08)(0,0)
  \put(0,0){\includegraphics{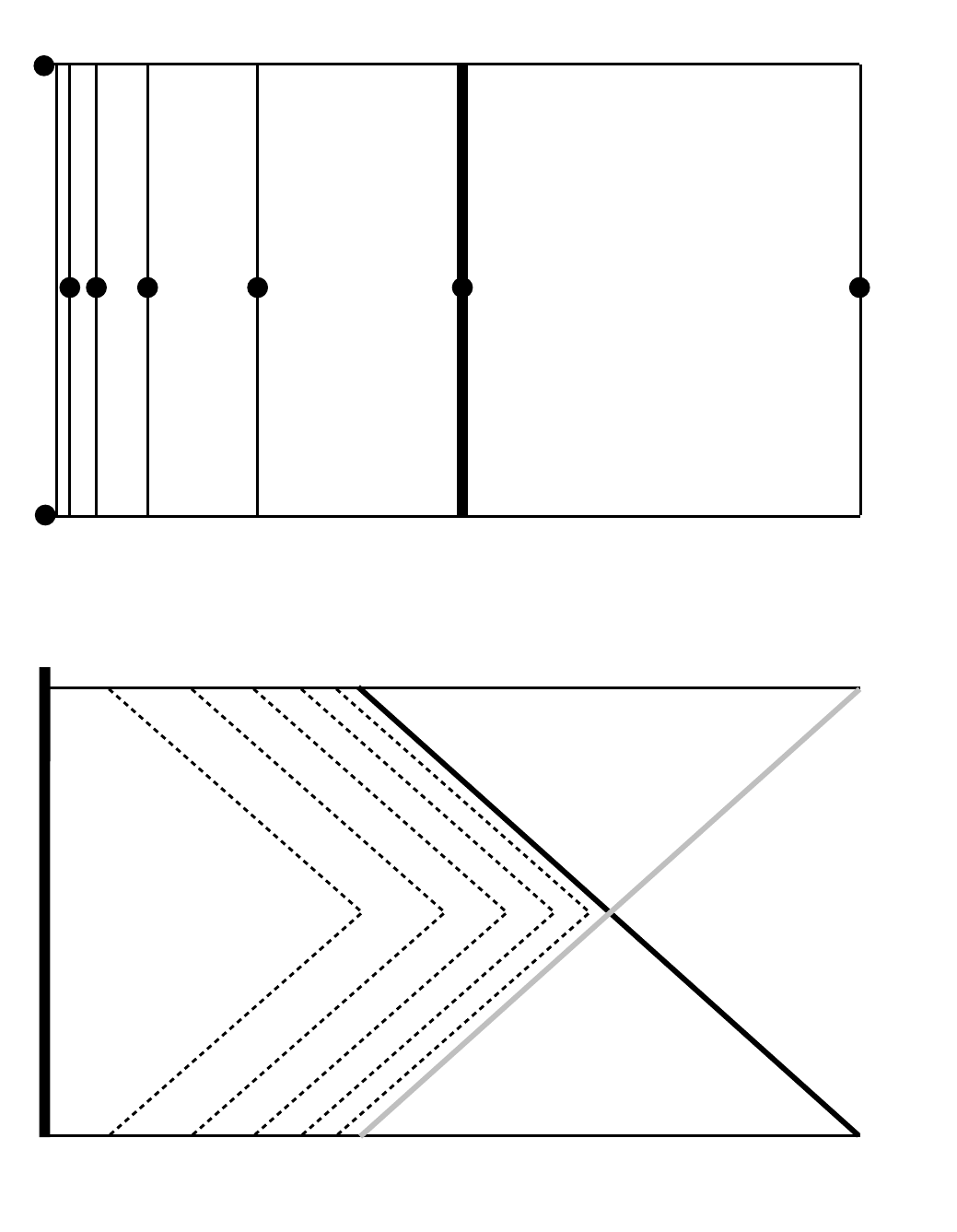}}
  \put(271.71,286.97){\fontsize{8.54}{10.24}\selectfont $x_1$}
  \put(147.54,287.34){\fontsize{8.54}{10.24}\selectfont $x_2$}
  \put(32.29,287.62){\fontsize{8.54}{10.24}\selectfont $x_n$}
  \put(16.91,173.86){\fontsize{8.54}{10.24}\selectfont $y$}
  \put(8.72,212.30){\fontsize{8.54}{10.24}\selectfont $(0,0)$}
  \put(5.67,371.74){\fontsize{8.54}{10.24}\selectfont $(0,1)$}
  \put(137.17,154.46){\fontsize{8.54}{10.24}\selectfont $-d(\cdot,(0,1))$}
  \put(137.50,40.99){\fontsize{8.54}{10.24}\selectfont $-d(\cdot,(0,0))$}
  \put(39.22,105.82){\fontsize{8.54}{10.24}\selectfont $-d(\cdot,x_n)+\frac{1}{2}$}
  \put(6.57,27.41){\fontsize{8.54}{10.24}\selectfont 0}
  \put(6.89,167.71){\fontsize{8.54}{10.24}\selectfont 1}
  \put(136.20,204.15){\fontsize{14.23}{17.07}\selectfont (A)}
  \put(135.08,8.73){\fontsize{14.23}{17.07}\selectfont (B)}
  \end{picture}%
\else
  \setlength{\unitlength}{1bp}%
  \begin{picture}(298.64, 384.08)(0,0)
  \put(0,0){\includegraphics{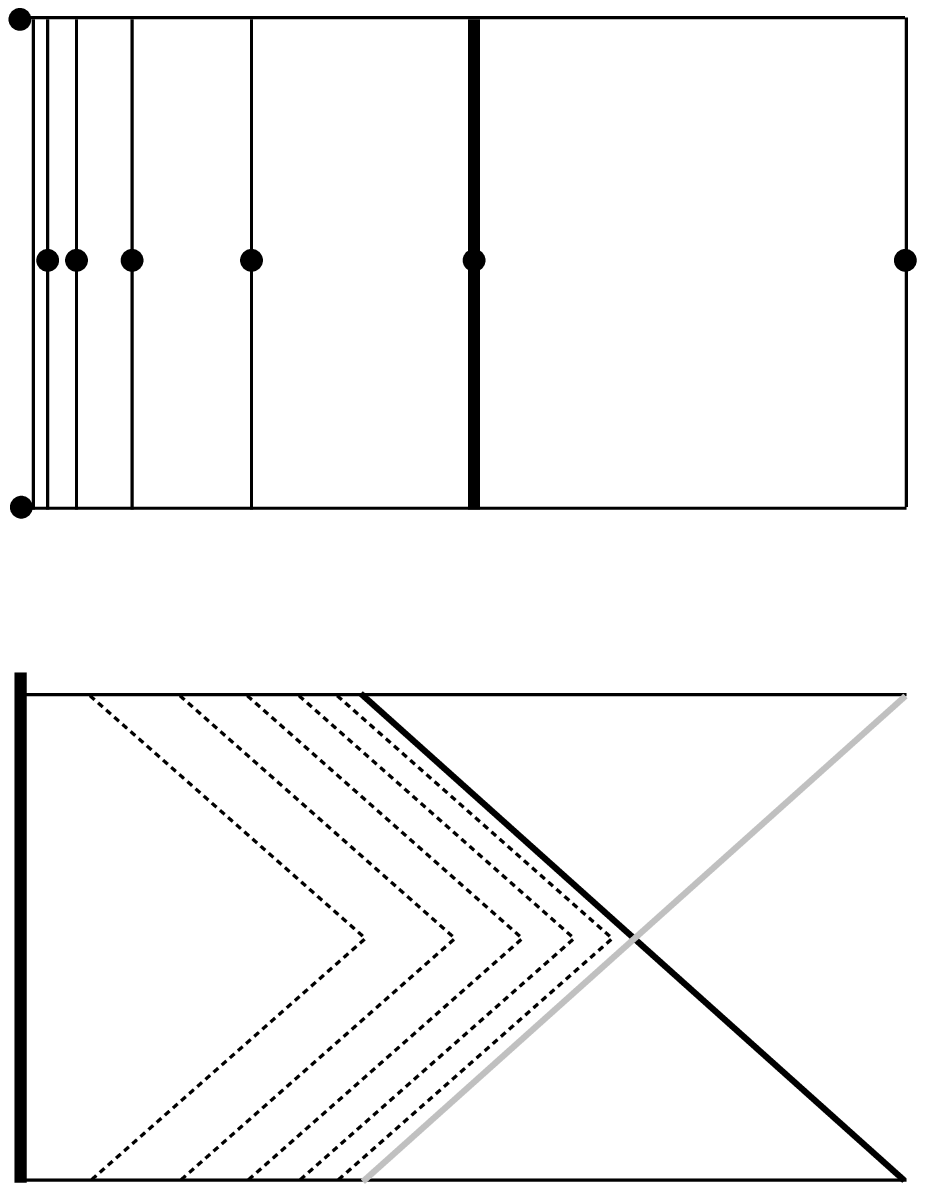}}
  \put(271.71,286.97){\fontsize{8.54}{10.24}\selectfont $x_1$}
  \put(147.54,287.34){\fontsize{8.54}{10.24}\selectfont $x_2$}
  \put(32.29,287.62){\fontsize{8.54}{10.24}\selectfont $x_n$}
  \put(16.91,173.86){\fontsize{8.54}{10.24}\selectfont $y$}
  \put(8.72,212.30){\fontsize{8.54}{10.24}\selectfont $(0,0)$}
  \put(5.67,371.74){\fontsize{8.54}{10.24}\selectfont $(0,1)$}
  \put(137.17,154.46){\fontsize{8.54}{10.24}\selectfont $-d(\cdot,(0,1))$}
  \put(137.50,40.99){\fontsize{8.54}{10.24}\selectfont $-d(\cdot,(0,0))$}
  \put(39.22,105.82){\fontsize{8.54}{10.24}\selectfont $-d(\cdot,x_n)+\frac{1}{2}$}
  \put(6.57,27.41){\fontsize{8.54}{10.24}\selectfont 0}
  \put(6.89,167.71){\fontsize{8.54}{10.24}\selectfont 1}
  \put(136.20,204.15){\fontsize{14.23}{17.07}\selectfont (A)}
  \put(135.08,8.73){\fontsize{14.23}{17.07}\selectfont (B)}
  \end{picture}%
\fi \caption{\label{peine} Bounded example with  $j(M_C)$
non-closed in $\lip_*$ (Example \ref{ex43} (b)) and non-Hausdorff
chr. topology  (Remark \ref{rjb}).
\newline
\indent $\hbox{$\quad$}$ (A) $\partial_C M$ includes only $(0,0)$
and $(0,1)$, but $\partial_G M$ also includes in a natural way
{\cambios the other  points $S=\{(0,y): 0<y<1\}$ of the segment between
that points. Then, $j:M_C\rightarrow M_G$ is an embedding but
$j(M_C)$ is not closed and $j(\partial_CM)\subsetneq
\partial_{CG}M$. Moreover, the point in $\partial_GM$ defined by the sequence $\{x_{n}\}_{n}$
(or any other point in $S\subset \partial_GM$)  cannot be reached
as the limit of any curve in $M$}.
\newline $\hbox{$\quad$}$ (B) Visualization of the chr. convergence
of $\{[-d(\cdot,x_n)]\}_n$ to $[-d(\cdot,(0,0))]$ and
$[-d(\cdot,(0,1))]$ in Figure (A). The bold vertical axis
represents the vertical segment through $x_2$ in (A). The dashed
lines correspond to the {\cambios graph of the functions
$-d(\cdot,x_n)+1/2$ on that vertical segment.} The limit is
${\rm max}\{-d(\cdot,(0,0)),-d(\cdot,(0,1))\}$.}
\end{figure}
}\end{example}

\bdefi Let $(M,\dri)$ be the metric space associated to a Riemannian manifold (or any space
under Convention \ref{cginebra}).

 The {\em Gromov completion} $M_G$ of $(M,\dri)$ is
the closure of $j(M)$ in $\lip_*$, and the {\em Gromov boundary}
$\partial_GM$ is 
$\partial_G M:=\overline{j(M)}\setminus j(M)$.

\edefi
The map obtained by restricting the co-domain of
the map $j$ defined in Prop. \ref{pemb} is also denoted  $$j:M_C\rightarrow M_G ,$$
and shows that the Cauchy boundary can be regarded as a point set
included in the Gromov one. \brema\label{r4.7} {\em Typically   in
Gromov's approach, one considers, instead of $\lip$, the space
${\cal C}(M)$ of all the continuous functions on $M$ endowed with
the topology of uniform convergence on compact sets.
However, in order to compute the closure of $j(M)$, this
difference in the ambient space becomes irrelevant. In fact, if
$[h]$ belongs to the closure of $j(M)$ in ${\cal C}(M)_*$ then any
$h\in [h]$ belongs to $\lip$. {\cambios Moreover, on this space, the topology of uniform convergence on compact sets and the topology of pointwise convergence coincide (see Remark \ref{rmetrizable}).}} \erema

\btheo \label{tgromov} The Gromov completion $M_G$ of any metric space $(M,\dri)$ associated to a
Riemannian metric (or space under Convention \ref{cginebra}) is a compact metrizable space.

The image of the manifold $j(M)\subset M_G$ is an open  dense
subset of $M_G$ and, thus,  the boundary $\partial_{G}M$ is also
compact.\etheo {\it Proof.} Take some $x_0\in M$. The map $i:
\lip_* \rightarrow \lip$ which maps each class $[f]$ to its
representative $f_0$ such that $f_0(x_0)=0$, is a topological
embedding. On the other hand, by Lemma \ref{l41}, $i(\lip_*)$ is
sequentially compact. Hence, the metrizable space $\lip_*$ is also
sequentially compact, and thus compact (see, for example,
\cite[Theor. 28.2]{Mu}). So, the first assertion follows from the
fact that $M_{G}$ is closed in $\lip_*$. Finally, the density of
$j(M)$ follows trivially from the construction of $M_G$,  and its
openness 
from the local compactness of $M$ as in
Proposition \ref{p329} (see also Corollary \ref{cembbusf} below).
\cvd

\smallskip

As $j|_M$ is an embedding, 
the continuous map $j$ can be dropped for $M$, and we will  write
simply $M \subset M_G$. However, $j$ will be maintained for $M_C$
and $\partial_CM$, as $j(M_C), j(\partial_CM)$ will be endowed
with the restriction of the Gromov topology of $M_G$, which, in
general, does not agree with the Cauchy completion topology
(Example \ref{ex43}, property (c)). In the Gromov boundary
$\partial_GM$, one can distinguish two disjoint parts: the {\em
Cauchy-Gromov boundary} $\partial_{CG}M$, which consists of all
the points which can be obtained as the limit of some sequence
of $M$ included in some bounded subset, and the {\em proper Gromov
boundary} $\bpgro =\partial_GM\setminus
\partial_{CG}M$, which
corresponds with the ``boundary at infinity''. The Cauchy boundary
$\partial_CM$ is naturally included in $\partial_{CG}M$ and, as
Example \ref{ex43} (property (b)) shows, the inclusion can be
strict, i.e. a {\em residual Gromov} boundary
$\partial_{RG}M=\partial_{CG}M\setminus j(\partial_CM)$ may
appear. Summing up:
$$M_G=M \cup \bgro \quad \quad
\bgro= \partial_{CG}M \cup \bpgro \quad \quad \partial_{CG}M =
j(\partial_{C}M) \cup
\partial_{RG}M.
$$
The topological subspace $M_{CG}:=M\cup
\partial_{CG}M$ of $M_G$ will be called the {\em Cauchy-Gromov
completion} and its properties are summarized as
follows.
\begin{corollary}\label{ccg} The Cauchy-Gromov completion $M_{CG}$
satisfies:
\begin{itemize}
\item[(1)] (Openness of $M$). The boundary $\partial_{CG}M$ is
closed in $M_{CG}$.

\item[(2)] (Heine-Borel). The closure in $M_{CG}$ of any bounded
set in $M$ is compact.

\item[(3)] (Consistency with the Cauchy boundary). If
{\cambios $M_C$ is locally compact}, then $j: M_C \rightarrow M_{G}$
is an embedding and $M_C\equiv M_{CG}$, i.e.
$\partial_{RG}M=\emptyset$.
\end{itemize}
\end{corollary}
{\em Proof.} The first assertion is straightforward from Theorem
\ref{tgromov}. For the second assertion recall that, from the
compactness of $M_G$, the closure of any set in $M_G$ is compact
and, from the definition of $M_{CG}$, the closure of any bounded
set belongs to $M_{CG}$. The first statement in (3) follows from
Remark \ref{mclc} and, for the last statement, it suffices to show
that, if $M_C$ is locally compact, then  the closed balls in
$M_C$, $\overline{B(x,r;M_C)}$, $x\in M_C$, $r>0$, are compact. In
fact, in this case any sequence defining a point of the
Cauchy-Gromov boundary will have a convergent subsequence in
$M_C$, and the result follows from the continuity of $j$ and the
Hausdorffness of $M_G$. So, in order to prove that
$\overline{B(x,r;M_C)}$ is compact, consider a sequence
$\{x_n\}\subset\overline{B(x,r;M_C)}$.
Define a sequence of curves $\{\gamma_n\}$,
$\gamma_n:[0,1]\rightarrow M$ joining $x$ and $x_n$, such that:
each $\gamma_n$ restricted to $(0,1)$ is smooth and contained in
$M$, it is parametrized with constant speed and  ${\rm
length}(\gamma_n)<2r$ (see Remark \ref{r1}). Clearly, the sequence
of curves $\{\gamma_n\}_n$ is
equicontinuous\footnote{\label{footequicon}For posterior
referencing, note that, if $t_1<t_2$ then
$\dri(\gamma_n(t_1),\gamma_n(t_2))\leq {\rm
length}(\gamma_n)|_{[t_1,t_2]}=(t_2-t_1){\rm length}(
\gamma_n) <(t_2-t_1)2r.$  This inequality, plus the symmetry of
$\dri$, implies the equicontinuity.}. Taking into account that $x$
is an accumulation point of the sequence $\{\gamma_n(0)\}$, and
$M_C$ is both, locally compact and complete, classical Arzela's
Theorem (see for example, Theorems \ref{tafe}, \ref{taf} below)
implies the existence of a curve $\gamma_\infty$ which is the
pointwise limit of a subsequence of $\{\gamma_n\}$. So, the
required limit is $\gamma_\infty(1)$. \cvd

\subsection{Busemann completion as a point set}\label{s42}


Next, we are going to focus our attention on certain subset of
$M_{G}$, which provides another compactification of $(M,d^R)$ when
is endowed with a different topology.

Let $\cur$ be the space of all the piecewise smooth curves
$c:[\alpha,\Omega) \rightarrow M$,
$-\infty<\alpha<\Omega\leq\infty$, with
$|\dot c|^{2}=g_R(\dot{c},\dot{c})< 1$ 
and consider the
associated function: 
\begin{equation}\label{ebus}
b_c(x) = \lim_{s\nearrow \Omega} (s-\dri(x,c(s))) \in \R \cup
\{+\infty\},\qquad x\in M. \end{equation} The following result
shows that $b_c$ is well-defined:
\begin{lemma}\label{l1} If $c\in \cur$ then the map
$s\mapsto s-\dri(x,c(s))$ is increasing for any $x\in M$.
\end{lemma}
{\it Proof.} 
Just notice:
$$\dri(x,c(s_{2}))-\dri(x,c(s_{1}))\leq \dri(c(s_{1}),c(s_{2}))\leq
{\rm length}(c\mid_{[s_{1},s_{2}]})< s_{2}-s_{1}.$$
\cvd


\bprop \label{p1} Let $b_c$ be 
associated to some $c\in \cur$.
\begin{itemize}
\item[(1)] If $b_c$ reaches the value $+\infty$ at  some $x\in M$
then $b_c \equiv \infty$; otherwise, $b_c\in \lip$ and its
corresponding class $[b_{c}]$ belongs to $M_G$. \item[(2)] If
$\Omega<\infty$ then: (i) there exists $\overline{x}\in M_C$ such
that $\lim_{s\nearrow \Omega}c(s)= \overline{x}$ in the metric
topology of $M_C$; (ii) $b_c(\cdot)= \Omega -
\dri_{\overline{x}}(\cdot)$; (iii) $j(\overline x)= [b_c] \in
j(M_C)(\subset M_G$) and $\lim_{s\nearrow \Omega}j(c(s))=[b_c]$.

Conversely, for any $\bar x\in M_C$, there exists some $c\in C(M)$
with $\Omega<\infty$ such that  $j(\bar x)=[b_c]$. In particular,
for $x\in M$, one has $j(x)=[b_{c_x}]$, where $c_x: [-1,0)
\rightarrow M$, $s\mapsto x$.

 \item[(3)] If $\Omega = \infty$ and $b_c\not\equiv
\infty$, then $[b_c] \in
\partial_{\cal G} M$.
\end{itemize}
\eprop {\it Proof.} It is a consequence of previous definitions
(see also \cite[Sect. 2]{H}):

(1) 
Both assertions follow from 
$\dri(x,c(s))-\dri(y,c(s))\leq \dri(x,y)$ for all $s$.

(2) As the length of $c$ is finite, some endpoint $\bar x\in M_C$
is determined (recall Remark \ref{r1}). This  proves (i), and the
pointwise convergence of $\dri_{c(s)}$ to $\dri_{\bar x}$ proves (ii).
The remainder follows from the definitions above.

(3) As $b_c<\infty$ and $\Omega=\infty$, the curve $c$ must leave
any bounded subset of $M$, and so, its class cannot belong to
$M_{CG}$. \cvd

\begin{definition}\label{db} A {\em Busemann function} $b:M\rightarrow
(-\infty,\infty]$ is a function which can be written as $b=b_c$
for some $c\in \cur$ as in (\ref{ebus}). The set of all the
finite-valued Busemann functions will be denoted $B(M)(\subset
\lip)$.

A {\em properly Busemann function} $b$ is a finite-valued one which
 is written as $b=b_c$ for some $c\in C(M)$ with $\Omega
=\infty$. The set of all the properly Busemann functions will be
denoted ${\cal B}(M)(\subset B(M)$).

As a point set, the {\em Busemann completion} $\bus$ of $(M,\dri)$
is defined as the subset $\bus= B(M)/\R$ of the Gromov completion
$M_G$. The {\em Busemann boundary} is then $\bbus= \bus\backslash
j(M)$ and the {\em properly Busemann boundary} $\bpbus= {\cal
B}(M)/\R$ ($=\bus\backslash j(M_C)\subset
\partial_{\cal G}M$).
\end{definition}

\begin{remark}\label{remm1}{\em About the notion of Busemann function, the
following comments are in order:

(1)  There is no restriction if, in the definition of $B(M)$, one
considers curves in $C(M)$ with unit velocity.  In fact, if
$\overline{c}$ is the arc-reparametrization of a curve $c\in C(M)$
with $b_c\in B(M)$, then $b_{c}=k_0+b_{\overline{c}}$, where
$k_0=\lim_t (t-$length$|_{[\alpha,t]}c)$ (which is finite as
$k_0<\lim_t (t-\dri(c(\alpha),c(t)))=b_c(c(\alpha))<\infty$). So,
$\tilde{c}(t)=\overline{c}(t-k_0)$ satisfies $b_{\tilde c}=b_c$,
as required. Therefore, we choose $|\dot{c}|^2< 1$ only for
convenience.

 (2) Moreover, according to previous paragraph, one can also
consider that the curves in $C(M)$ have velocity lower or equal than one.

(3) In any metric space $(X,d)$, one can say that $b:X\rightarrow
(-\infty,\infty]$ is a {\em Busemann function} when there exists a
sequence $\sigma=\{x_n\}_{n}\subset X$ such that $b=b_\sigma$,
being $b_\sigma(x)=\lim_n\left(\sum_{k=1}^n
d(x_{k-1},x_k)-d(x_n,x)\right)$, $x\in X$. This avoids the use of
length spaces.
However, we will need other properties of length spaces, and this
general definition will not be considered in this paper.}
\end{remark}
Summing up,   Definition \ref{db} and Proposition \ref{p1} yield
the disjoint unions
$$M_B= j(M_C) \cup \bpbus , \;  \quad \quad \bbus= j(\partial_CM) \cup
\bpbus$$ and the relations
$$\partial_{\cal B}M \subset
\partial_{\cal G}M , \; \quad \quad \partial_BM \cap \partial_{RG}M=\emptyset
$$
complemented with $\bbus \subset
\partial_GM$, $j(\partial_CM)\subset \partial_{CG}M$,
 where all the inclusions may be strict (see
Figure \ref{escalera2}).

\begin{figure}
\centering
\ifpdf
  \setlength{\unitlength}{1bp}%
  \begin{picture}(262.52, 313.19)(0,0)
  \put(0,0){\includegraphics{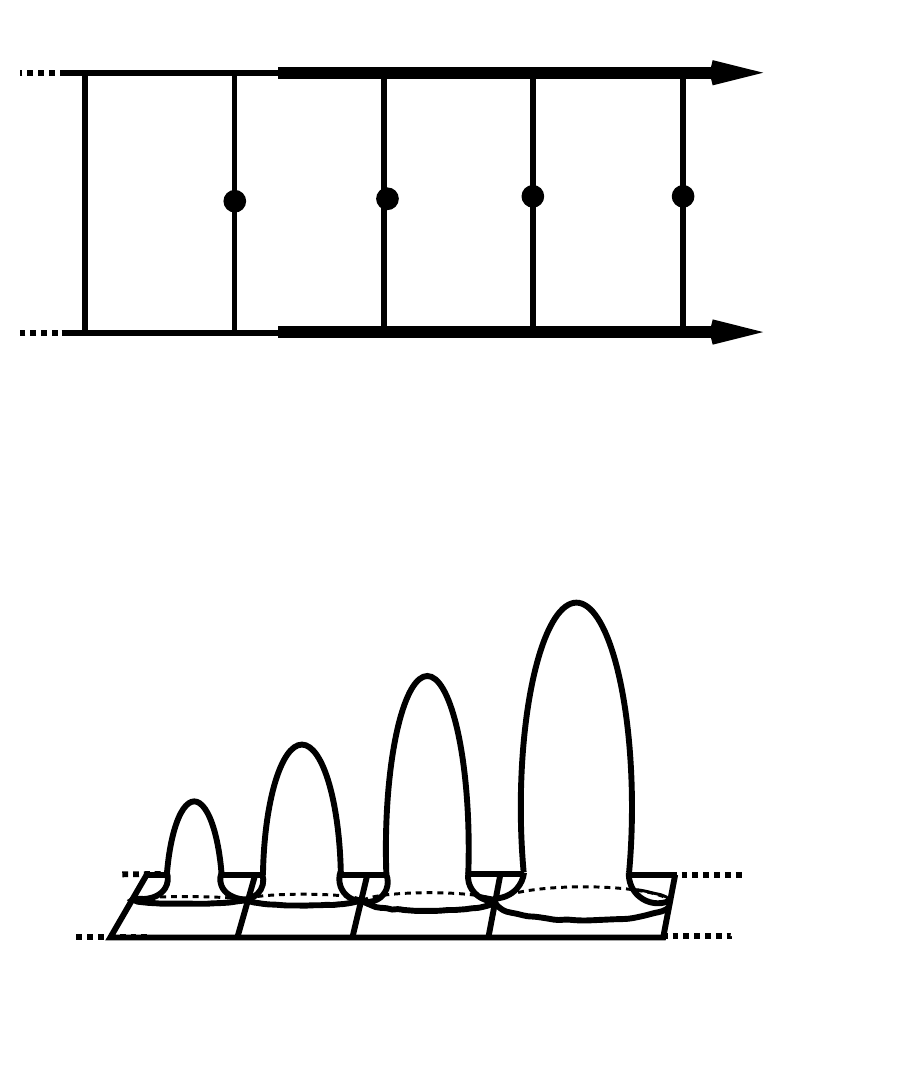}}
  \put(108.96,182.66){\fontsize{14.23}{17.07}\selectfont (A)}
  \put(108.92,8.73){\fontsize{14.23}{17.07}\selectfont (B)}
  \put(154.86,206.51){\fontsize{11.38}{13.66}\selectfont $c_1$}
  \put(155.56,298.63){\fontsize{11.38}{13.66}\selectfont $c_2$}
  \put(73.72,255.36){\fontsize{11.38}{13.66}\selectfont $x_1$}
  \put(116.99,255.36){\fontsize{11.38}{13.66}\selectfont $x_2$}
  \put(202.13,255.36){\fontsize{11.38}{13.66}\selectfont $x_n$}
  \end{picture}%
\else
  \setlength{\unitlength}{1bp}%
  \begin{picture}(262.52, 313.19)(0,0)
  \put(0,0){\includegraphics{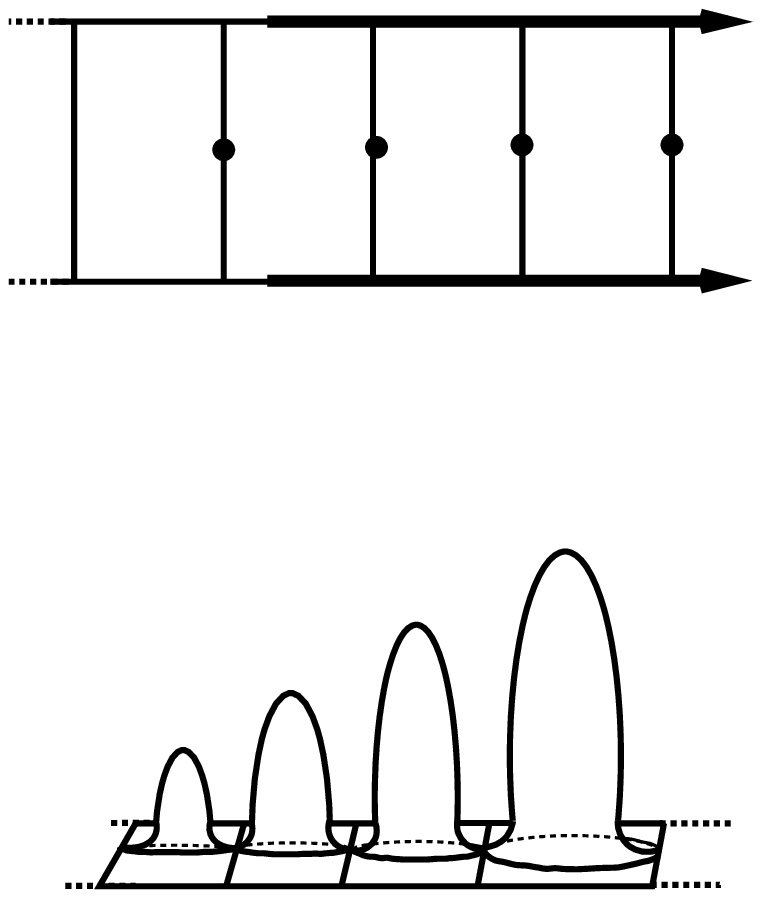}}
  \put(108.96,182.66){\fontsize{14.23}{17.07}\selectfont (A)}
  \put(108.92,8.73){\fontsize{14.23}{17.07}\selectfont (B)}
  \put(154.86,206.51){\fontsize{11.38}{13.66}\selectfont $c_1$}
  \put(155.56,298.63){\fontsize{11.38}{13.66}\selectfont $c_2$}
  \put(73.72,255.36){\fontsize{11.38}{13.66}\selectfont $x_1$}
  \put(116.99,255.36){\fontsize{11.38}{13.66}\selectfont $x_2$}
  \put(202.13,255.36){\fontsize{11.38}{13.66}\selectfont $x_n$}
  \end{picture}%
\fi \caption{\label{escalera2}(A) A simple modification of Figure
\ref{peine} which satisfies $\partial_{\cal B} M\subsetneq
\partial_{\cal G} M$. The Busemann functions for the
horizontal curves $c_1 , c_2$ yield two boundary points $z_1, z_2$
in $\partial_{\cal B}M$. The sequence $\{x_n\}$  converges to
both, $z_1$ and $z_2$, in the chr. topology, and to a different
point $z$ of $\partial_{\cal G}M$ in the pointwise topology.
\newline $\hbox{$\quad$}$ (B) Filling the holes in (A) by using diverging (or just big enough) ``bubbles'',  a complete
simply-connected space which behaves qualitatively as the previous
one ($\partial_{\cal B} M\subsetneq
\partial_{\cal G} M$) is constructed. Positive sectional curvature is required in some regions (see Remark \ref{rgrape}).
}
\end{figure}
\begin{remark}\label{rgrape} {\em In some sense, $\bpbus$ represents the set of
possible ``directions towards infinity''. In a Hadamard manifold,
$\partial_CM=\emptyset$ and $\bpbus=\partial_GM$ {\cambios (see
Corollary \ref{new} below);} moreover, in this case one can
construct $\bpbus$ by using only rays, instead of all $\cur$ (in
fact, the classical construction by Gromov shows that $M_G$, and
thus $M_B$, can be obtained by using only rays). However, it is
not difficult to find even a complete, 1-connected Riemannian
manifold where $\bpbus \varsubsetneq
\partial_GM$. This is suggested by the example in
Figure \ref{escalera2} (B);
 a different example is the universal covering of the ``grapefruit
on a stick'', studied in detail in \cite[Sect. 2.1]{FH}.
}\end{remark}




If $\bbus\neq
\partial_GM$ then $M_B$ is not compact with the pointwise
convergence topology, and the name {\cambios ``Busemann completion
or compactification'' may seem} misleading. However, the
natural topology for this completion, the {\em chronological
topology} (which is suggested by the causal boundary), makes
$\bus$ and $\bbus$ {\cambios sequentially compact} (but at the price of
being non-Hausdorff). Such a topology is studied below, in the
more general (non-necessarily symmetric) Finslerian case,
Subsection \ref{sectbusf}.

\section{Finslerian completions}\label{s5}

Our aim in this section is to study previous compactifications for
a (connected) smooth manifold $M$ endowed with the generalized
distance $d$ associated to a Finsler metric $F$ on $M$. As in the
Riemannian case, this can be done in somewhat more general spaces.

\begin{convention}\label{cginebra2}{\em Along the present section, $(M,d)$ will denote
any generalized metric space satisfying: $M$ is locally compact
and contains a countable dense set, and $d$ is derived from a
length space (except if specified otherwise). Recall that {\cambios
most of the notions of the previously studied  symmetric case can
be translated to the present non-symmetric case in a}
straightforward way; in particular, for a smooth curve, the
Finslerian norm $F(\dot{c})$ must be interpreted as the dilatation
(see for example \cite[Appendix]{CJS}, and also Convention
\ref{cginebra}). The extended quasi-distance on the Cauchy
completions $M_C^\pm$ is denoted by $d_Q$.
}
\end{convention}

In order to translate previous constructions to a generalized
metric space $(M,d)$, first we need to generalize the notion of
Lipschitz function.
\begin{definition}\label{dflips}
 A function
$f:M\rightarrow \R$ on a generalized metric space $(M,d)$ is:
\begin{itemize}
\item {\em 1-Lipschitz} or, simply, {\em Lipschitz}, if
$f(y)-f(x)\leq d(x,y)$ for all $x,y\in M$.
\item {\em Max-Lipschitz} if
$|f(x)-f(y)|\leq max\{d(x,y),d(y,x)\}$ for all $x,y\in M$.
\end{itemize}
The sets of all the Lipschitz and Max-Lipschitz functions on
$(M,d)$ will be denoted $\lipd$ and $\lipdM$, resp.
\end{definition}

\begin{remark}\label{rrr1}{\em (1) In principle, the space $\lipdM$ uses $d$ in a symmetric way and
seems closer to the one used in the Riemannian case. However, this
space (as well as the notion of Lipschitz function for the
symmetrized distance $d^s$, which is defined in the obvious way)
has the drawback that the involved distance does not come from a
length space (see \cite{CJS}, especially its Appendix, for a
discussion and examples in the class of Randers metrics). This and
other reasons of compatibility with the causal boundary justify
that, in what follows, the space $\lipd$ will be preferred.
However, the space $\lipdM$ will be useful for technical
purposes.

%
%
%

(2) The reasonings in Remark \ref{rmetrizable} are applicable to $\lipdM$. So, this space is metrizable and the pointwise, compact-open and uniform
topologies coincide on it. Moreover, Lemma \ref{l41} remains true if we replace $\lip$ by
$\lipdM$ and, then, a first Gromov-type compactification as in
Theorem \ref{tgromov} would be defined
 by using the space $\lipdM_*=\lipdM/\R$.} \end{remark} Now, for $\lipd$:
\begin{proposition}\label{llll1}
(1) $\lipd$ is a closed subset of $\lipdM$.

(2) Let $\{f_n\}$ in $\lipd$ with $\{f_n(x_0)\}$ converging to
$k_0\in [-\infty,\infty]$ for some $x_0\in M$. If $k_0 \in \R$
then it admits a pointwise convergent subsequence; otherwise, it
converges uniformly on compact sets to $\pm \infty$.

(3) The uniform, compact-open and pointwise topologies are equivalent on $\lipd$.
\end{proposition}
{\it Proof.} For (1), let $f$ be the pointwise limit of a sequence
$\{f_n\}\subset \lipd$. Note that, if $f(y)>f(x)$, then
$f_n(y)>f_n(x)$ for $n$ big enough, thus $f(y)-f(x)=\lim_n
(f_n(y)-f_n(x))<d(x,y)$.

The proof of (2) follows from the facts that Lemma \ref{l41} is
true for $\lipdM$, and $\lipd$ is a closed subset of $\lipdM$.
Item (3) is proved analogously. \cvd

\subsection{Gromov completions for the non-symmetric case}\label{secgro}

Endow $\lipd$ with the pointwise convergence topology and
$M_C^\pm$ with its natural quasi-distance $d_Q$ and topologies
(Convention \ref{conv}). Clearly, the maps
\begin{equation}\label{equuu1}\begin{array}{l}\hat{j}^{+}: M_C^+ \rightarrow \lipd, \quad x\mapsto - d^{+}_{x}, \quad \hbox{ where} \;
d^{+}_{x}=d_Q(\cdot,x)\\ \hat{j}^{-}: M_C^- \rightarrow \lipd,
\quad x\mapsto +d^{-}_{x}, \quad \hbox{ where} \;
d^{-}_{x}=d_Q(x,\cdot)\end{array}\end{equation} define two
injective  functions. Up to a sign, the second one can be regarded
as the first one for $d_{Q}^{{\rm rev}}$.
In what follows, we will consider {\cambios only the map} $\hat{j}^+$
(being the case for $\hat{j}^-$ completely analogous). The map
$\hat{j}^{+}$ is not automatically continuous, as happened in the
Riemannian case (Prop. \ref{p41}). However:
\begin{lemma}\label{lcf}
The map $\hat{j}^+:M_C^+\rightarrow \lipd$ is continuous if and only if $d_Q$ satisfies
the following property:
\begin{itemize} \item[(a4')] For any sequence $\{x_n\}\subset
M_C^+$ and $x\in M_C^+$ such that $\lim_n d_Q(x_n,x)=0$,
necessarily $\lim_n d_Q(x,x_n)=0$, i.e.
the topology generated by the backward balls is finer than the one generated by the forward balls.
\end{itemize}
Moreover, if $d_Q$ is a generalized distance, then $\hat{j}^+$ is
an embedding.
\end{lemma}
{\it Proof.} For the implication to the right, suppose that
$\hat{j}^+$ is continuous and consider a sequence
$\{x_{n}\}\subset M^{+}_{C}$ and a point $x$ in $M_C^+$ such that
$\lim_n d_{Q}(x_n,x)=0$. By the continuity of $\hat{j}^+$,
necessarily $\{\hat{j}^+(x_n)\}\rightarrow \hat{j}^+(x)$, i.e.
$\lim_n d_{Q}(\cdot,x_n)=d_{Q}(\cdot,x)$. Then, the result follows
by evaluating the last expression in $x$.

For the implication to the left, suppose that condition (a4') is
satisfied and take $\{x_{n}\}\subset M_C^+$, $x\in M_C^+$ with
$\{x_n\}\rightarrow x$, i.e. $\lim_n d_{Q}(x_n,x)=0$. Then:
\[
\begin{array}{rl}d_{Q}(\cdot,x)\leq & \lim_n\left(
d_{Q}(\cdot,x_n)+d_{Q}(x_n,x)\right)\\ \leq & \lim_n\left(
d_{Q}(\cdot,x)+d_{Q}(x,x_n)+d_{Q}(x_n,x)\right).\end{array}
\]
Taking into account that $\lim_n d_{Q}(x_n,x)=\lim_n
d_{Q}(x,x_n)=0$, the intermediate equality holds, i.e.
$d_{Q}(\cdot,x)=\lim_n d_{Q}(\cdot,x_n)$, and so,
$\hat{j}^+(x_n)\rightarrow \hat{j}^+(x)$.

For the last assertion, recall that any sequence
$\{d_Q(\cdot,x_n)\}$ converging pointwise to $d_Q(\cdot,x)$
satisfies that $\lim_n d_Q(x,x_n)=d_Q(x,x)=0$. So, from (a4) in
Definition \ref{gms}, $\lim_n d_Q(x_n,x)=0$, i.e., $\{x_n\}$
converges to $x$. \cvd

\begin{remark}\label{rcf} {\em Observe that condition (a4) in Definition \ref{gms} is more restrictive than (a4') here (see also Remark
\ref{ra4}). So, if $d_Q$ is a generalized distance then
$\hat{j}^+$ is continuous (and analogously for $\hat{j}^{-}$).}
\end{remark}
From Lemma \ref{lcf}, and reasoning as
in Proposition \ref{pemb} {\cambios and Remark \ref{mclc}} one deduces:
\begin{proposition}\label{pembf}
Consider the natural map
$$j^+ :M_C^+ \rightarrow \lipd _*,\qquad x\mapsto [-d^{+}_{x}].$$
Then, $j^+|_{M}$ is a topological embedding. Moreover  $j^+$ is
continuous if and only if $d_Q$ satisfies (a4'). If, in addition,
$d_Q$ is a generalized distance and $M_C^+$ is locally compact
then $j^+$ is also an embedding.
\end{proposition}
\begin{definition}\label{def1}
Let $(M,d)$ be the generalized metric space associated to a
Finsler manifold (or any space under Convention \ref{cginebra2}).
The {\em forward (resp. backward) Gromov completion} $\gromp$
(resp. $\gromm$) of $(M,d)$ is the closure of $j^+(M)$ (resp.
$j^-(M)$) in $\lipd _*$.

Then, the {\em forward (resp. backward) Gromov boundary}
$\partial_G^+ M$ (resp. $\partial_{G}^{-}M$) of $M$ is defined as
$\partial_G^+ M=M_{G}^{+}\setminus j^{+}(M)$ (resp. $\partial_G^-
M=M_{G}^{-}\setminus j^{-}(M)$).
\end{definition}
\begin{theorem}\label{tgromovf}
$M_G^+$ and $\partial_G^+ M$ are compact metrizable
topological spaces.

The image of the manifold $j^+(M)\subset M_G$ is an open dense
subset of $M_G^+$ and, thus, the boundary $\partial_{G}^+M$ is also compact.
\end{theorem}
{\it Proof.} Notice that the map $i:{\lipd}_*\rightarrow \lipdM$,
which sends each class $[f]$ to its representative $f_0$ such that
$f_0(x_0)=0$ for some fix $x_0\in M$, is an embedding. So, the
sequential compactness follows Prop. \ref{llll1} (2), and the
result follows taking into account that $\lipdM$ is
metrizable (Remark \ref{rrr1} (2)).

As in Theorem \ref{tgromov}, the second assertion follows from
standard arguments. \cvd

\smallskip

In analogy with the Riemannian case, the Gromov boundary (and the
Gromov completion analogously) can be divided in two disjoint
parts. One of them is the {\em Cauchy-Gromov boundary},
$\partial_{CG}^+ M$, which contains all the points which are limit
of some sequence of $M$ included in some bounded subset, i.e.
which are contained in {\cambios some} (forward) ball
$B^+(x,r)$ with $x\in M$ and $r>0$. The other part is the {\em
proper Gromov boundary}, $\partial_{\cal G}^+ M=\partial_G^+
M\setminus
\partial_{CG}^+ M$. The Cauchy-Gromov boundary contains the Cauchy
boundary as a point set, so we define the {\em residual Gromov
boundary}, $\partial_{RG}^+ M:=
\partial_{CG}^+ M\setminus j^{+}(\partial_{C}^+ M)$. Consequently, the {\em Cauchy-Gromov
completion}
is defined as $M_{CG}^+ = M\cup \partial_{CG}^+ M$ and $M_{G}^+=
M_{CG}^+ \cup \partial_{{\cal G}}^+ M$.

In order to understand the structure of $M^{+}_{CG}$, a Finslerian
version of Corollary \ref{ccg} is required. Recall however that,
in the proof of that corollary, Arzela's Theorem was applied to a
sequence of curves. One can prove directly a version of Arzela's
Theorem for generalized distances. But this would be insufficient
for our purposes here, as we will also need to include a
generalization of the notion of equicontinuity\footnote{As an
alternative to Arzela's Theorem, the extended version of
Hopf-Rinow Theorem for spaces under Convention \ref{cginebra} (see
\cite{Gr}) could have been used. In order to follow such an
approach here, a Finslerian version of this extended Hopf-Rinow
Theorem would be required. Our approach below also allows to prove
this version, and makes apparent the difficulties associated to
non-symmetry.}.

\begin{definition}\label{doreq}
Let $(M,d)$ be a generalized metric space associated to a Finsler
manifold. A sequence of functions $\{f_n\}$, $f_n:[a,b]\rightarrow
M$, is {\em oriented-equicontinuous} if, for any $\epsilon>0$,
there exists $\delta>0$ such that, for any $t_1,t_2\in [a,b]$ with
$0\leq t_2-t_1<\delta$, then $d(f_n(t_1),f_n(t_2))<\epsilon$ for
all $n$. If the requirement $0\leq t_2-t_1<\delta$ is weakened
into $0\leq |t_2-t_1|<\delta$, then $\{f_n\}$ is {\em
equicontinuous}.
\end{definition}
\begin{remark}{\em (1) If the sequence $\{f_n\}$ is oriented-equicontinuous, then each $f_n$ is continuous. In fact, fix $c\in
[a,b]$ and consider a sequence $\{t_m\}\rightarrow c$. With no
loss of generality, consider the cases: (i) $t_m\leq c$ and (ii)
$t_m\geq c$. For (i), as $f_n$ belongs to an
oriented-equicontinuous sequence, the definition yields directly
$\{d(f_{n}(t_m),f_{n}(c))\}_{m}\rightarrow 0$. For (ii), the
oriented-equicontinuity ensures that
$\{d(f_n(c),f_n(t_m))\}_{m}\rightarrow 0$ and then, by using that
$d$ is a generalized distance, we deduce
$\{d(f_{n}(t_m),f_{n}(c))\}_{m}\rightarrow 0$  from condition (a4)
in
Definition \ref{gms}, as required. 

(2) Oriented-equicontinuity does not imply equicontinuity. In
fact, one can easily construct an example of this situation by
endowing the space in Figure \ref{peine2} with a Finslerian metric
such that the vertical segments have unitary lengths when they are
{\cambios parameterized} by a curve from down to up, and diverging
lengths when they are {\cambios parameterized} from up to down.
}
\end{remark}


The following result is the (local) version of Arzela's Theorem
for equicontinuity in generalized metric spaces, which will be
refined later for oriented-equicontinuity (Theorem \ref{taf}).
\begin{theorem}\label{tafe}
Let $(M,d)$ be a generalized metric space and $\{f_n\}$,
$f_n:[a,b]\rightarrow M$, a sequence of equicontinuous functions.
Assume that $x\in M$ is an accumulation point of $\{f_n(c)\}$,
with $c\in [a,b]$, admitting a compact neighborhood $K$. Then,
there exist $r>0$, a subsequence $\{f_{n_k}\}\subset \{f_n\}$ and
a continuous function $f:[c-r,c+r]\cap [a,b]\rightarrow M$ such
that $\{f_{n_k}(t)\}_k\rightarrow f(t)$ for each $t\in
[c-r,c+r]\cap [a,b]$.
\end{theorem}
{\it Proof.} It is carried out by standard arguments, which we
include for the sake of completeness\footnote{Recall that here not
all the assumptions in Convention \ref{cginebra2} are required.}.
By the equicontinuity, and the facts that $x$ is an accumulation
point and $d$ a generalized distance, there exists $r>0$ such that
$f_n([c-r,c+r]\cap [a,b])\subset K$ for all $n$ big enough. Let
$A=\{t_m\}_{m=0}^{\infty}\subset [c-r,c+r]\cap [a,b]$ be a
countable dense subset, and consider the sequence of sequences
$\{f_n(t_m)\}_{n}\subset K$. A standard diagonal argument (see the
proof of Lemma \ref{l41}) ensures the existence of a subsequence
$\{f_{n_k}\}_{k}$ such that, for each $m$, $\{f_{n_k}(t_m)\}_{k}$
converges to some point, denoted $f(t_m)$.

Now, we are in conditions to prove that $\{f_{n_k}(t)\}$ is
convergent for each $t\in [c-r,c+r]\cap [a,b]$. In fact, as the
sequence $\{f_{n_k}(t)\}_{k}$ is contained in the compact $K$, it
suffices to show that it is a Cauchy sequence. Take $\epsilon>0$,
from the equicontinuity of $\{f_{n_k}\}$ and the density of $A$,
there exists $t_m\in A$ such that
$d(f_{n_k}(t_m),f_{n_k}(t)),d(f_{n_k}(t),f_{n_k}(t_m))<\epsilon/3$
for all $k$. On the other hand, as $\{f_{n_k}(t_m)\}$ is a Cauchy
sequence, there exists $k_0$ such that
$d(f_{n_{k_1}}(t_m),f_{n_{k_2}}(t_m))<\epsilon/3$ for $k_2>k_1\geq
k_0$
and, so: 
\[
\begin{array}{rl} d(f_{n_{k_1}}(t),f_{n_{k_2}}(t))\leq & d(f_{n_{k_1}}(t),f_{n_{k_1}}(t_m))+d(f_{n_{k_1}}(t_m),f_{n_{k_2}}(t_m))\\ & +d(f_{n_{k_2}}(t_m),f_{n_{k_2}}(t)) < \frac{\epsilon}{3}+\frac{\epsilon}{3}+\frac{\epsilon}{3}=\epsilon,\end{array}
\]
as required. Let $f(t)$ be the limit point of the sequence
$\{f_{n_k}(t)\}_{k}$. The continuity of the so obtained $f$ follows  from the
equicontinuity of the sequence $\{f_{n_k}\}$. \cvd

\smallskip

In order to deal with oriented-equicontinuity, let us consider the following result.

\begin{lemma}\label{oeqeq}
Let $(M,d)$ be a generalized metric space associated to a Finsler
manifold, and $\{f_n\}$, $f_n:[a,b]\rightarrow M$, a sequence of
oriented-equicontinuous functions satisfying $f_n([a,b])\subset
K$, with $K$ compact. Then, the sequence $\{f_n\}$ is
equicontinuous.
\end{lemma}
{\it Proof.}
Assume by contradiction that $\{f_n\}$ is not equicontinuous.
Then, there exist $\epsilon>0$, two sequences
$\{t_m^1\},\{t_m^2\}\subset [a,b]$ and a subsequence
$\{f_{n_m}\}_m$ such that $0<t^{2}_{m}-t^{1}_{m}\searrow 0$ and
\begin{equation}\label{auxi1}
d(f_{n_{m}}(t^2_m),f_{n_{m}}(t^1_m))>\epsilon\quad \hbox{for all
$m$.}\end{equation} Even more, since $\{f_{n}\}$ is
oriented-equicontinuous, we can also assume that
\begin{equation}\label{auxi2}
d(f_{n_{m}}(t^1_m),f_{n_m}(t^2_m))<1/m\quad \hbox{for all
$m$.}\end{equation} As $\{f_{n_m}(t_m^i)\}_{m}$, $i=1,2$, is
contained in the compact $K$, up to a subsequence, we can assume
that $\{f_{n_m}(t_m^i)\}_m\rightarrow x^i$ for some $x^i\in K$,
$i=1,2$. From (\ref{auxi1}), $d(x^2,x^1)\geq\epsilon$, but from
(\ref{auxi2}), $d(x^1,x^2)=0$, in contradiction to the fact that
$d$ is a generalized distance. \cvd

%

\begin{remark}\label{remorequi}{\em
As a consequence, if  $(M,d)$ is locally compact and $\{f_n\}_n$
converges pointwise to a continuous function $f$, then the family
$\{f_n\}_{n\geq n_0}$ is equicontinuous for some $n_0$. In fact,
the compactness of  $f([a,b])$ implies the existence of a compact
neighborhood $K$ of $f([a,b])$, and the inclusion
$f_{n}([a,b])\subset K$ follows from the
uniform convergence of $\{f_n\}$ to $f$.} 
\end{remark}
The required version of Arzela's Theorem is the following:
\begin{theorem}\label{taf}
Let $(M,d)$ be a  generalized metric space and $\{f_n\}$,
$f_{n}:[a,b]\rightarrow M$, a sequence of
oriented-e\-qui\-con\-ti\-nuous functions. Suppose that
$\{f_n(a)\}$ has an accumulation point and $(M,d)$ is locally
compact and complete. Then, there exist a continuous function
$f:[a,b]\rightarrow M$ and an equicontinuous subsequence
$\{f_{n_k}\}\subset \{f_n\}$ such that $\{f_{n_k}(t)\}\rightarrow
f(t)$ for each $t\in [a,b]$.
\end{theorem}
{\it Proof.} By Remark \ref{remorequi}, it suffices to prove the
existence of a continuous function $f:[a,b]\rightarrow M$ and a
subsequence $\{f_{n_k}\}\subset \{f_n\}$ such that
$\{f_{n_k}(t)\}\rightarrow f(t)$ for each $t\in [a,b]$. To this
aim, consider the set $$B=\left\{r\in (a,b]:\;
\begin{array}{l}\exists \hbox{ subsequence of $\{f_n\}$ which converges
(pointwise)}\\ \hbox{to a continuous function $f:[a,r]\rightarrow
M$}\end{array} \right\}.$$ Reasoning as in
Theorem \ref{tafe}, there exists $r>0$
such that $f_n([a,a+r])\subset K$ for all $n$ and for some compact
set $K$. From Lemma \ref{oeqeq}, the family $\{f_n\}$ is
equicontinuous in $[a,a+r]$, and from Theorem \ref{tafe},
necessarily $B\neq \emptyset$. Let $r_0\leq b$ be the supremum of
$B$. By a canonical diagonal argument, there exists a continuous
function $f:[a,r_0)\rightarrow M$ which is the pointwise limit of
a partial $\{f_{n_k}\}$ of $\{f_n\}$. Then, define recursively a
sequence $\{s_m\}$ in the following way:
$$s_m\geq s_{m-1},\qquad0<r_0-s_m<\delta_m<\min\{\delta_{m-1},1/m\},$$ where $\delta_0=1$ and $\delta_m$ is obtained from the oriented-equicontinuity of $\{f_{n_k}\}$ applied to $\epsilon=1/m$. The sequence $\{s_m\}$ satisfies that $\{s_m\}\nearrow r_0$ and $\{f(s_m)\}_m$ is a Cauchy sequence. From the completeness of $M$, $\{f(s_m)\}_m$ converges to some point $f(r_0)$. As $M$ is locally compact, there exists $r>0$ such that $B^+(f(r_0),r)\subset K'$ with $K'$ compact. Since $f(s_m)\rightarrow f(r_0)$ and $d$ is a generalized distance, there exists $m$ such that $d(f(r_0),f(s_m))<r/3$ and $1/m<r/3$. Moreover, as $f_{n_k}(s_m)\rightarrow f(s_m)$, for $n_k$ big enough $d(f(s_m),f_{n_k}(s_m))<r/3$. Therefore, we have that, for all $s\in [s_m,s_m+\delta_m]$ and $n_k$ big enough,
\[
\begin{array}{rl}d(f(r_0),f_{n_k}(s))\leq & d(f(r_0),f(s_m))+d(f(s_m),f_{n_k}(s_m))+d(f_{n_k}(s_m),f_{n_k}(s))\\ < &
\frac{r}{3}+\frac{r}{3}+\frac{r}{3}=r.
\end{array}
\]
So, for $n_k$ big enough, $f_{n_k}([s_m,s_m+\delta_m])\subset K'$ and then, by Lemma \ref{oeqeq}, the family $f_{n_k}$ is equicontinuous in $[s_m,s_m+\delta_m]$. Moreover, as $r_0-s_m<\delta_m$, necessarily $r_0\in [s_m,s_m+\delta_m]$. Therefore, from the equicontinuity and the fact that $f(s_m)\rightarrow f(r_0)$, we deduce that $f_{n_k}(r_0)\rightarrow f(r_0)$. Finally, by Theorem \ref{tafe}, $r_0$ must be a maximum, and thus, $r_0=b$. \cvd

\smallskip

Now we are in conditions to prove:
\begin{corollary}\label{ccg2}
The Cauchy-Gromov completion $M_{CG}^+ =M\cup
\partial_{CG}^+ M$ satisfies:
\begin{itemize}
\item[(1)] (Openness of $M$). The boundary $\partial_{CG}^+ M$
is closed in $M_{CG}^+$.

\item[(2)] (Heine-Borel). The closure in $M_{CG}^+$ of any
bounded set in $M$ is compact.

\item[(3)] (Consistency with the Cauchy boundary). If {\cambios $M_C^+$ is locally compact} and $d$ is a
generalized distance, then $j^+: M_C^+ \rightarrow M_{G}^+$ is an
embedding and $M_C^+\equiv M_{CG}^+$, i.e. $\partial_{RG}^+
M=\emptyset$.
\end{itemize}
\end{corollary}
{\it Proof.} The first assertion is straightforward from
Theorem \ref{tgromovf}. The second assertion follows from the compact character of $M_G^+$ and the definition of $M_{CG}^+$. For
the first assertion in (3), as $d$ is a generalized distance,
condition (a4') holds, and thus, $j^+$ is continuous. The
continuity of the inverse of $j^+$ follows reasoning as in Remark
\ref{mclc}. For the last assertion, in analogy to the Riemannian
case, it is sufficient to show that $\overline{B^+(x,r;M_C^+)},
x\in M_C^+,r>0$, are compact. Consider a sequence $\{x_n\}$ in
$B^+(x,r;M_C^+)$ and define a sequence of curves $\{\gamma_n\},
\gamma_n:[0,1]\rightarrow M$ joining $x$ and $x_n$, such that:
each $\gamma_n$ is contained in $M$, it is parametrized with
constant speed and length$(\gamma_n)<2r$. Reasoning as in the
Riemannian case (see footnote \ref{footequicon}), the family of
curves is oriented-equicontinuous. So, taking into account that
$M_C^+$ is locally compact and complete,
and $x$ is an
accumulation point of $\{\gamma_n\}$, Theorem \ref{taf} ensures
the existence of a curve $\gamma_\infty$ which is the pointwise
limit of a subsequence of $\{\gamma_n\}$, and so the required
limit is $\gamma_\infty (1)$.\cvd

\subsection{Busemann completions}\label{sectbusf}

Next, we are going to develop the Busemann completion $M_B^+$ for
any Finsler manifold $(M,F)$. {\cambios First, in Subsection
\ref{s521}, we will deal with its definition as a point set by
making definitions analogous to the Riemannian ones. The topology
is developed in  the following three subsections.} As always in
the non-symmetric case, a second completion $M_B^-$ is possible,
and some non-trivial conventions about it will be stressed at
Subsection \ref{s525}.

\subsubsection{Busemann completion as a point set} \label{s521} Let $(M,d)$ be a generalized metric space associated to a Finsler manifold $(M,F)$ or any space under Convention \ref{cginebra2}.
We will introduce the analogous objects to those ones in Section
\ref{s42} for the Riemannian case.

Let $C^+(M)$ be the set of curves $c:[\alpha,+\Omega)\rightarrow
M$, $\Omega\leq \infty$,
such that $F(\dot{c})< 1$ (recall Remark \ref{remm1}).
 For $c\in C^+(M)$
the associated {\em (forward) Busemann function} $b_c^{+}: M\rightarrow (-\infty, \infty]$ is:
\begin{equation}
 \label{ebusemann+}
b_c^{+}(\cdot)=\lim_{s\rightarrow\Omega}(s- d(\cdot,c(s)),
\end{equation}
which is well-defined as
\begin{equation}\label{equdecre}s\mapsto
s-d(x,c(s))\;\; \hbox{ is increasing for any} \, x\in M
\end{equation}
(extend Lemma \ref{l1}). 
Following steps analogous to those in Proposition \ref{p1}, first,
we have  that if a Busemann function is not identically infinity
then it is finite-valued everywhere, and we denote:
$$B^+(M):=\{b_{c}^{+}<\infty: c\in C^{+}(M)\}.$$
Second, we also have
$$
B^+(M) \subset \lipd, \quad \hbox{and so,} \quad b_c^+\in B^+(M) \Rightarrow [b_c^+]\in M^+_G .
$$
Moreover, when $\Omega<\infty$ there exists some $\bar x\in M_C^+$ such that
$b_c^+(x)=\Omega-d(x,\bar x)$ for all $x\in M$ (here, $d$ is extended to the elements of $M_C^+$
according to Prop. \ref{Srelat}), and we will write
\begin{equation}\label{equuu2}
b_c^+ =d_p^+ :=\Omega-d^+_{\bar x},
\end{equation}
where $d^+_{\bar x}(x):= d(x,\bar x)$ for all $x\in M$ and 
$p=(\Omega,\bar x)\in \R\times M_C^+ $
(recall that the domain of $d^+_{\bar x}$ can be also extended to any point of $M_C^+$).
If, otherwise, $\Omega=\infty$, then $b^{+}_c$ is called a {\em
properly Busemann} function
and we write: $${\cal
B}^+(M):=\{b^+_c<\infty:\, c\in C^{+}(M),\, \Omega=\infty\}.$$

\begin{definition} \label{ddd0} As a point set, the {\em (forward) Busemann completion} of $(M,d)$
 is the quotient $M^{+}_{B}:= B^+(M)/\R \subset M_G^+$, the {\em (forward) Busemann boundary} is
$\partial^{+}_{B}M:=M^{+}_{B}\backslash M \; (\subset
\partial^{+}_GM)$, and the {\em (forward) properly Busemann boundary} is $\partial^{+}_{{\cal B}}M:= {\cal B}^+(M)/\R$.
\end{definition}
From previous discussion, clearly,  $\partial^{+}_{B}M
=j^+(\partial_C^+ M) \cup \partial^{+}_{{\cal B}}M.$

\subsubsection{Chronological topology on $M_B^+$}

The chronological topology on $B^+(M)$ is introduced defining  the
closed sets by means of a limit operator $\hat L$. Such an
operator has its origins in an operator introduced by Flores and
Harris (\cite{FH}; see also \cite{H,H2}) in order to study the
future causal boundary of a static spacetime, and it has been
developed further, including some consequences for the Riemannian
case (see \cite{FH,F,AF, Ha07, Orl-M, FHS0}). Here, we introduce
{\cambios directly the chronological topology} for the Busemann
completion of a Finsler manifold, $M_B^+$. The definition is
somewhat involved, and its motivation will appear more clearly in
stationary spacetimes. However, one can appreciate its nice
properties and its relation with Gromov topology
---even in Finsler manifolds which are not involved {\cambios in
constructions associated to  boundaries of spacetimes.}

{\cambios The limit operator $\hat L$ will allow to define the
closed sets for the topology in $B^+(M)$.} Let
$\sigma=\{f_n\}\subset B^+(M)$, the subset $\hat{L}(\sigma)
\subset B^+(M)$ is defined as:
\begin{equation}\label{ecu1}
f\in \hat{L}(\sigma)\iff \left\{\begin{array}{l}(a)\; f\leq
\liminf_{n}f_{n} \hbox{ and}\\ (b)\; \forall g\in B^{+}(M) \hbox{
with } f\leq g\leq \limsup_{n} f_{n}, \hbox{ it is } g=f.
\end{array}\right. \end{equation}
\begin{definition}\label{ddd1}
The {\em (forward) chronological topology} ({\em chr. topology}
for short) on $B^+(M)$ is the one such that $C\subset B^+(M)$ is
closed if and only if $\hat{L}(\sigma)\subset C$ for any sequence
$\sigma=\{f_n\}\subset C$,

The induced quotient topology
on the Busemann completion 
$M_B^+=B^+(M)/\R$ will be also called {\em Busemann} or {\em chr.
topology} on $M_B^+$.
\end{definition}

The Busemann completion $M_B^+$ will be always regarded as a
topological space endowed with the chr. topology, except if
otherwise is explicitly said.

\begin{remark}\label{rsequential} 
{\em It is easy to check that this notion of closed set is
consistent and, moreover, if $f\in \hat L(\sigma)$ then $\sigma$
converges to $f$. In fact, this can be seen from the following
more general viewpoint (see \cite[Appendix 3.6]{FHS0} for details;
recall also Remark \ref{propsimplepunt}(1) above).

A {\em limit operator} $L$ on a set $X$ is any map $L : {\cal
S}(X) \rightarrow {\cal P}(X)$ (where ${\cal S}(X)$ is the set of
all the sequences in $X$ and ${\cal P}(X)$ is the set of parts of
$X$) such that if $\bar \sigma$ is a subsequence of some
$\sigma\in S(X)$ then $L(\sigma)\subset L(\bar{\sigma})$. This
property ensures that an {\em associated topology} $\tau_L$ is
defined on $X$ just declaring that $ C \subset X$ is closed if and
only if $L(\sigma) \subset C$ for any sequence $\sigma$ in $C$.
For this topology, $x\in L(\sigma)$ implies that $\sigma
\rightarrow x$. The converse does not hold in general and, when it
holds, $L$ (and, then, $\tau_L$) is called {\em of first order}.

 Such a topological space $(X, \tau_L)$ is a {\em sequential space}, i.e.,
it satisfies that any {\em sequentially closed} subset ({\cambios which
is} a subset $S\subset X$ such that if a {\cambios sequence in} $S$
converges to a point $x\in X$ then $x\in S$) is closed {\cambios
---recall that} the converse always holds.  In the class of
sequential spaces (which includes the class of first countable
ones) the continuity of a function can be characterized as
sequential continuity, i.e., in terms of the preservation of
converging sequences, as in the standard case of metric spaces.
This fact will be used frequently for maps between spaces such as
$M_B^+$, $B^+(M)$, $M_G^+$ or $M_C^+$, as all these spaces are
sequential (the last two ones because they are metrizable, and
$B^+(M)$, $M_B^+$ by construction,
see also Proposition \ref{gh}).

Finally, consider the constant sequence $\hat{x}=\{x\}_n$, for
each $x\in X$. Clearly, if $L(\hat{x})=\{x\}$ for all $x\in X$,
then the points in $X$ are closed for $\tau_L$, i.e. $\tau_L$ is
$T_1$ (conversely, if $\tau_L$ is $T_1$ then $L(\hat{x}) \subset
\{x\}$). In particular, the chr. topology is always $T_1$.
}
\end{remark}


An example of result (to be used later) which is valid for the
chr. topology as well as any topology constructed from a limit
operator $L$, is the following. Recall that it is applicable to
the limit operator $\hat L$ for the chr. topology in $M_B^+$.

\begin{proposition}\label{ppo} Let $X$ be a set endowed with a
topology associated to a limit operator $L$. Assume that
$L(\hat{x})=\{x\}$ for any $x\in X$, where $\hat{x}=\{x\}_{n}$.
Given a sequence $\sigma\subset X$, if
$L(\overline{\sigma})=\{x\}$ for any subsequence
$\overline{\sigma}\subset \sigma$ then {\cambios $\sigma$ only
converges to $x$.}
\end{proposition}
{\it Proof.} From Remark \ref{rsequential}, $x$ is a limit of
$\sigma$ and we will prove that it is unique. Assume by
contradiction that $x'\neq x$ is also limit of $\sigma$. In
particular, $x,x'$ are limits of any subsequence
$\overline{\sigma}$ of $\sigma$. Notice  that  {\cambios
$\overline{\sigma}\cup \{x\}$ is closed. In fact, if $\tau$ is any
sequence constructed from the elements in $\overline{\sigma}\cup
\{x\}$ then some subsequence $\overline{\tau}$ of $\tau$ must be
either a subsequence of $\overline{\sigma}$ or constantly equal to
some element of $\overline{\sigma}\cup \{x\}$. In both cases, our
assumptions imply $L(\overline{\tau})\subset \overline{\sigma}\cup
\{x\}$, and then $L(\tau)\subset \overline{\sigma}\cup \{x\}$ as
required. As a consequence, $x'\in\overline{\sigma}$ for any
subsequence $\overline{\sigma}\subset \sigma$. This means that
$\sigma$ is constantly equal to $x'$ up to a finite number of
terms, in contradiction to  $L(\sigma)=\{x\}$.}\cvd

The definition of the chr.
topology makes necessary to study properties of convergent
sequences.

\begin{lemma} \label{ld} (i) For $p_i=(t_i,x_i)\in \R\times M$, $i=1,2$, the inequality
$d_{p_1}^+<d_{p_2}^+$ (i.e. $d_{p_1}^+(x)<d_{p_2}^+(x)$ for all
$x\in M$) holds if and only if $d(x_1,x_2)<t_2-t_1$.

(ii) For any $f\in B^+(M)$ there exists a sequence $\{p_n\}\subset \R\times M$ such that $d^+_{p_n}<d^+_{p_{n+1}}$ and
$f =\sup_n d_{p_n}^+$. Moreover, if $D\subset \R\times M$ is a dense subset, $\{p_n\}$ can be chosen included in $D$.

(iii) Let $g\in \lipd$ and $p_0=(t_0,x_0)\in \R\times M$ such that
$t_0<g(x_0)$. Then, $d_{p_0}^{+}<g$.
\end{lemma}
{\em Proof}. (i) The implication to the right follows by
evaluating the inequality $d_{p_1}^+<d_{p_2}^+$ at $x_{1}$. For
the implication to the left, just use that, from the triangle
inequality, $d^{+}_{x_{2}}(y)-d^{+}_{x_{1}}(y)\leq d(x_{1},x_{2})$
for all $y$.

(ii) Let $f=b_c^+$ for $c: [\alpha, \Omega)\rightarrow M$.
Any sequence $\{p_n=c(t_n)\}$ with $\{t_n\}\nearrow \Omega$ fulfills the
required properties. For the last assertion, recall that each $W_n:=\{p\in \R\times M: d^+_{p_n}<d^+_{p}<d^+_{p_{n+1}}\}$ is open
and non-empty (use {\em (i)}), and  replace the original $p_n$  by any point in $D\cap W_n$.

(iii) Notice that $d_{p_0}^+(x) < g(x_0)-d(x,x_0)\leq g(x)$, the
last inequality from the Lipschitz condition. \cvd

{\cambios The following result exhibits some natural open subsets.}
\begin{proposition} \label{lopen} For each $f\in B^+(M)$ the subset
\begin{equation}\label{eimasf}
{\cambios I^+(f)}:= \{h\in {B}^{+}(M): \;\exists\; p\in\R\times M\;
\hbox{such that}\; f\leq d^+_{p}<h\}\subset B^{+}(M)
\end{equation}
is open for the chr. topology.
\end{proposition}
{\em Proof}. In order to prove that $B^+(M)\backslash I^+(f)$ is
closed, it is enough to check that, for any sequence $\sigma =
\{g_n\}\subset B^{+}(M)$ and function $h\in I^+(f)$ such that
$h\in \hat L(\sigma)$, the sequence $\sigma$ is contained in
$I^+(f)$, up to a finite subset.
Let $p_0=(t_0,x_0)\in \R\times M$
such that $f\leq d_{p_0}^+<h$, recall that
$t_0=d_{p_0}^+(x_0)<h(x_0)$
{\cambios and, as $h\leq \liminf g_n$, $d_{p_0}^+(x_0)<g_n(x_0)$
for $n$ big enough}. So, from Lemma \ref{ld} (iii)
 $d_{p_0}^+<g_n$  and
$g_n\in I^+(f)$ for $n$ large enough. \cvd

\smallskip {\cambios
In the particular case of $f=d_{p_0}^+$, for $p_0\in \R\times M$,
the notation will be replaced by $I^+(p_0,B^+(M)):=
I^+(d_{p_0}^+)$. Alternatively, by using the Lipschitz property of
$h$ in (\ref{eimasf}):
\begin{equation}\label{ee}
I^{+}(p_0,B^{+}(M)):=\{h\in B^+(M): t_0< h(x_0)\} \quad \hbox{for
any} \; p_0=(t_0,x_0)\in \R\times M.
\end{equation}
The notation suggests that $I^{+}(p_0,B^{+}(M))$ is related with
the chronological future of $p_0$ in some spacetime, as shown
explicitly in Section \ref{s6}.}

 {\cambios The
natural topological properties of $B^+(M)$ are transmitted to
$M_B^+$}. In fact, as the maps $T_k: B^+(M) \rightarrow B^+(M)$,
$f\mapsto f+k$ are homeomorphisms for all $k\in\R$, we have:

\begin{proposition}\label{gh}  The natural projection  $\Pi: B^+(M)\rightarrow M^+_B=B^+(M)/\R$, is an open map.


\end{proposition}

%

\subsubsection{The inclusion of $M_C^+$ in $M_B^+$ from the topological viewpoint}

In order to obtain a true completion, the original space $M$ has
to be embedded as a dense subset in $M_B^+$ in a natural way. In
fact, we know that $j^+:M_C^+\rightarrow M_G^+$ is injective, and
$j^+(M_C^+)\subset M_B^+$. Denote
\begin{equation}\label{jbmas}
j_B^+:M^+_C\rightarrow M_B^+
 \end{equation}
the map obtained by  restricting the codomain of $j^{+}$ to
$M_B^+$ (which is regarded now as a topological space with the
chr. topology).
We will see that $j^+_B|_{M}$ yields the required embedding and, as in the case of the embedding in $M_G^+$, the situation is subtler
for $j^+_B$. 

\begin{proposition}\label{lemmabusemanM}
Consider a sequence $\{b_{c_n}^+\}\subset B^+(M)$ and
$p=(\Omega,x)\in \R\times M$. {\cambios If
$d_p^+=\Omega-d_x^+ (\in B^+(M))$ belongs to
$\hat{L}(\{b_{c_n}^+\})$, then $\{b_{c_n}^+\}$ converges pointwise
to $d_p^+$.} Moreover, for $n$ big enough,
$b^+_{c_n}=\Omega_n-d^+_{x_n}$ with $\Omega_n\in \R$ and $x_n\in
M$ satisfying $\Omega_n\rightarrow \Omega$ and $x_n\rightarrow x$.

{\cambios If $d_Q$ is a generalized distance and $M_C^+$
(necessarily equal  to $M_C^-$ and $M_C^s$, Prop. \ref{genersym})
is locally compact, then the same conclusion holds for any
$p=(\Omega,x)\in \R\times M_C^+$.}
\end{proposition}
{\it Proof.} Consider $c_n:[\alpha_n,\Omega_n)\rightarrow M$ and
suppose by contradiction that {\cambios $\Omega-d^+_x\in
L(\{b^+_{c_n}\})$ but $\{b^+_{c_n}\}$ does not converge pointwise
to $\Omega-d^+_x$.} By the chronological convergence, one has
$\Omega\leq \liminf b^+_{c_n}(x)$. So, taking a subsequence
$\{\Omega_{k_{n}}\}_{n}$ if necessary, and re-naming
$\Omega_{k_{n}}\equiv \Omega_{n}$ for all $n$, we can assume
$\Omega_n\geq \Omega-1/n$. Notice also that there exists some
$\epsilon>0$ such that:
\begin{equation}\label{lemeq1}\Omega_n>{\cambios \Omega+\epsilon, }\quad\hbox{up to a subsequence.}
\end{equation}
In fact, otherwise $\Omega_n\rightarrow \Omega$ and for $n$ large
enough, $b^+_{c_n}=\Omega_n-d^+_{x_n}$ for some $x_n\in M^+_C$
(Prop. \ref{p1}). Since $\Omega\leq \liminf_{n}b^{+}_{c_{n}}(x)$,
{\cambios for each $m\in\N$ there exists some $m_0$ such that $\Omega -
1/m\leq b^+_{c_n}(x)=\Omega_n-d(x,x_n)$} for all $n\geq
n_{0}\geq m$, that is:
$$
0\leq d(x,x_n)\leq (\Omega_n-\Omega ) +\frac{1}{m}.$$
Thus, $d(x,x_{n})\rightarrow 0$, which joined to
$\Omega_{n}\rightarrow\Omega$, implies that $\{b^+_{c_n}\}$
converges pointwise to $\Omega-d^+_x$, a contradiction.

Next, assuming that (\ref{lemeq1}) holds, we are going to obtain a
contradiction with the maximality of $\Omega-d^+_x$ as {\cambios a
$\hat L$-limit according to (b) in (\ref{ecu1})}. The local
compactness of $M$ allows to take a compact neighborhood $U$ of
$x$ and define the following subsets of Busemann functions, which
are also compact with the pointwise topology:
\begin{equation}\label{lemeq3}\begin{array}{c} C:=\{t-d^+_y\in B^+(M): t\in [\Omega-\frac{\epsilon}{2},\Omega+\frac{\epsilon}{2}], y\in U\}\\
\partial C:=\{t-d^+_y\in C:\;\hbox{either}\;\; t=\Omega-\frac{\epsilon}{2},\Omega+\frac{\epsilon}{2}\;\; \hbox{or} \;\; y\in \partial U\}.\end{array}\end{equation}
For the proof, we will need the following {\it claim} which will
be proved later:
\begin{itemize}
\item[] {\bf Claim:} For any  $s\in
[\Omega-\frac{\epsilon}{2},\Omega]$ and $z\in U$ such that
$s<\Omega-d(z,x)$, there exist $t^z\in
[\Omega-\epsilon/2,\Omega+\epsilon/2]$, $y^z\in U$ and a
subsequence $\{b^+_{c_{n_k}}\}$ such that: (i) $t^z-d^+_{y^z}\in
\partial C$, (ii) $s\leq t^z-d^+_{y^z}(z)$ and (iii) $\liminf
b^+_{c_{n_k}}\geq t^z-d^+_{y^z}$.
\end{itemize}
Let us apply the claim to each $s_m$ of a sequence $\{s_m\}\subset
[\Omega-\frac{\epsilon}{2},\Omega]$ with $s_m\nearrow \Omega$ and
$z=x$.
We find some $t^m\in
[\Omega-\frac{\epsilon}{2},\Omega+\frac{\epsilon}{2}]$, $y^m\in U$
and a subsequence $\{b_{c_n}^{+ [m]}\}_n$ of $\{b^+_{c_n}\}$ satisfying
(i), (ii) and (iii). Making each $\{b_{c_n}^{+ [m+1]}\}_n$ be a
subsequence of $\{b_{c_n}^{+ [m]}\}_n$, by a diagonal argument we can
suppose that the same subsequence $\{b^+_{c_{n_k}}\}$ is valid for
all $m$. Up to a subsequence, $\{t^m\}$ and $\{y^m\}$ converge to
some $t^0\in
[\Omega-\frac{\epsilon}{2},\Omega+\frac{\epsilon}{2}]$, $y^0\in U$
and, so,
\begin{equation}\label{lemeq4}
\{t^n-d^+_{y^n}\} \quad\hbox{converges pointwise to}\quad
t^0-d^+_{y^0}\in
\partial C.
\end{equation}
By statements (ii) and (iii) of the claim, the function
$t^{0}-d^+_{y^{0}}$ also satisfies $\Omega\leq t^0-d^+_{y^0}(x)$ and
$\liminf b^+_{c_{n_k}}\geq t^0-d^+_{y^0}$. The first inequality (plus
the triangle one) yields
 $t^0-d^+_{y^0}\geq \Omega-d^+_x$ and, as
$t^0-d^+_{y^0}\in \partial C$, the inequality is strict at some
point. Thus, $t^0-d^+_{y^0}$ contradicts the maximality of
$\Omega-d^+_x$ in $\limsup b^+_{c_n}$, as required.
%

{\it Proof of the claim.} As $\liminf b^+_{c_n}\geq
\Omega-d^+_{x}$, {\cambios and so, $s<\Omega-d^+_x(z)\leq \liminf
b^+_{c_{n}}(z)$,} 
there exists $t_n\in (\Omega+\epsilon, \Omega_n)$ (observe
(\ref{lemeq1})) satisfying:
$$b^+_{c_{n}}(z)>t_{n}-d(z,c_n(t_n))\quad \hbox{and}\quad s<t_{n}-d(z,c_n(t_n))\;\;\hbox{for all $n$ big enough.}$$

That is, $d(z,c_{n}(t_n))<t_n-s$, and there exists a curve
$r:[s,t_n]\rightarrow M$ with $F(\dot{r})<1$ joining $z$ and $c_n(t_n)$. By Lemma
\ref{l1}, $\tau-d^+_{c(\tau)}$ is increasing, and so,
 $$s-d^+_z< t-d^+_{r(t)}<t_n-d^+_{c_n(t_n)}\quad \forall t\in (s,t_n).$$
As $s-d^+_z\in C$, $t_n-d^+_{c_n(t_n)}\notin C$ and the functions
$t-d^+_{r(t)}$ vary continuously with $t$ in the pointwise
convergence topology,
\begin{equation}\label{lemeq5} \exists\;\; t^z_n: t^z_n-d^+_{c_n(t_n^z)}\in \partial C.\end{equation}
Summing up, we have found $t^z_n,y^z_n(=r(t^z_n))$ such that:
$t_n^z-d^+_{y_n^z}\in \partial C$ {\cambios and $s-d^+_z<
t^z_n-d^+_{y_n^z}<t_n-d^+_{c_n(t_n)}\leq b^+_{c_n}$,
 {\cambios the last inequality because of the Lipschitz condition
(moreover, it is strict because of our choice $F(\dot c)<1$,
recall Remark \ref{remm1})}.
As $t^{z}_{n}-d^+_{y^{z}_{n}}\in\partial C\subset C$ compact,
there exists $t^z, y^z$ such that $t^z-d^+_{y^z}$ is the pointwise limit of a subsequence $\{t^z_{n_k}-d^+_{y^z_{n_k}}\}_k$.
From the construction, $t^z, y^z$ satisfy the conditions (i) and
(ii). And, as $t^z-d^+_{y^z}$ is the pointwise limit {\cambios  of
$t^z_{n_k}-d^+_{y^z_{n_k}}\leq b^+_{c_{n_k}}$,} the condition
(iii) is also satisfied, which ends the proof of the claim.

For the last assertion, observe that the essential properties used
from $(M,d)$ have been: {\cambios (1) the fact that $(M,d)$ is a
length space;} (2) the local compactness of $M$; (3) in order to
prove (\ref{lemeq3}), (\ref{lemeq4}), the sentence below
(\ref{lemeq1}) and the convergence of
$\{t^z_{n_k}-d^+_{y^z_{n_k}}\}_k$ above, we have used the
continuity of the function $d$ in the second variable; (4) in
order to prove (\ref{lemeq5}), we have essentially used that any
curve $c$ going from inside to outside the compact neighborhood
$U$ of $x$ (which can be assumed to be a ball) must intersect the
boundary of $U$, i.e. the continuity of $d$ in the first variable
(recall that we consider backward balls in general). So, if we
impose that {\cambios $d_Q$ is a generalized distance and $M_C^+$  is
locally compact} (thus, {\cambios $M_C^+=M_C^s$ is a length
space and} the continuity in both variables is ensured; recall
{\cambios Proposition \ref{genersym} and} Lemma \ref{lcf} applied
to $d$ and $\rd$), the last assertion is obtained. \cvd

\smallskip

\noindent The required result is then:

\begin{corollary}\label{cembbusf} Let $(M,d)$ be under Convention \ref{cginebra2} and $j_B^+: M_C^+\rightarrow M_B^+$ as in (\ref{jbmas}).
Then, $j^{+}_B$ is injective and it is continuous if and only if
the condition (a4') holds. The restriction $j^+_B|_M$ is an   {\cambios open} embedding {\cambios and $j^+_B(M)$ is dense in
$M_B^+$.}

Moreover, if   $d_Q$ is a generalized distance and $M^+_C$ is
locally compact, then $j^+_B$ is also an embedding.  In
particular, this happens if $d$ is a (symmetric) distance and $M_C
(=M^+_C)$ is locally compact.
\end{corollary}
{\it Proof.} The continuity and injectivity of $j^{+}_{B}$ follows
as in Proposition \ref{pembf},  {\cambios the density of $j^+_B(M)$
from Lemma \ref{ld}(ii), and $j^+_B(M)$ is open because
$\Pi^{-1}(\partial_B^+M)\subset B^+(M)$ is closed for the chr.
topology (for any sequence $\sigma$ in $\Pi^{-1}(\partial_B^+M)$,
$\hat{L}(\sigma)$ is also included in $\Pi^{-1}(\partial_B^+M)$ by
Proposition \ref{lemmabusemanM}).
The continuity of the inverse of {\cambios of
$j^+_B$} follows as a consequence of previous Proposition
\ref{lemmabusemanM} {\cambios and the fact that $M_{B}^{+}$ is
sequential (recall Remark \ref{rsequential})}. \cvd

\begin{remark}\label{rem} {\em
Notice that analogous results are also being obtained for
$B^+(M)$. In particular, the map $\R\times M\hookrightarrow
B^+(M)$, $(t,x)\mapsto t-d_x^+$ is an open embedding with dense
image.
}\end{remark}

\subsubsection{Compactness of the Busemann completion}

Our aim is to prove that $M_{B}^+$ provides a {\cambios sequential
compactification}
of $M$. To this objective, again we need a previous technical result.

\begin{lemma}\label{legg} Assume that, for some $\phi\in \lipd$, $x_{0}\in M$ and $t_{0}\in\R$, the following subset of $B^{+}(M)$ is non-empty:
\begin{equation}\label{equf}
{\cal F}=\{g\in B^{+}(M): g\leq \phi,\;\; g(x_{0})\geq t_{0}\}.
\end{equation}
Then, ${\cal F}$ admits some maximal element (for the partial
order $\leq$) in $B^{+}(M)$.
\end{lemma}
{\it Proof.} Applying Zorn's Lemma to ${\cal F}$, it suffices to
show that any totally ordered subset admits an upper bound in
${\cal F}$. By contradiction, suppose that there exists a totally
ordered subset ${\cal F'}$ which has no upper bound. Associated to
this set, there exists a countable family of functions
$\{d_{(t_n,x_n)}^+\}$ such that any $b_c^+\in {\cal F'}$ is
chronological limit of some increasing subsequence of
$\{d_{(t_n,x_n)}^+\}$ and each element $d_{(t_n,x_n)}^+$ of the
family admits some $b_{c}^+\in {\cal F'}$ satisfying
$d_{(t_n,x_n)}^+ <b_{c}^+$. In fact, from Lemma \ref{ld} (ii),
$b_c^+=\lim_n d_{(t_n^c,x_n^c)}^+$ (pointwise, and then,
chronologically), with
$d^{+}_{(t^{c}_{n},x^{c}_{n})}<d^{+}_{(t^{c}_{n+1},x^{c}_{n+1})}$,
$(t_{n}^{c},x_{n}^{c})\in \mathbb{Q}\times D$ for all $n$, and $D$
a numerable, dense subset of $M$.

Construct a subsequence $\{d_{(s_n,y_n)}^+\}$ of
$\{d_{(t_n,x_n)}^+\}$ by the following inductive way. Take
$d_{(s_1,y_1)}^+:=d_{(t_1,x_1)}^+$, and given $d_{(s_n,y_n)}^+$,
take $d_{(s_{n+1},y_{n+1})}^+$ such that $d_{(s_{n+1},y_{n+1})}^+>
d_{(s_{n},y_{n})}^+$ and
$d_{(s_{n+1},y_{n+1})}^+>d_{(t_i,x_i)}^+$, $i=1,...,n-1$. For the
existence of $d_{(s_{n+1},y_{n+1})}^+$, note that, as
$d_{(t_i,x_i)}^+ <b_{c_i}^+\in {\cal F'}$, $i=1,...,n-1$,
$d_{(s_{n},y_{n})}^+< b_{c_n}^+\in {\cal F'}$ and ${\cal F'}$ is
totally ordered without upper bound, there exists
$b_{c_{n+1}}^+\in {\cal F'}$ such that $d_{(t_i,x_i)}^+<
b_{c_{n+1}}^+, i=1,...,n-1$ and $d_{(s_{n},y_{n})}^+<
b_{c_{n+1}}^+$; hence, the density of $\{d_{(t_{n},x_{n})}^+\}$ in
${\cal F'}$ ensures that such a $d_{(s_{n+1},y_{n+1})}^+$ exists.

Finally, from Lemma \ref{ld} (i) one has that
$d_{(s_{n+1},y_{n+1})}^+>d_{(s_{n},y_{n})}^+$ implies
$s_{n+1}-s_{n}>d^+(y_n,y_{n+1}),$ and so, one can define a curve
$\tau:[s_1,\Omega)\rightarrow M$, with $\Omega=\lim_{n}s_n$, such
that $F(\dot{\tau})<1$ and $\tau(s_{n})=y_n$. Moreover, the
function $b_{\tau}^+$ satisfies: $b_{\tau}^+\leq \phi$ (as
$d_{(s_{n+1},y_{n+1})}^+\leq \phi$), $b_{c}^+\leq b_{\tau}^+$ for
all $b_{c}^{+}\in {\cal F'}$ (as $b_{\tau}^+\geq
d_{(t_{n},x_{n})}^+$ for all $n$) and $b_{\tau}^+(x_0)\geq
b_{c}^+(x_0)\geq t_0$ (where $b_{c}^+\in {\cal F'}$). Hence,
$b_{\tau}^+\in {\cal F}$ is an upper bound of ${\cal F'}$, in
contradiction with our initial hypothesis. \cvd

\smallskip

Next, the required result on compactification is the following:
\begin{theorem}\label{agnadida1} The (forward) Busemann completion $M_B^+$ of $(M,d)$ endowed
with the chronological topology is a {\cambios sequentially} compact
topological space.
\end{theorem}
{\it Proof.}
Consider any sequence $\sigma=\{[b^+_{c_{n}}]\}$ in {\cambios $M_{B}^+$}. Pick a
point $x_{0}\in M$ and, with no loss of generality, assume that
$b^+_{c_{n}}(x_0)=0$.
By Lemma \ref{l41}, there exists a subsequence
$\{b^+_{c_{n_{k}}}\}_{k}$ with some pointwise limit $\phi\in {\cal
L}_{1}(M,d)$ (in particular, $\phi(x_{0})=0$). Putting $t_0=0$ in
Lemma \ref{ld} (iii),  $-d^+_{x_0}\leq \phi$
and the set ${\cal F}$ defined in (\ref{equf}) is non-empty. From
Lemma \ref{legg}, ${\cal F}$ must contain some maximal function
$b^+_{c}$. This function is a chr. limit of
$\{b^+_{c_{n_{k}}}\}_{k}$ according to (\ref{ecu1}). In fact, let
$h\in B^+(M)$ with $b^+_{c}\leq h\leq \phi$. Necessarily, $h\in
{\cal F}$ and by the maximality of $b^+_{c}$, $h=b^+_{c}$, i.e.
$b^+_{c}\in \hat{L}(\{b^+_{c_{n_{k}}}\}_{k})$. Therefore,
$[b^+_{c}]$ is also a limit of
$\{[b^+_{c_{n_{k}}}]\}_{k}$, as required. \cvd

\subsubsection{The backward Busemann completion}\label{s525}

As always in the non-symmetric case, we can introduce ``backward''
definitions and results, just by using the reversed elements $\F$,
$\rd$ instead of $F$, $d$. However, there are some conventions
(especially about signs) which suggest a different choice for the
backward elements. These conventions will be very useful for the
relations between the forward and backward Busemann boundaries
provided by the c-boundary (see Section \ref{s6}).

So, for the construction of the space $M_B^-$, we start with the
map $j^-:x\mapsto +d(x,\cdot)$ instead of the map $j^+:x\mapsto
-d(\cdot,x)$. The space $C^-(M)$ is the set of all the piecewise
smooth curves $c:[\alpha,\overline{\Omega})\rightarrow M$,
$\overline{\Omega}\leq \infty$, such that $\F (\dot c)\leq 1$.
Here, the endpoint of the interval has been written $\overline{\Omega}$, and we will put typically $\overline{\Omega}=-\Omega$. Now, the {\em backward Busemann
function} is defined as $b_c^-(\cdot)=\lim_{s\rightarrow
\overline{\Omega}}(-s+d(c(s),\cdot))\in \R\cup \{-\infty\}$, and
$B^-(M)$ denotes the space of these functions which are finite.
The {\em backward chronological topology} on $B^-(M)$ is
introduced by means of the following limit operator $\check{L}$,
which replaces $\hat L$ in (\ref{ecu1}):
\[
f\in \check{L}(\sigma)\iff \left\{\begin{array}{l} (a) f\geq
\limsup_{n}f_{n} \hbox{ and}\\ (b) \forall g\in B^{-}(M) \hbox{
with } f\geq g\geq \liminf_{n} f_{n}, \hbox{ it is } g=f.
\end{array}\right.
\]

Recall that the
corresponding notions of Gromov  and Busemann
completions, according to this alternative conventions, are totally
equivalent, and so are the corresponding results.

\subsection{Chronological topology vs Gromov topology}\label{s53}

The Busemann completion and its different parts are regarded as
topological spaces with the chr. topology. However, we have
emphasized that $M_B^+\subset M_G^+$ as a point set and, in fact,
$j_B^+:M_C^+\rightarrow M_B^+$ was just $j^+:M_C^+\rightarrow
M_G^+$ with the co-domain restricted at the point set level. We
have already analyzed the relation between the Cauchy topology on
$M_{C}^{+}$ (induced from the backwards balls of the
quasi-distance $d_Q$) and both, the Gromov and Busemann topologies
on $j^+(M_C^+) (=j_B^+(M_C^+))$, see Corollaries \ref{ccg2} and
\ref{cembbusf}. Next, we are going to study the relations between
the Gromov and Busemann topologies when the former is restricted
to $M_B^+$.


\subsubsection{The inclusion of $M_B^+$ in $M_G^+$ from the topological viewpoint}

 The following result shows that the chr. topology is coarser than
the Gromov one. As the latter is coarser than the Cauchy topology
on $j^+(M_C^+)$, the chr. topology is the coarsest one of the
three topologies.

\begin{proposition}\label{icont}
Consider a sequence $\{f_n\}\subset B^+(M)$ which converges
pointwise to a function $f\in B^+(M)$. Then, $f$ is the unique
chronological limit of $\{f_n\}$.

In particular, the inclusion $i: M_B^+\rightarrow M_G^+$ satisfies
that $i^{-1}: i(M_B^+) \rightarrow M_B^+$ is continuous.
\end{proposition}
{\it Proof.}
For the first assertion observe that, from the pointwise
convergence of $\{f_n\}$, we have $\limsup f_n=\liminf f_n=f$. In
particular, $\sigma=\{f_{n}\}$ and $f$ lie under the hypotheses
of
Proposition \ref{ppo}. Thus, $f$ is the
unique limit of $\sigma$ in the chr. topology.

{\cambios For the last one, as $M_G^+$ is metrizable (Theorem
\ref{tgromovf}) the continuity of $i^{-1}$ is characterized by
sequences, so let $\{[f_n]\}\subset M^+_B$ which converges
pointwise to $[f]\in M^+_B$.}  With no loss of generality, we
can assume  that {\cambios $\{f_n\}$ converges  pointwise to $f$. Then,
$f$ is the (unique) chr. limit of $\{f_n\}$ from the first part}, and the required chr. convergence of $\{[f_n]\}$ to $[f]$
follows from the continuity of the natural projection $\Pi:
B^+(M)\rightarrow M^+_B$.
\cvd

\begin{remark}\label{rjb}
{\em The differences between the
Busemann and Gromov topologies at the point set level appeared
clearly even in the Riemannian case, see Remark \ref{rgrape}. It
is easy to check that, in these cases, there is also a difference
at the topological level. In fact, in Figure \ref{peine} (A) the
sequence $\{[-\dri_{x_n}]\}_{n}$ converges to both,
$[-\dri_{(0,0)}]$ and $[-\dri_{(0,1)}]$ (recall Fig. \ref{peine}
(B)). As we will see, there are differences at the point set level
occur if and only if the Busemann topology is different (strictly
coarser) than the induced by Gromov one in $M_B^+$ (see Remark
\ref{rgrape2}).

Therefore, in general the chr. topology is strictly coarser than
Gromov one. Recall also that Example \ref{ex43} also showed that
the Gromov topology (and then, the Busemann one) is strictly
coarser than the Cauchy topology.}
\end{remark}

\subsubsection{The case of separating topologies}

Most of the difficulties of the chr. topology come from the fact
that the limit operator $\hat L$ may be not of the first order. We
will see next that a simple hypotheses prevents this case and
ensures automatically quite a few of nice properties. These
hypotheses are adapted from a natural one for the c-boundary of
spacetimes, which was studied in  \cite[Defn. 3.42, Prop.
3.44]{FHS0}. Recall that the open subsets $I^{+}(p_0,B^{+}(M))$
were defined in (\ref{ee}). 
\begin{definition}\label{dseparating} The (forward) chr. topology on $B^+(M)$ is {\em
separating} if the following property holds: for any two functions
$f,f'\in B^+(M)$, with $f'\leq f$, $f'\neq f$, there exists some
$p_0=(t_0,x_0)\in \R\times M$ which {\em separates} them, i.e.,
$p_0$ satisfies $f'\not\in \overline{I^{+}(p_0,B^{+}(M))}$ but
$f'(x_0)<t_0<f(x_0)$ (and thus, $f\in I^{+}(p_0,B^{+}(M))$).
\end{definition}
\begin{proposition}\label{fg} If $B^+(M)$ (endowed with the chr. topology) is
separating, then $B^+(M)$ {\cambios is  second countable and $\hat L$ is
of first order}. Therefore, $M_{B}^{+}$ is also second
countable.
\end{proposition}
{\it Proof.} First, let us prove that the following family of
subsets in $B^{+}(M)$ is a (countable) topological basis:
\begin{equation}\label{ebasetop}
\Omega=\{I^{+}(p_{n},B^{+}(M))\cap \overline{\cup_{i=1}^{k} I^{+}(
p_{j_{i}},B^{+}(M))}^{c}: n, j_{i}\in {\mathbb N}\},
\end{equation}
where $D=\{p_{n}\}_{n}$ is a dense countable subset of $\R\times
M$. Consider $f'\in B^{+}(M)$, choose a chain $\{d^+_{p'_{n}}\}$
converging pointwise to $f'$ with $\{p'_{n}\}\subset D$ {\cambios
(recall Lemma \ref{ld}(ii)) and construct a basis $\{U_n\}\subset
\Omega$ of open neighborhoods of $f'$ as follows.  Let
$\{s_{i}=p_{n_i}\}_i\subset D$  be the (countable) subset composed
by those points in $D$ which separate $f'$ and $f$, for some $f\in
B^{+}(M)$ such that $f'\leq f$, $f'\neq f$. Recall about this
subset that, from the separation property assumed in the
hypotheses, for any $f$ as above there exists some
$p_0=(t_{p_0},x_{p_0})\in \R\times M$ which separates $f'$ and
$f$. Thus,  any $s=(t_{s},x_{s})\in D$ satisfying
$d(x_{p_0},x_s)<t_s-t_{p_0}$ and $t_{s}<f(x_{s})$ provides an
element of $\{s_i\}$ associated to $f$ (recall that necessarily
$d^+_{p_0}<d^+_s$, then $I^+(s,B^+(M))\subset I^+(p_0,B^+(M))$ and
this inclusion holds for the closures, i.e., $f'\not\in
\overline{I^+(s,B^+(M))}$). The density of $D$ allows to find
(infinitely many) such a $s$ with $d(x_{p_0},x_s)$ small enough,
and the full $\{s_i\}$ is constructed by taking all such $s$ for
all such functions $f$.
 Now, define each $U_{n}\in \Omega$ as
$U_{n}:=I^{+}(p'_{n},B^{+}(M))\cap A_n$, where
$A_n=\overline{\cup_{i=1}^{n}I^{+}(s_{i},B^{+}(M))}^{c}$. Since
$f'\not\in \overline{I^{+}(s_i,B^{+}(M))}$ for all $i$, the open
subsets $U_n$ contain $f'$.}

{\cambios
 To check that  $\{U_n\}$ constitute the required basis around
$f'$, let us prove first that any sequence $\{g_n\}$ such that
$g_n\in U_m$ for all $m$ and all $n\geq n(m)$, satisfies $f'\in
\hat{L}(\{g_n\})$. In fact, to check the first condition in
(\ref{ecu1}), recall that $g_n\in I^+(p'_m,B^+(M))$, and so,
$g_n>d^+_{p'_m}$. As $d^+_{p'_m}\rightarrow f'$ pointwise, the
inequality $f'\leq \liminf_{n}g_n$ holds. For the second condition
in (\ref{ecu1}), if $f'\leq f\leq\limsup_{n}g_n$, $f'\neq f$, then
$f\in I^{+}(s_{i_{0}},B^{+}(M))$ for some $s_{i_0}\in \{s_i\}$
and, thus, $t_{s_{i_{0}}}<f(x_{s_{i_{0}}})$. As $f\leq
\limsup_{n}g_n$, then $t_{s_{i_{0}}}<g_n(x_{s_{i_{0}}})$, i.e.,
$g_n\in I^+(s_{i_{0}},B^+(M))$, for infinitely many $n$, in
contradiction with the hypothesis $g_n\in U_m$.  Now, it is clear
that if some neighborhood $U$ of $f'$ satisfied $U_{n}\not\subset
U$ for all $n$, the sequence $\{g_n\}$ obtained by choosing some
$g_n\in U_{n}\setminus U$ for all $n$ would not converge to $f'$
with the chr. topology but $f'\in \hat L(\{g_n\})$, an absurd (see
Remark \ref{rsequential}). }

{\cambios The first order property of $\hat L$ is also
straightforward now, because if $\sigma=\{g_{n}\}$ converges to
$f'$ with the chr. topology, necessarily $g_{n}\in U_{m}$ for all
$m$ and all $n\geq n(m)$, and then $f'\in \hat{L}(\{g_n\})$. }

For the last assertion,  just use Proposition \ref{gh} and project
the basis (\ref{ebasetop}). \cvd

\begin{remark}\label{rduda} {\em As we will see, if $B^+(M)$ is
Hausdorff then it is separating (just apply Remark \ref{rULS} and
Proposition \ref{key}), but the converse does not hold.  In the
framework of spacetimes, it is easy to find examples of causal
completions such that its (chr.) topology is not separating and,
even more, non-first countable, see \cite[Example 3.43, Remark
3.40]{FHS0}. So, it is conceivable that these properties may occur
for $B^+(M)$ (and, then, be transmitted to $M_B^+$), even though
we do not have any example where (\ref{ebasetop}) is not a basis
for the chr. topology.}\end{remark}

\subsubsection{Topologies with unique limits of sequences}

As we have seen, $M_B^+$ may be non-Hausdorff and this is a clear
difference with $M_G^+$. Our aim will be to show that only in this
case {\cambios there are} differences between both topologies. However,
we will study the non-Hausdorff character from a more general
viewpoint (compare with \cite{Douwen}).

\begin{definition} Let $X$ be a set.
A topology on  $X$ is  {\em unique limit for sequences} (ULS) if
any sequence in $X$ has at most one limit. A topology $\tau_L$
associated to a limit operator $L$ on  $X$ is {\em unique limit
for $L$} (ULL) if, for any sequence $\sigma$ in $X$, $L(\sigma)$
contains at most one element.
\end{definition}
\begin{remark}\label{rULS}{\em  Clearly, Hausdorff $\Rightarrow$ ULS $\Rightarrow $ ULL
(recall Remark \ref{rsequential}) and, for {\cambios first} countable
spaces, ULS $\Rightarrow$ Hausdorff. }\end{remark} {\cambios Now, let us
consider again $B^+(M)$ with the chr. topology ---always regarded
as the topology associated to $\hat L$.}

\begin{lemma}\label{jh} Assume that $B^+(M)$ is {\cambios ULL}. If a sequence
$\{b^+_{c_{n}}\}\subset B^{+}(M)$ converges pointwise to some
Lipschitz function $\phi\in {\cal L}_{1}(M,d)$ and $b^+_{c}\in
\hat{L}(\{b^+_{c_{n}}\})$, then $\phi=b^+_{c}$.
\end{lemma}
{\it Proof.}  Since $b^+_{c}\in
\hat{L}(\{b^{+}_{c_{n}}\})$, necessarily $b^+_{c}\leq\phi$. 
Assuming by contradiction that $\phi\neq b^+_{c}$, there is some
$x_{0}\in M$ and $t_{0}\in\R$ with
$b^+_{c}(x_{0})<t_{0}<\phi(x_{0})$. Thus, for $p_{0}=(t_0,x_0)$,
necessarily
$d_{p_0}^+<\phi$ (recall Lemma \ref{ld} (iii)), and so
$d_{p_0}^+\in {\cal F}$, i.e. the set ${\cal F}$ constructed in
(\ref{equf}) is non-empty. Then, Lemma \ref{legg} ensures that
${\cal F}$ must contain a maximal function $b^+_{\tau}$ {\cambios and,
even more}, $b^+_{\tau}\in \hat{L}(\{b^{+}_{c_{n}}\})$ (for any
$h\in B^{+}(M)$ with $b^+_{\tau}\leq h\leq \phi$, we have $h\in
{\cal F}$; since $b^+_{\tau}$ is maximal, necessarily
$h=b^+_{\tau}$, and thus, $b^+_{\tau}\in
\hat{L}(\{b^{+}_{c_{n}}\})$). But since $b^+_{\tau}(x_{0})\geq
t_0>b^+_{c}(x_{0})$, we know $b^+_{\tau}\neq b^+_{c}$, in
contradiction to the {\cambios ULL} character of $B^+(M)$. \cvd

\begin{proposition}\label{key} If $B^{+}(M)$ is
ULL, then it is separating. So, $\hat L$ is of first order and both,
the chr. topologies on $B^+(M)$ and $M_{B}^{+}$, are {\cambios ULS,
second countable and Hausdorff.}
\end{proposition}
{\it Proof.} To prove that $B^{+}(M)$ is separating, let $f',f\in
B^+(M)$, with $f'\leq f$ and $f'(x_0)<t_0<f(x_0)$ for some
$p_0=(t_0,x_0)$, and let us prove that $f'\not\in
\overline{I^{+}(p_0,B^{+}(M))}$ (recall (\ref{ee})).
Note that it suffices to
check that the subset $\{g\in B^+(M): t_0\leq g(x_0)\}$ (which
does not contain $f'$), is chr. closed. Reasoning  by
contradiction,
assume that there exists some sequence $\sigma=\{f_n\}$, with
$t_0\leq f_n(x_0)$ for all $n$, and $h\in\hat{L}(\sigma)$, such
that $h(x_0)<t_0$. As $M_G^+$ is  compact, {\cambios there exist
constants $t_n\in\R$ and a subsequence $\{f_{n_k}+t_k\}$
converging pointwise to some Lipschitz function $\phi'$. It is not
a restriction to suppose that this subsequence satisfies
$t_k\rightarrow t_\infty$ for some $t_\infty \in \R$.  In fact,
otherwise, for some subsequence $\{t_{k_{l}}\}_l$ we have either
$t_{k_{l}}\rightarrow\infty$, and so, as $t_0\leq f_n(x_0)$,
$\phi'\equiv\infty$ (an absurd), or $t_{k_{l}}\rightarrow -\infty$
and, then, $f_{n_{k_{l}}}(x_0) \rightarrow\infty$, which, by
Remark \ref{rrr1}(2), implies liminf$\{f_{n_{k_{l}}}\}=\infty$, in
contradiction with the maximality of $h$ in $\limsup \{f_n\}$
required for the elements of $\hat{L}(\sigma)$.} Therefore,
$\overline{\sigma}=\{f_{n_k}\}$ converges pointwise to
$\phi=\phi'-t_\infty$. In particular, $h(x_0)< t_0\leq \phi(x_0)$.
Since $h\in \hat{L}(\overline{\sigma})$, this  contradicts Lemma
\ref{jh}.

{\cambios For the last assertions recall that, as $B^+(M)$ is
separating, Proposition \ref{fg} ensures that the chr. topology is
second countable and $\hat{L}$ is of first order. Moreover, as
$B^+(M)$ is ULL, the chr. topology is ULS and, so, Hausdorff
(recall Remark \ref{rULS}). } \cvd

\subsubsection{Main results}

The following theorems show the precise relation between the
Busemann and Gromov completions. The following technical lemma is
required first.
\begin{lemma}\label{lee} $B^{+}(M)$ is ULL if and only if $M_{B}^{+}$ is Hausdorff.
\end{lemma}
{\it Proof.} The implication to the right  is included in
Proposition \ref{key}. For the implication to the left, assume
that $B^{+}(M)$ is not ULL. {\cambios Then, $\hat{L}(\sigma)$ contains
two distinct  elements for some sequence $\sigma\subset B^{+}(M)$.
From property (b) in (\ref{ecu1}), the difference between these
elements is not a constant. So, they define two different classes
of $M_{B}^{+}$, and thus, $M_{B}^{+}$ is not Hausdorff}. \cvd
\begin{theorem}\label{below}
The following statements are equivalent:
\begin{itemize}
\item[(a)] No sequence in $M_B^+$ converges to more than one point
in $\partial_{B}^{+}M$ with the forward chronological topology.
For any sequence $\sigma\subset M_{B}^{+}$, $L(\sigma)$ contains
at most one point in $\partial_{B}^{+}M$. \item[(b)] $M_B^+$ is
Hausdorff with the forward chronological topology. \item[(c)]
$i:M_B^+\rightarrow M_G^+$ is continuous (and thus, an embedding).
\item[(d)] $M_B^+=M_G^+$ (both, as a point set and topologically).
\item[(e)] {\cambios $M_B^+$ is equal to $M_G^+$ as a
point set}.
\end{itemize}
\end{theorem}
{\it Proof.} $(a)\Rightarrow (b)$.  Suppose by contradiction that
$M_B^+$ is not Hausdorff, and thus, $B^{+}(M)$ is not ULL (recall
Lemma \ref{lee}). {\cambios So, there exists a sequence
$\{b^+_{c_n}\}_n$ such that $\hat L (\{b^+_{c_n}\})$ contains two
distinct elements of $B^{+}(M)$}. By hypothesis (a), one of
these functions is of the form $\Omega-d^+_x$ with $x\in M$. From
Proposition \ref{lemmabusemanM}, $\{b^+_{c_n}\}_n$ converges
pointwise to $\Omega-d^+_x$. But, then, Proposition \ref{icont}
ensures that this is the unique limit, a contradiction.

$(b)\Rightarrow (c)$. Since $M_{B}^{+}$ is a sequential space {\cambios the continuity of $i$ is characterized by means of sequences
(recall Remark \ref{rsequential}); by Proposition
\ref{icont}, $i$ will be then also an embedding. So, assume that
$\{[b^+_{c_{n}}]\}$ converges with the forward chronological
topology to $[b^+_{c}]$. It is not a restriction to assume that
$\{b^+_{c_{n}}\}$ converges with the forward chronological
topology to $b^+_{c}$. As $M_{B}^{+}$ is Hausdorff, $B^{+}(M)$ is
ULL (Lemma \ref{lee}), and thus, $\hat{L}$ is of first order
(Proposition \ref{key}). So, we can also assume that $b^+_{c}\in
\hat{L}(\{b^+_{c_{n}}\})$ and, reasoning  by contradiction, that
$\{b^+_{c_{n}}\}$ does not converge pointwise to $b^+_{c}$.
Moreover, there exists a subsequence $\{b^+_{c_{n_{k}}}\}_{k}$
converging pointwise to some $\phi\neq b^+_{c}$ (reason as in
Proposition 5.37, namely,  for any $x_0\in M$, the sequence
$\{|b^+_{c_{n_{k}}}(x_0)|\}_{k}$ must be bounded because,
otherwise, a contradiction with $\hat{L}(\{b^+_{c_{n}}\})\neq
\emptyset$ would appear).
This is in contradiction with Lemma \ref{jh}, since
$b^+_{c}$ belongs also to $\hat{L}(\{b^+_{c_{n_{k}}}\}_{k})$. }


$(c)\Rightarrow (d)$. {\cambios Any point $[f]\in M_G^+$ is  the
pointwise limit of some sequence $\{j^{+}(x_n)\}\subset M_B^+$,
with $\{x_n\}\subset M$. As $M_B^+$ is sequentially compact, up to
a subsequence, $\{j^{+}(x_n)\}$ converges with the forward
chronological topology to some $[\overline{f}]\in M_B^+$. By the
continuity of $i$, $\{j^{+}(x_n)\}$ also converges pointwise to
$[\overline{f}]$, and so, $[\overline{f}]=[f]$. So, $i$ is
onto and, under our hypotheses, a homeomorphism. }

$(d)\Rightarrow (e)$. Trivial.

$(e)\Rightarrow (a)$. By contradiction, assume that
$\{[b^{+}_{c_{n}}]\}\subset M_{B}^{+}$ converges with the chr.
topology to distinct elements $[b_{c}^{+}]\neq [b_{c'}^{+}]$ of
$\partial_B^+ M$. By the compactness of $M^{+}_{G}$ and hypothesis
(e), there exists some subsequence $\{b^{+}_{c_{n_{k}}}\}_{k}$
converging pointwise to some $[f]\in M^{+}_{G}=M_{B}^{+}$. From
Proposition \ref{icont}, $[f]$ must be the unique chronological
limit of $\{[b^{+}_{c_{n_{k}}}]\}_{k}$, and thus of
$\{[b^{+}_{c_{n}}]\}$, in contradiction with the initial
assumption. \cvd

Recall that Gromov and Eberlein and O'Neill's boundaries coincide,
both as a point set and topologically, when $(M,g)$ is a Hadamard
 manifold \cite{BGS,BH}. So, Theorem \ref{below} also ensures
 the equality with the Busemann boundary.
\begin{corollary}\label{new} Let $(M,g)$ be a Hadamard manifold.
Then the Busemann and Gromov compactifications (as well as
Eberlein and O'Neill's one) coincide.
\end{corollary}
{\cambios {\em Proof.} As Eberlein and O'Neill's compactification is
carried out by means of Busemann functions of rays, then it is
always included in $M_B$ as a point set. So, as that
compactification agrees with Gromov's one in a Hadamard manifold,
$M_G=M_B$ as point sets and, from Theorem \ref{below}, also
topologically.} \cvd

\begin{remark}\label{rgrape2} {\em
As emphasized in
Remark \ref{rgrape} and Figure \ref{escalera2}, it is not
difficult to find examples with $M_B\neq M_G$ as a point set
---and, as a consequence of  Theorem \ref{below}, where $i:M_B^+\rightarrow
M_G^+$ is neither an embedding. This happens when there exists
some point in $M_{G}$  that  cannot be reached from $M$ by means
of an {\em asymptotically ray-like curve}, i.e. a curve with
finite Busemann function. Perhaps the most dramatic difference
happens when some points in $\partial_GM$ cannot be reached as
limits of any type of curves in $M$ (recall Figure \ref{peine} (A)).
But, even when they can, it is not difficult to construct examples
where the Gromov completion differs from Busemann one.

The example depicted in Figure \ref{figchim1}  may be especially
illustrating. Here, if one modifies the usual metric of a cylinder
by putting  the depicted ``chimney'' (Figure \ref{figchim1} (A)),
then a point appears for both, the Gromov and Busemann boundaries
(as well as for the Cauchy one). If two such chimneys are included
(Figure \ref{figchim1} (B)), two points appear for $\partial_BM$
(and $\partial_CM$), but {\em a continuum of points} (a
segment)
appears for
$\partial_GM$. The two points in $\partial_BM$ are not Hausdorff
separated for the chr. topology. This reflects the fact that
$\partial_B M$ is a {\em compactification}, not only a {\em
completion} (as $\partial_CM$ was). The ``surprising'' appearance
of a continuum of points for $\partial_G M$ is a consequence of
the Hausdorffness of the Gromov compactification.

It is also worth pointing out that, in the Figure
\ref{figchim1} (B), the induced topology in $\partial_BM$ from
$M_B$ is  the discrete one {\cambios (this is a consequence of the fact
that the chr. topology is $T_1$ and $\partial_BM$ contains only
two points)}. So $\partial_BM$ is Hausdorff, even though the two
points of $\partial_BM$ are not Hausdorff related as points of
$M_B$. Analogously, local compactness is not satisfied (neither in
the case (A) nor in (B) of Figure \ref{figchim1}) by the points of
$\partial_CM$ regarded as points of $M_C$, but it is satisfied
 in $\partial_CM$ regarded as a topological space. }
\end{remark}

\begin{figure}
\centering \ifpdf
  \setlength{\unitlength}{1bp}%
  \begin{picture}(503.47, 203.87)(65,0)
  \put(0,0){\includegraphics{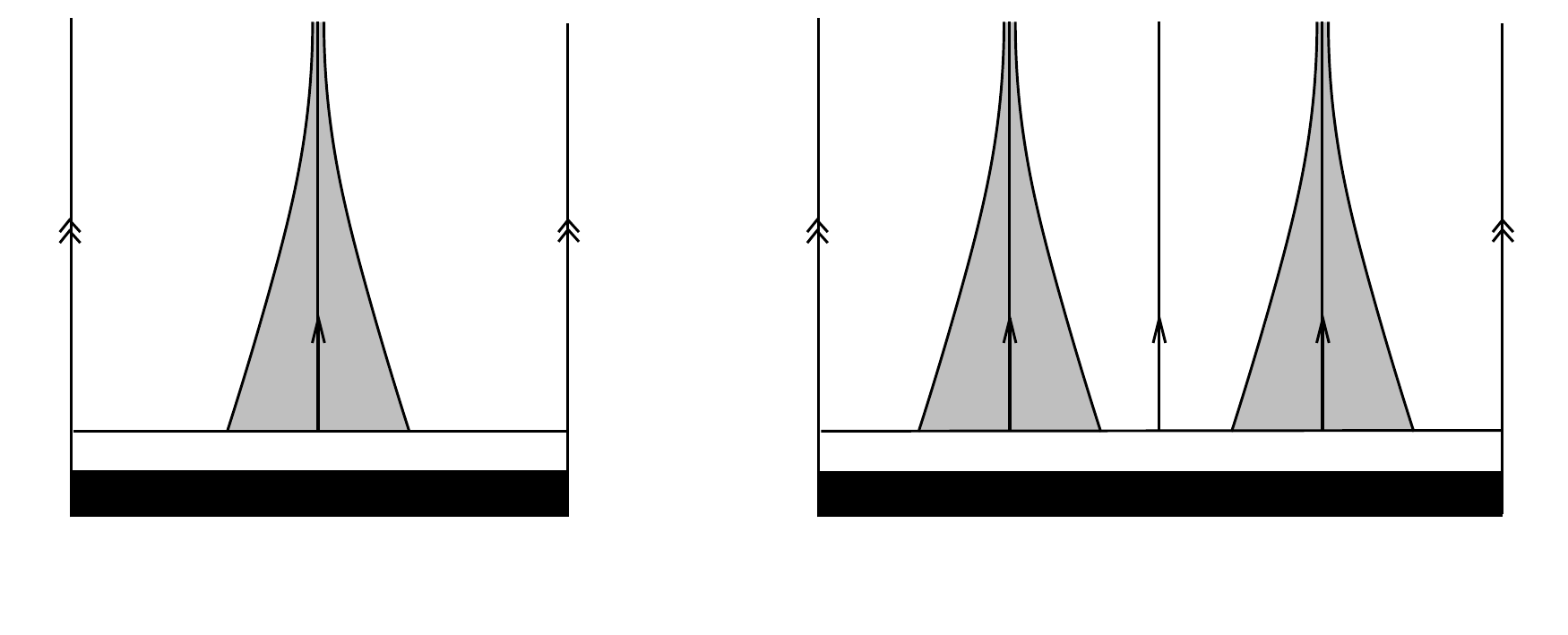}}
  \put(88.05,79.81){\fontsize{8.54}{10.24}\selectfont $c_0$}
  \put(307.73,75.10){\fontsize{8.54}{10.24}\selectfont $c_{-3}$}
  \put(428.56,75.10){\fontsize{8.54}{10.24}\selectfont $c_3$}
  \put(375.21,82.95){\fontsize{8.54}{10.24}\selectfont $c_a$}
  \put(5.67,63.33){\fontsize{8.54}{10.24}\selectfont $-6$}
  \put(246.53,63.33){\fontsize{8.54}{10.24}\selectfont $-6$}
  \put(183.77,63.33){\fontsize{8.54}{10.24}\selectfont $6$}
  \put(485.04,63.33){\fontsize{8.54}{10.24}\selectfont $6$}
  \put(96.68,56.27){\fontsize{8.54}{10.24}\selectfont $0$}
  \put(317.93,57.06){\fontsize{8.54}{10.24}\selectfont $-3$}
  \put(420.71,57.06){\fontsize{8.54}{10.24}\selectfont $3$}
  \put(368.14,57.06){\fontsize{8.54}{10.24}\selectfont $a$}
  \put(88.92,8.73){\fontsize{14.23}{17.07}\selectfont (A)}
  \put(365.79,8.73){\fontsize{14.23}{17.07}\selectfont (B)}
  \end{picture}%
\else
  \setlength{\unitlength}{1bp}%
  \begin{picture}(503.47, 203.87)(65,0)
  \put(0,0){\includegraphics{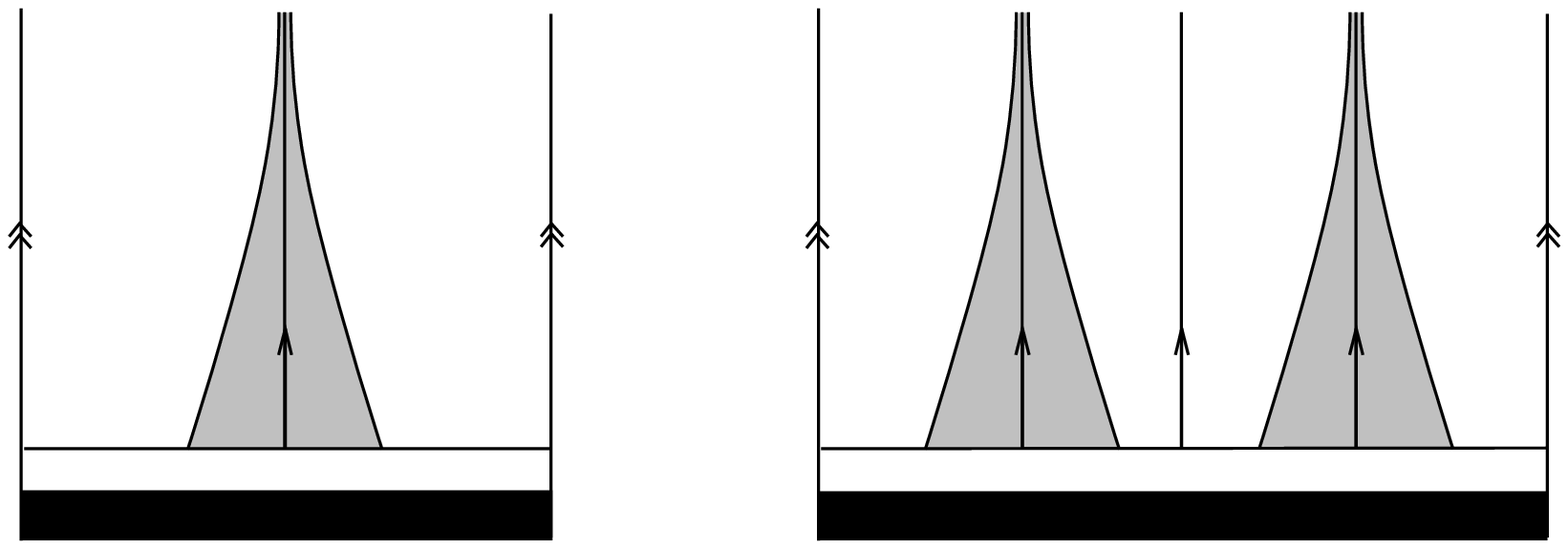}}
  \put(88.05,79.81){\fontsize{8.54}{10.24}\selectfont $c_0$}
  \put(307.73,75.10){\fontsize{8.54}{10.24}\selectfont $c_{-3}$}
  \put(428.56,75.10){\fontsize{8.54}{10.24}\selectfont $c_3$}
  \put(375.21,82.95){\fontsize{8.54}{10.24}\selectfont $c_a$}
  \put(5.67,63.33){\fontsize{8.54}{10.24}\selectfont $-6$}
  \put(246.53,63.33){\fontsize{8.54}{10.24}\selectfont $-6$}
  \put(183.77,63.33){\fontsize{8.54}{10.24}\selectfont $6$}
  \put(485.04,63.33){\fontsize{8.54}{10.24}\selectfont $6$}
  \put(96.68,56.27){\fontsize{8.54}{10.24}\selectfont $0$}
  \put(317.93,57.06){\fontsize{8.54}{10.24}\selectfont $-3$}
  \put(420.71,57.06){\fontsize{8.54}{10.24}\selectfont $3$}
  \put(368.14,57.06){\fontsize{8.54}{10.24}\selectfont $a$}
  \put(88.92,8.73){\fontsize{14.23}{17.07}\selectfont (A)}
  \put(365.79,8.73){\fontsize{14.23}{17.07}\selectfont (B)}
  \end{picture}%
\fi \caption{\label{figchim1} Non-Hausdorff $\partial_BM$ and
{\em surprising} points in $\partial_GM$ (Remark \ref{rgrape2}).
\newline Starting at a strip on the plane, consider the half
cylinder $M=[-6,6]\times (-1,\infty)/\sim$ where $\sim$ identifies
each $(-6,y)$ with $(6,y)$, for all $y>-1$ (the black regions are
irrelevant for our problem, and will not be taken into account).
Put a Riemannian metric (which can be chosen conformal to the
usual Euclidean one) in such a way that: (i) out of the
``chimneys'' (the regions drawn in grey) the metric is the
Euclidean one, and (ii) the central curves of the chimneys ($c_0$
in the case (A) and $c_{-3},c_3$ in the case (B)) have finite
length.\newline $\hbox{$\quad$}$ Case (A): $M_B$ and $M_G$
coincide {\cambios (as, obviously, no sequence in $M$ can converge
to two distinct points of  $\partial_BM$). Moreover, the
corresponding boundary is naturally identified with
$\{x_0\}$, where $x_0$ is the Cauchy boundary point associated to
the curve $c_0$.
\newline $\hbox{$\quad$}$ Case (B): The Cauchy and Busemann
boundaries coincide as point sets and are identifiable with
$\{x_{-3}, x_3\}$, where $x_{-3},x_3$ are the Cauchy boundary
points associated to $c_{-3}$ and $c_3$ resp. Nevertheless, the
Gromov boundary include also a point for any curve of the form
$c_a(t)=(a,t)$ with $a\in [-3,3]$, i.e., the Gromov boundary
naturally includes a ``segment'' which connects $x_{-3}$ and
$x_3$.}}
\end{figure}

\subsection{Proof of Theorem \ref{tresumen}}\label{5.4}

Now, we are in conditions to
summarize the proof of Theorem \ref{tresumen}.

\smallskip

{\it Proof of Theor. \ref{tresumen}.} {\em Point 1.} It follows
from Theorem \ref{tgromovf}.

{\em Point 2.} The {\cambios sequential} compactness is proved in
Theorem \ref{agnadida1}. For the $T_1$ character of the chr.
topology, recall Remark \ref{rsequential}. The non-Hausdorff
possibility (as well as the possible non-continuity of the
inclusion $i$ in the second paragraph), are illustrated by Figure
\ref{peine} (Remark \ref{rjb}). The continuity of $i^{-1}$ in the
second paragraph is proved in Proposition \ref{icont}. The third
paragraph  is included in Theorem \ref{below}.

{\em Point 3.} (First paragraph). The character of quasi-distance
is ensured by Proposition \ref{ppp1'} and the $T_0$ property by
Theorem \ref{tcompltcompl}.

{\em Point 3.} (A) The first paragraph follows from Proposition \ref{pembf}.
 For the second paragraph, the local compactness of $M^+_{CG}$
follows from the compactness of $M^+_G$ and the Heine-Borel
property is proved in Corollary \ref{ccg2} (2). The third
paragraph collects  Corollary \ref{ccg2} (3).

{\em Point 3.} (B)  follows from Corollary \ref{cembbusf}.

{\em Point 3.} (Last paragraph). {\cambios Direct consequence of
Theorem \ref{tcompltcompl}, Proposition \ref{p329}, Theorem
\ref{tgromovf} and Corollary \ref{cembbusf}  (for the last
assertion, use also Theorem \ref{below}).}
\cvd

\section{C-boundary of standard stationary spacetimes}\label{s6}

Next, we are going to use the structures introduced in previous
sections to describe the c-boundary of standard stationary
spacetimes. To this aim, we will consider spacetimes as in
(\ref{ggg}) with the choice of conformal factor $\Lambda=1$, i.e.
\begin{equation}\label{h'}
V=\R\times M\qquad\hbox{and}\qquad g=-dt^2+\pi^*\omega\otimes
dt+dt\otimes \pi^*\omega +\pi^*g_{0}.
\end{equation}
The natural projection $t: \R\times M \rightarrow \R$ will be called the {\em (standard) time},
---and we will use the letter $t$ with some obvious abuse of notation.
From this metric, we consider the
Fermat metrics $F^{\pm}=\sqrt{g_0 +\omega^2} \pm \omega$ as in
(\ref{h}). We will denote with subscripts $+$ (resp. $-$) the
objects associated to the metric $F^{+}$ (resp. $F^{-}$).
 So, the associated generalized distance will be denoted $\dpl$ and its reverse distance $\dm$.
$(M,d^{\pm})$ will be {\cambios the {\em generalized metric spaces}}
associated to the standard stationary spacetime.

\subsection{Chronological relations and Lipschitzian functions}

In this subsection we are going to establish the chronological
relations of $(V,g)$ in terms of the Finsler metrics $F^{\pm}$.
The arguments are inspired by the seminal work of Harris for the
standard static case \cite{H}.
\begin{proposition}\label{k}
Consider two points $(t_{0},x_{0}),(t_{1},x_{1})\in V$. Then:
$$(t_{0},x_{0})\ll (t_{1},x_{1})\iff \dpl(x_{0},x_{1}) < t_{1}-t_{0}.$$
\end{proposition}
{\it Proof.} If $(t_{0},x_{0})\ll (t_{1},x_{1})$ then there exists
a future-directed timelike curve $\gamma(s)=(t(s),c(s))$ such that
$\gamma(0)=(t_{0},x_{0})$, $\gamma(1)=(t_{1},x_{1})$. Since
$\gamma$ is future-directed timelike, necessarily
$\dot{t}(s)>\sqrt{g_{0}(\dot{c}(s),\dot{c}(s))+\omega(\dot{c}(s))^{2}}+\omega(\dot{c}(s))=F^{+}(\dot{c}(s))$
for all $s$. In particular,
\[
t_{1}-t_{0}=\int_{0}^{1}\dot{t}(s)ds>\int_{0}^{1}F^{+}(\dot{c}(s))={\rm
length}_{+}(c)\geq \dpl(x_{0},x_{1}).
\]
Reciprocally, if $t_{1}-t_{0}>\dpl(x_{0},x_{1})$ then there
exists a curve $c(s)$ in $M$ with $c(0)=x_{0}$, $c_{1}(1)=x_{1}$
and constant speed $F^{+}(\dot{c}(s))<t_{1}-t_{0}$. Therefore, the curve
$\gamma(s)=(t(s),c(s))$, with $t(s)=(t_{1}-t_{0})s+t_{0}$, is
timelike, future-directed and satisfies $\gamma(0)=(t_{0},x_{0})$,
$\gamma(1)=(t_{1},x_{1})$, as required. \cvd

\smallskip

Notice that we can reparametrize the timelike curves by using the standard time
$t$. Therefore, from the expression of the standard
stationary metric $g$ in (\ref{h'}), we  deduce directly:
\begin{proposition}
A curve $\gamma(t)=(t,c(t))$ (resp. $\gamma(t)=(-t,c(t))$), $t\in I$ (interval),
is future (resp. past) directed timelike if and only if $F^{+}(\dot{c})<1$ (resp.
$F^-(\dot{c})<1$).
\end{proposition}
As a consequence, the past and future sets of $V$ (according to
the definitions in Section \ref{fff'}) can be characterized by
means of Lipschitz functions. More precisely:
\begin{definition}\label{d63} Let $f$ be a Lipschitz
function
 for $\dpl$ (resp. $\dm$). Its {\em past} $P(f)$ (resp. {\em future} $F(f)$) is the following subset of $V=\R \times M$:
\[
P(f):=\{(t,x)\in V: t<f(x)\}\subset V\quad (\hbox{resp.} \,\;
F(f):=\{(t,x)\in V: t>f(x)\}\subset V).
\]
\end{definition}
\begin{proposition}\label{jdl}
A (non-empty)  subset $P$, $P\varsubsetneq V$ (resp. $F$, $F\varsubsetneq V$) is a past (resp. future) set if
and only if $P=P(f)$ (resp. $F=F(f)$) for some Lipschitz function $f$ for $\dpl$ (resp. for $\dm$).
\end{proposition}
{\it Proof.} For the implication to the right, given $P$  define
\begin{equation} \label{eps} f(x):={\rm sup}\{t\in {\mathbb R}:(t,x)\in P\}, \quad \forall x\in
M,\end{equation} and choose two points $(t_0,x_0)\notin P$,
$(t_1,x_1)\in P$. As $\dpl(x,x_1)<t_1-t$ for large $-t$,
necessarily $f(x)>-\infty$. Moreover, $f(x)\leq t_0+\dpl(x_0,x)$
because, otherwise, any $t\in (t_0+\dpl(x_0,x),f(x))$ would yield
the contradiction $(t_0,x_0)\ll (t,x)\in P$. In particular,
$f(x)<\infty$ and (\ref{eps}) defines a function $f:M\rightarrow
\R$. As $P$ is a past set, observe that any pair $(t,x)$ with
$t<f(x)$ belongs to $P$ (if $t<f(x)$, there exists $t'$ such that
$t<t'<f(x)$ and $(t,x)\ll (t',x)\in P$). Furthermore, $f$ is
Lipschitz (so,  we can write $P=P(f)$ according to Defn.
\ref{d63}). In fact, otherwise
 $f(y)>f(x)+\dpl(x,y)$ for some $x, y$. Picking some $t\in (f(x)+\dpl(x,y),f(y))$
we have $(f(x),x)\ll (t,y)\in P$, in contradiction to the facts
that $P$ is a past set and $(f(x),x)\not\in P$ ($P$ must be open).

For the implication to the left, it suffices to prove that
$P(f)=I^{-}[P(f)]$. For the inclusion $\subseteq$, take any
$(t,x)\in P(f)$. Then, $t<f(x)$, and so, there exists $t'$ such
that $t<t'<f(x)$. In particular, $(t,x)\ll (t',x)\in P(f)$, since
$t'-t>0=\dpl(x,x)$. For the inclusion $\supseteq$, suppose that
$(t,x)\in I^{-}[P(f)]$. Then, there exists some point $(t',x')\in
P(f)$ such that $(t,x)\ll (t',x')$, or, equivalently,
$\dpl(x,x')<t'-t$. Assume by contradiction that $(t,x)\not\in
P(f)$. Then, $t\geq f(x)$. This joined to the fact that
$(t',x')\in P(f)$, i.e. $t'<f(x')$, implies
\[
t'-t<f(x')-f(x)\leq \dpl(x,x'),
\]
the latter inequality because $f$ is Lipschitz.\cvd

\begin{convention}
{\em As $V$ is also a past and future set, we put
$P(\infty)=F(-\infty):=V$.
}
\end{convention}

\subsection{Future and past c-boundaries as point sets}\label{mapStatRand}

Next, we are going to construct the partial (future and past)
c-boundaries. Consider a future-directed timelike curve
$\gamma:[\alpha,\Omega)\rightarrow V$ parametrized as
$\gamma(t)=(t,c(t))$. By using that function $t\mapsto
t-\dpl(\cdot,c(t))$ is increasing (recall (\ref{equdecre})), the
past of $\gamma$ can be expressed in the following way:
\begin{equation}\label{eimas}
\begin{array}{lll}
I^{-}[\gamma] & = & \{(t',x')\in V:(t',x')\ll \gamma(t) \hbox{ for some }t\in [\alpha,\Omega) \}\\
& = &\{(t',x')\in V:t'<t-\dpl(x',c(t)) \hbox{ for some }t\in [\alpha,\Omega)\}\\
 &= &\{(t',x')\in V:t'<\lim_{t\rightarrow \Omega}(t-\dpl(x',c(t)))\}\\
 & = &\{(t',x')\in V: t'<b_{c}^{+}(x')\}.
\end{array}
\end{equation}
Therefore, the set of all the IPs, $I^{-}[\gamma]\neq V$, can be
identified with $B^{+}(M)$. When $b^{+}_c\equiv \infty$ then
$I^{-}[\gamma]=V$, and this TIP will be denoted $i^{+}$. Summing
up:
\begin{proposition}\label{caractIP} A subset $P\subset V$ (resp. $P\subsetneq V$) is an IP
if and only if $P=P(b^{+}_{c})$ for some Busemann function
$b^{+}_{c}$ (resp. for some finite Busemann function $b_c^+\in
B^{+}(M)$).
\end{proposition}
Obviously, such a $b^{+}_{c}$ is univocally determined by $P$.
Therefore, the study of the future c-boundary of $(V,g)$ is
reduced to the study of Busemann functions for the Finsler
manifold $(M,F^{+})$. More precisely, choose
any $x_0$ in $M$, and notice that the following map is
bijective,
\begin{equation} \label{ezzz}
\pi^+_{x_0}:B^+(M)\rightarrow \R \times M_B^+ \quad b^{+}_c\mapsto
(b^{+}_c(x_0),[b^{+}_c])
\end{equation}
($\pi^+_{x_0}$ may be non-continuous, see Section \ref{topcaubo}).
Identifying each element of $B^+(M)$ with its image under
$\pi^+_{x_0}$, we can establish the following identification at a
point set level:
\[
\begin{array}{lll} \hat{V} &= & B^{+}(M)\cup \{\infty\}=(\R\times M_B^+)\cup \{i^{+}\}=(\R\times (M\cup \partial_{B}^{+}M))\cup \{i^{+}\}\\
&=&(\R\times (M\cup \partial_{C}^{+}M\cup \partial_{\cal B}^{+}M))\cup \{i^{+}\}\\
&=&V\cup (\R\times\partial_{C}^{+}M)\cup (\R\times\partial_{\cal
B}^{+}M)\cup \{i^{+}\}, \\ & & \\ \hat{\partial}V
&=&(\R\times\partial^{+}_B M)\cup \{i^{+}\}=
(\R\times\partial_{C}^{+}M)\cup (\R\times\partial_{\cal B}^{+}M)\cup \{i^{+}\}.
\end{array}
\]
This suggests a structure of cone for $\hat{\partial}V$ described
below (Defn. \ref{defcono}).
\smallskip

We proceed analogously for the past c-boundary. Consider a
past-directed timelike curve $\gamma:[\alpha,-\Omega)\rightarrow
V$ parametrized as $\gamma(t)=(-t,c(t))$. Then:
\[
\begin{array}{lll}
I^{+}[\gamma] & = & \{(t',x')\in V: \gamma(t)\ll (t',x') \hbox{ for some }t\in [\alpha,-\Omega) \}\\
& = &\{(t',x')\in V: -t<t'-\dpl(c(t),x') \hbox{ for some }t\in [\alpha,-\Omega)\}\\
 &= &\{(t',x')\in V: t'>\lim_{t\rightarrow -\Omega}(-t+\dm(x',c(t)))\}\\
 & = &\{(t',x')\in V: t'>b_{c}^{-}(x')\}.
\end{array}
\]
\begin{proposition}\label{caractIF} A subset $F\subset V$ (resp. $F\subsetneq V$) is an IF if and only if $F=F(b^{-}_{c})$ for
some Busemann function $b_{c}^{-}$ (resp. for some finite Busemann
function $b^{-}_{c}\in B^{-}(M)$).
\end{proposition}
Defining a bijection $\pi^-_{x_0}$ analogous to $\pi^+_{x_0}$
above, we also obtain the point set identifications:
\[
\begin{array}{lll}
\check{V}&=&V\cup (\R\times\partial_{C}^{-}M)\cup (\R\times\partial_{\cal
B}^{-}M)\cup \{i^{-}\} \\ \check{\partial}
V&=&(\R\times\partial_{C}^{-}M)\cup (\R\times\partial_{\cal B}^{-}M)\cup \{i^{-}\}.
\end{array}
\]
\begin{convention}\label{coco}
{\em (A) When a point  $x_0\in M$ has been chosen as above, any
set $P\in \hat{\partial} V\setminus \{i^+\}$ in the future causal
boundary will be determined univocally by a pair $(\Omega,x^+)\in
\R\times \partial_B^+ M$, according to (\ref{ezzz}).
Analogously, any terminal set $F\in \check{\partial} V\setminus
\{i^-\}$ in the past causal boundary will be represented by a pair
$(\Omega,x^-)\in \R\times\partial_B^-M$.

(B)  According to (\ref{equuu2}) and the conventions in Subsection
\ref{s525}, any non-proper forward or backward Busemann function
$b^\pm_c$ determines univocally a pair $p=(\Omega,x^\pm)\in\R\times
\partial_{C}^{\pm}M$ such that
$$
b_c^\pm(\cdot)= d_{p}^{\pm}(\cdot):=\Omega\mp
d^{\pm}(\cdot,x^\pm)=\Omega\mp d^\pm_{x^\pm}(\cdot).
$$
(By Prop. \ref{Srelat}, the domain of these functions can be
extended to points of the form $\R\times (\partial_{C}^{+}M\cup
\partial_{C}^{-}M)$.) Therefore, we can write

\begin{equation}\label{eident}
B^{\pm}(M)\setminus \mathcal{B}^{\pm}(M) \equiv \R\times
\partial_C^{\pm}M
\end{equation}
independently of the choice of any point $x_0$ as above.

(C) In what follows, we will use the intrinsic decomposition
(\ref{eident}) explained in (B), except when a point $x_0$ has
been chosen explicitly, and the decomposition (\ref{ezzz}) is
used. Thus,  the pair $(\Omega,x^{\pm})$  will be considered as
explained in (A) (for all the boundary points but $i^{\pm}$) or in
(B) (for points associated to the Cauchy boundary) depending on
the context. In each case, we will say that the TIP $P$ or the TIF
$F$ is {\em associated to} the pair $(\Omega,x^{\pm})$. Moreover,
as  the projection (\ref{ezzz}) on the second factor is
independent of the chosen $x_0$, we can also say that
$(\Omega,x^{\pm})$ is {\em associated to} $x^\pm$, and that
$(\pi^\pm)^{-1}_{x_0}(x^\pm)$ is the {\em line} on $x^{\pm}$.}
\end{convention}
In order to express previous results more compactly, we introduce
the following definition:
\begin{definition}\label{defcono}
For  any topological space $T$,  the {\em forward cone with base}
$T$ is the quotient topological space $((-\infty,\infty]\times
T)/\sim$ where the unique non-trivial identifications are
$(\infty,x)\sim (\infty,x')$ for all $x,x'\in T$. The class of
$(\infty,x)$ will be called the
{\em apex} of the cone. 

For notational
convenience, we will consider when necessary the {\em backward cone}, obtained by means of
an analogous construction for
 $[-\infty,\infty)\times T$.
\end{definition}
The results in this section have been obtained at a point set
level, and can be summarized as follow:
\begin{proposition}\label{pepe} Let $(V,g)$ be a standard stationary spacetime.
Fixing a point $x_0\in M$, the future (resp. past) c-completion
$\hat{V}$ ($\check{V}$) is identified naturally as a point set
with a forward (resp. backward) cone with basis $M_B^+$ (resp.
$M_B^-$) and apex $i^+$ (resp. $i^-$).

Moreover, the future (resp. past) c-boundary $\hat{\partial} V$
(resp. $\check{\partial} V$) is formed by a cone with basis
$\partial^{+}_{B}M=\partial^{+}_{C}M\cup \partial^{+}_{{\cal B}}M$
(resp. $\partial^{-}_{B}M=\partial^{-}_{C}M\cup
\partial^{-}_{{\cal B}}M$) and apex $i^{+}$ (resp. $i^{-}$).
\end{proposition}
The extension of this result to a topological level will be
studied in Section \ref{topcaubo}.

\subsection{The (total) c-boundary as a point set}\label{s6.3}

Next, we are going to construct the (total) c-boundary $\partial
V$ from the partial (future and past) ones. To this aim, it is
essential to study which TIPs and TIFs are S-related, according to
Section 2 (recall (\ref{eSz}), (\ref{eSz2})).

Consider a timelike curve
$\gamma_{\pm}:[\alpha,\pm\Omega)\rightarrow V$ parametrized as
$\gamma_{\pm}(t)=(\pm t,c(t))$ with $c\in C^{\pm}(M)$. If
$\pm\Omega=\infty$, the terminal set associated to $\gamma_{\pm}$
is S-related to the empty set, since the common future/past is empty (see Prop. \ref{k}). So,
we will assume $\pm\Omega<\infty$. The following result determines
the common future/past of $\gamma_{\pm}$ in terms of $\Omega\in\R$
and the limit point of $c$ in $M_{C}^{\pm}$.
\begin{lemma}\label{l}
Consider a timelike curve
$\gamma_{\pm}:[\alpha,\pm\Omega)\rightarrow V$ parametrized as
$\gamma_{\pm}(t)=(\pm t,c(t))$, with $\pm\Omega<\infty$, $c\in
C^{\pm}(M)$ and $c(t)\rightarrow x^\pm\in M^\pm_C$ with the
extended distance topology (Conv. \ref{conv}) and denote
$p=(\Omega,x^\pm)$ (Conv. \ref{coco} (B)). Then:
\[
\begin{array}{l}
I^-[\gamma_+]=\{(t',x')\in V: t'<d_p^+(x'):=\Omega-\dpl(x',x^+)\}=P(d_{p}^{+})\\
\uparrow I^-[\gamma_+]=\{(t',x')\in V:
t'>d_p^-(x'):=\Omega+\dm(x',x^+)\}=F(d_{p}^{-})
\end{array}
\]
\[
\left(\begin{array}{l}
I^+[\gamma_-]=\{(t',x')\in V: t'>d_{p}^-(x'):=\Omega+\dm(x',x^-)\}=F(d_{p}^{-})\\
\downarrow I^+[\gamma_-]=\{(t',x')\in V:
t'<d_{p}^+(x'):=\Omega-\dpl(x',x^-)\}=P(d_{p}^{+})
\end{array}
\right).
\]
\end{lemma}
{\it Proof.} The expressions for $I^-[\gamma_+]$ are
straightforward from (\ref{eimas}), and we focus on:
\[
\begin{array}{rl}
\uparrow I^{-}[\gamma_{+}] & =I^{+}[\{(t',x')\in V:\gamma_{+}(t)\ll (t',x') \;\;\hbox{for all}\;\; t\}] \\
& =I^{+}[\{(t',x')\in V:\dm(x',c(t))<t'-t \;\;\hbox{for all}\;\; t\}] \\
& =I^{+}[\{(t',x')\in V:t+\dm(x',c(t))<t'
\;\;\hbox{for all}\;\; t\}].
\end{array}
\]
Note that the function $t\mapsto t+\dm(x',c(t))$  increases
monotonically. In fact, since $c\in C^+(M)$, and so,
$F^+(\dot{c})< 1$, it is $$\dm(x',c(t_1))-\dm(x',c(t_2))\leq
\dpl(c(t_1),c(t_2))<t_2-t_1\quad\hbox{for}\;\; t_1<t_2.$$
Therefore,
\[
\begin{array}{rl}
\uparrow I^{-}[\gamma_{+}]
&
=I^{+}[\{(t',x')\in V:\lim_{t\nearrow
\Omega}(t+\dm(x',c(t)))\leq t'\}] \\ & =\{(t',x')\in
V:t'> d_{p}^{-}(x')\},
\end{array}
\]
the last equality for $p=(\Omega,x^+)$, recalling that
length$_+(c)<\infty$, $c(t)\rightarrow x^+$ and Prop.
\ref{Srelat}. Finally, the expressions for $\gamma_{-}$ are
deduced similarly. \cvd
\smallskip

\begin{remark}\label{remprocau}
{\em Observe that Propositions \ref{caractIP}, \ref{caractIF} and
Lemma \ref{l} remain true if we extend the set $C^{\pm}(M)$ with
curves of velocity lower or equal to one (recall Remark
\ref{remm1}). Therefore, there is no difference if we consider
timelike or causal curves for the computation of the c-completion
 i.e. the c-boundary of any standard stationary spacetime is {\em
properly causal} (any inextensible future-directed causal curve
will also have an endpoint), according to {\cambios Defn. 3.33 and Th. 3.35}
in \cite{FHS0}.}
\end{remark}

With previous lemma at hand, the following result can be
established.
\begin{proposition}\label{po} If $P\in\hat\partial V$, $F\in\check\partial V$ are represented by the same point
$(\Omega,x^s)\in \R\times \partial_{C}^{s} M$ then they are
S-related, $P\sim_{S} F$. In particular, they define a unique
element $(P,F)\in\partial V$.
\end{proposition}
{\it Proof.} Just observe that, from previous proposition,
$$P=\downarrow F, \qquad F=\uparrow P, $$
which provides the desired result. \cvd

\vspace{1mm}

This {\cambios proposition} establishes identifications between
those TIPs and TIFs associated to the same point of $\R\times
\partial_{C}^{s}M$. Examples \ref{ex1'} and \ref{ex2'} illustrate situations
where additional identifications between terminal sets appear;
concretely, between TIPs and TIFs associated to points in $\bcaum
\setminus \bcaums$ and $\bcaumr\setminus\bcaums$, respectively.
Remarkably, these identifications cannot appear in static spaces
as $F^+=F^-$ (see also \cite{AF}), and so, they are characteristic
of the standard stationary case. In the following proposition we
provide some restrictions for these additional identifications.
\begin{proposition}
Let $P,F$ be two terminal sets associated to
$(\Omega^+,x^+)$,
$(\Omega^-,x^-)$ resp., where $\Omega^+,\Omega^-\in \R$ and $x^\pm\in
\partial_C^\pm M\setminus \partial_C^sM$. If $P\sim_S F$, then
$\Omega^{-}-\Omega^+=\dm(x^-,x^+)$.
\end{proposition}
{\it Proof.} 
Assume by contradiction
that $\Omega^{-}-\Omega^+\neq \dm(x^-,x^+)$ and consider the cases: 


(A) $\Omega^{-}-\Omega^+<\dm(x^-,x^+)$. Take
$\Omega^{-}<t<\Omega^{+}+\dm(x^-,x^+)$ and consider a sequence
$\{(t,x_n)\}_n\subset V$ such that $x_n\rightarrow x^-$ with the
topology associated to the extended distance $\dm$. Then:
$$\lim_n (\Omega^{-}+\dm(x_n,x^-))=\Omega^{-}<t<\Omega^{+}+\dm(x^-,x^+)=\lim_n (\Omega^{+}+\dm(x_n,x^+)).$$
Therefore, if we take $n_0\in \N$ such that
$\Omega^{-}+\dm(x_{n_0},x^-)<t<\Omega^{+}+\dm(x_{n_0},x^+)$, then
$(t,x_{n_0})\in F$ and $(t,x_{n_0})\notin \uparrow P$, in
contradiction with $P\sim_S F$.

(B) $\Omega^{-}-\Omega^+>\dm(x^-,x^+)$. Take
$\overline{\Omega}=\Omega^{+}+\dm(x^-,x^+)$ and denote by
$\overline{F}$ the TIF associated to $(\overline{\Omega},x^-)$. As
$d^-_{(\Omega^{-},x^-)}>d^-_{(\overline{\Omega},x^-)}$, one has
$F\varsubsetneq \overline{F}$. Therefore:
\[\begin{array}{lll} (t,x)\in \overline{F} & \Rightarrow & t>\overline{\Omega}+\dm(x,x^-) \\ & \Rightarrow & t>\Omega^{+}+\dm(x^-,x^+)+\dm(x,x^-)\geq \Omega^{+}+\dm(x,x^+) \\ & \Rightarrow &(t,x)\in \uparrow P.\end{array}\]
Hence $(F\varsubsetneq)\overline{F}\subset \uparrow P$, in
contradiction to the maximality of $F$ in $\uparrow P$ required
from $P\sim_S F$. \cvd

\smallskip

Summing up:
\begin{theorem}\label{tiden}
Let $V=\R\times M$ be a standard stationary
spacetime. Then, those pairs $(P,F)$ in the causal boundary $\partial V$ which have $P\neq \emptyset$, satisfy $P=P(b_{c}^+)$ for some $c\in C^+(M)$, and:
\begin{itemize}
\item[(a)] If $b_c^+\equiv \infty$ then $P=V, F=\emptyset$.
\item[(b)] If $b_{c}^+\in {\cal B}^+(M)$
then $F=\emptyset$. \item[(c)] If $b_c^+\in
B^+(M)\setminus {\cal B}^+(M)(\equiv\R\times
\partial_C^+M)$, then $b_c^+=d_p^+$ with $p=(\Omega^+,x^{+})\in \R\times
\partial_C^+M$ and $F\subseteq F(d_{p}^-)$. In this
case, there are two exclusive possibilities:
\begin{itemize}
\item[(c1)] either $F=\emptyset$, \item[(c2)] or $F=F(d_{p'}^-)$
with $p'=(\Omega^-,x^{-})\in \R\times \partial_C^-M$ and satisfying
$\Omega^- -\Omega^+=\dm(x^{-},x^{+})$ (in this case, $p'$ is not
necessarily unique).
\end{itemize}
Moreover, if $x^{+}\in \partial_C^s M$, then (c2) holds with
$p'=p$, $\uparrow P=F(d_{p}^{-})$ and $P$ is univocally S-related
with $F=F(d_p^-)$.
\end{itemize}
A dual result holds for pairs $(P,F)\in \partial V$ with $F\neq
\emptyset$.
\end{theorem}

Observe that, for $P\neq \emptyset$, previous theorem determines a
unique $F$ such that $P\sim_S F$ except in the particular case
when $P=P(d_p^+)$ and $p=(\Omega^+,x^{+})\in \R\times
(\partial_C^+M\setminus \partial_C^sM)$. The presence of these
exceptions is related with the boundedness of the quasi-distance
$\dqm$ associated to the generalized metric space $(M,\dpl)$
(recall Remark \ref{motrem}). In fact, these exceptions will not
appear if we suppose that $\dqm$ is a generalized distance:
%
\begin{corollary}\label{corolarionew}
If $(M,\dpl)$ is a generalized metric space with either forward or
backward {\ras}, then the {\cambios c-completion} of $V$ is simple
as a point set (according to Def. \ref{simpletop}). In particular,
this happens if $d_Q^+$ is a generalized distance.
\end{corollary}
{\it Proof.} The proof follows by analyzing the case {\it (c2)} in
Theorem \ref{tiden}. In fact, if $P=P(d_p^+)$ with
$p=(\Omega,x^+)\in \R\times (\partial_C^+M\setminus\partial_C^s
M)$, then $\uparrow P=\emptyset$ (recall Lemma \ref{l} and
Definition \ref{dtotsep}). For the last assertion recall Corollary \ref{genersym2}.\cvd


\smallskip

Finally, recall that each point $x^\pm = [b_c^\pm]\in
\partial^\pm_B M$ determines in a natural way a {\em line}
$\{b_c^\pm+k: k\in\R\}$ in $\hat{\partial} V$, $\check{\partial}
V$ (Convention \ref{coco} (C)) and then a (point set) cone. If $V$
has a {\cambios c-completion} which is simple as a point set, such lines also
yield a simple picture of the c-boundary: $\partial V$ is composed
by  all the lines in $\hat{\partial} V$ and $\check{\partial} V$.
But, when the point in $\partial^\pm_B M$ belongs to
$\partial_C^sM$, then the lines for $\hat{\partial}V$ and
$\check{\partial}V$ must be identified as a single one ---they
represent $S$-paired elements. If the {\cambios c-completion} is not
simple, the pairings are not so trivial. However, the structure of
lines can be also defined, as we do  now for future referencing.
\begin{definition}\label{deflin}
Consider a point $(P,F)\in \partial V$ with $\emptyset\neq
P=P(b_c^+)$. Then, we define the line over
$(P,F)$, which will be denoted
$L(P,F)$, as:
\begin{itemize}
\item If $F=\emptyset$ then
$L(P,\emptyset):=\{(P',\emptyset):P'=P(b_c^+ +k),\; k\in \R\}$.
\item If $F\neq\emptyset$, and so $P=P(d_p^+), F=F(d_{p'}^-)$
(recall Theorem \ref{tiden}), then $$L(P,F):=\{(P',F'):P'=P(d_p^+
+k), F'=F(d_{p'}^- +k), k\in \R\}. $$
\end{itemize}
A dual definition is assumed for  $(P,F)\in
\partial V$ with $\emptyset\neq F=F(b_{c}^{-})$.
\end{definition}

\begin{remark}{\em
Of course, if $P\neq \emptyset\neq F$, then the two definitions of
line agree. Moreover, from the definition of the $S$-relation, if
$P\sim_S F$ with $P=P(d_{p}^{+})$, $F=F(d_{p'}^{-})$, then
necessarily $P'\sim_S F'$ with $P'=P(d_p^+ +k)$, $F'=F(d_{p'}^-
+k)$ and $k\in \R$. So, the definition of line for $P\neq
\emptyset$ is consistent.
 }
\end{remark}
Finally, we provide a pair of examples which illustrate the
non-obvious possibilities in Theorem \ref{tiden}.

\begin{example}\label{ex1'}{\em
(The possibility  $p'\neq p$ in case {\it (c2)} of Theorem
\ref{tiden} can hold.) Consider a standard stationary spacetime as
explained in Figure \ref{fig2}. Recall that this example satisfies
$$\partial_{B}^{+}M=\partial_{C}^{+}M=\{z^{+}\},\quad
\partial_{B}^{-}M=\partial_{C}^{-}M=\{z^{-}\}\quad\hbox{and}\quad
0<\dpl(z^{+},z^{-})<\infty.$$  These properties imply that, if $P$
is a TIP associated to $z^+$, then $\uparrow P$ is non-empty
(recall Lemma \ref{l}), and one can find $S$-related non-empty
sets as follows. Take a maximal IF, $F$, into $\uparrow P$. Notice
that $F$ must be terminal: otherwise, any chain generating $P$,
and so converging to $F=I^{+}(p)$ with the past chronological
topology (recall that $F$ is maximal in $\uparrow P$), must also
converge to $p$ with the topology of the manifold (since both
topologies coincide on $V$), in contradiction to the terminal
character of $P$. Hence, the TIF $F (\neq \emptyset )$ is
associated to $p^-=(\Omega^-,z^-)$ for some $\Omega^-$. As
$P\subset \downarrow F$, let $P'\supset P$ be some maximal IP into
$\downarrow F$. Analogously, $P'$ is associated to
$p^+=(\Omega^+,z^+)$ for some $\Omega^+$. By construction,
$P'\sim_S F$ with $p^+\neq p^-$.}
\end{example}
\begin{example}\label{ex2'}{\em (The suggested non-uniqueness of $p'$ in Theorem \ref{tiden} {\it
(c2)} can hold.) Consider a standard stationary spacetime
resulting from sticking two spaces as in previous example (see
Figure \ref{fig3}). This
spacetime is defined in such a way the following properties hold:

$$\partial_{C}^{+}(M)\equiv \{z_{1}^{+},z_{2}^{+}\},\;
\partial_{C}^{-}(M)\equiv \{z_{1}^{-},z_{2}^{-}\},\;
\dpl(z_1^+,z_1^-)=\dpl(z_2^+,z^-_2)=r$$ with $0<r<\infty$ and then
$P\sim_{S}F$, $P'\sim_{S}F$, $P\sim_{S}F'$, $P'\sim_{S}F'$, where
\[
\begin{array}{ll}
P=\{(t,x):t<d^{+}_{(\Omega,z_{1}^{+})}(x)\} &
 F=\{(t,x):t>d^{-}_{(\Omega+r,z_{1}^{-})}(x)\} \\
P'=\{(t,x):t<d^{+}_{(\Omega,z_{2}^{+})}(x)\} &
F'=\{(t,x):t>d^{-}_{(\Omega+r,z_{2}^{-})}(x)\}.
\end{array}
\]
So, we have the existence of four terminal sets, $P,P',F,F'$
satisfying
$$(P,F),(P,F'),(P',F),(P',F')\in \partial V.$$

\smallskip

\begin{figure}

\centering \ifpdf
  \setlength{\unitlength}{1bp}%
  \begin{picture}(410.68, 129.39)(30,0)
  \put(0,0){\includegraphics{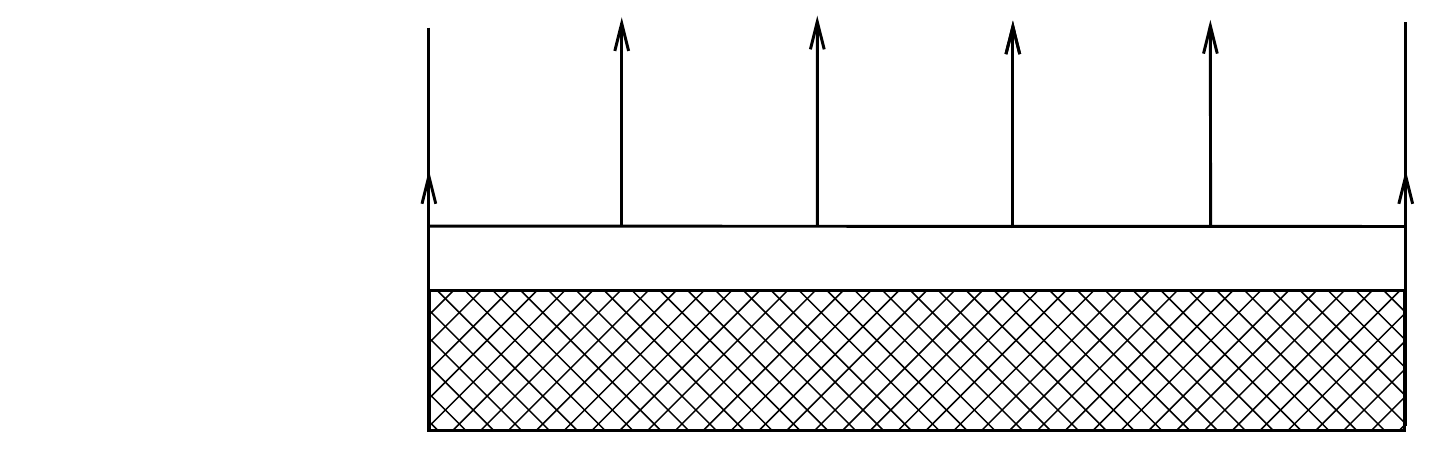}}

  \put(231.89,55.07){\fontsize{9.72}{11.66}\selectfont -3}
  \put(288.67,55.07){\fontsize{9.72}{11.66}\selectfont 3}
  \put(6.13,104.47){\fontsize{9.72}{11.66}\selectfont $\partial^{+}_{C} M\equiv \{z^{+}_{1},z^{+}_{2}\}$}

\put(5.67,86.00){\fontsize{9.72}{11.66}\selectfont
$\partial^{-}_{C} M\equiv \{z^{-}_{1},z^{-}_{2}\}$}

\put(5.67,65.23){\fontsize{9.72}{11.66}\selectfont
$\partial_{C}^s M=\emptyset$}

  \put(346.84,55.07){\fontsize{9.72}{11.66}\selectfont 9}
  \put(393.84,55.07){\fontsize{9.72}{11.66}\selectfont 12}
  \put(175.10,55.07){\fontsize{9.72}{11.66}\selectfont -9}
  \put(125.10,55.07){\fontsize{9.72}{11.66}\selectfont -12}
  \put(183.41,110.48){\fontsize{9.72}{11.66}\selectfont $c_{1}$}
  \put(239.74,111.40){\fontsize{9.72}{11.66}\selectfont $c_{2}$}

\put(295.60,111.86){\fontsize{9.72}{11.66}\selectfont $c_{3}$}
  \put(351.46,111.40){\fontsize{9.72}{11.66}\selectfont $c_{4}$}
  \end{picture}%
\else
  \setlength{\unitlength}{1bp}%
  \begin{picture}(410.68, 129.39)(30,0)
  \put(0,0){\includegraphics{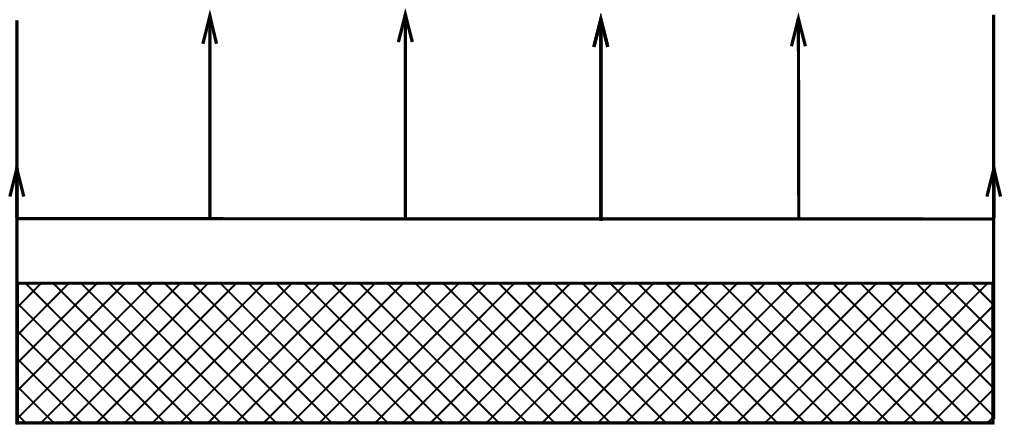}}

  \put(231.89,55.07){\fontsize{9.72}{11.66}\selectfont -3}
  \put(288.67,55.07){\fontsize{9.72}{11.66}\selectfont 3}
  \put(6.13,104.47){\fontsize{9.72}{11.66}\selectfont $\partial^{+}_{C} M\equiv \{z^{+}_{1},z^{+}_{2}\}$}

\put(5.67,86.00){\fontsize{9.72}{11.66}\selectfont
$\partial^{-}_{C} M\equiv \{z^{-}_{1},z^{-}_{2}\}$}

\put(5.67,65.23){\fontsize{9.72}{11.66}\selectfont
$\partial_{C}^s M=\emptyset$}

  \put(346.84,55.07){\fontsize{9.72}{11.66}\selectfont 9}
  \put(393.84,55.07){\fontsize{9.72}{11.66}\selectfont 12}
  \put(175.10,55.07){\fontsize{9.72}{11.66}\selectfont -9}
  \put(125.10,55.07){\fontsize{9.72}{11.66}\selectfont -12}
  \put(183.41,110.48){\fontsize{9.72}{11.66}\selectfont $c_{1}$}
  \put(239.74,111.40){\fontsize{9.72}{11.66}\selectfont $c_{2}$}

\put(295.60,111.86){\fontsize{9.72}{11.66}\selectfont $c_{3}$}
  \put(351.46,111.40){\fontsize{9.72}{11.66}\selectfont $c_{4}$}
  \end{picture}%
\fi \caption{\label{fig3} This example is constructed by sticking
two spaces as in Figure \ref{fig2}. In particular, $c_1,c_3$
(resp. $c_2,c_4$) determine two points $z_1^+,z_2^+\in
\partial_C^+ M$ (resp. $z_1^-,z_2^-\in \partial_C^- M$). Moreover,
the distances $d(z_i^+,z_j^-)$ are finite for $i,j=1,2$, but
$z_i^+,z_j^-\not\in
\partial_C^s M$ (in fact, $\partial_C^s M=\emptyset$).}
\end{figure}
}
\end{example}

\begin{remark}\label{r621}{\em
The examples above become also relevant for the following
discussions on the foundations of the notion of ``causal
boundary'':

(a) The idea that some TIP may be {\cambios S-paired with more than one
TIF} was stressed by Marolf and Ross \cite[Figure 3]{MR2}.
However, their example might be regarded as somewhat ``artificial'
(namely, some topological manipulations are carried out in their
construction). Example \ref{ex2'} yields the same effect without
topological manipulations.

(b) By modifying Marolf and Ross' example, Flores constructed a
new one which stressed that the set of all the S-related pairs
could yield a ``non-minimal'' completion. Again, this effect is
reproduced by the ``less artificial'' Example \ref{ex2'}. The
choice of non-minimal completions in the definition of c-boundary
is discussed in \cite{FHS0, Orl-M}.

(c) Even when the pairs $(P,F)$ with $P\neq \emptyset \neq F$ are
univocally determined by $P$ (or $F$), the highly non-trivial
identification in Example \ref{ex1'} stresses the limitations of
the conformal boundary \cite[Figure 11, 12]{FHS0}. Such
non-trivial identifications have been relevant to understand the
AdS-CFT correspondence in plane waves (see \cite{MR,FS} and
references therein).}\end{remark}

\subsection{Causality of the c-boundary}\label{s6.4}

In this section we are going to study the causal relations between
points of the c-boundary of a standard stationary spacetime.
Recall that in Section \ref{fff'} we defined the chronological
relation between points of the c-completion, but not a causal
relation which extends the relation $\leq$ on the spacetime to
pairs $(P,F),(P',F') \in \overline{V}$. For our purposes, it is
enough to consider the following obvious sufficient criterium for
$(P,F),(P',F')\in \overline{V}$ with either $P\neq\emptyset$ or
$F'\neq\emptyset$: $$P\subset P' \hbox{ and } F'\subset F
\Rightarrow (P,F)\leq (P',F').$$ In the particular case of
standard stationary spacetimes (which have a {\em properly causal}
c-boundary, see Remark \ref{remprocau}), this criterium can be
also regarded as the definition of $\leq$, but we will not discuss
this point here.
 Simply, we point out that a general definition of $\leq$ in $\overline{V}$ should take
 into account that the chronological
relation always induces a causal relation (see \cite[Defn. 2.22,
Th. 3.69]{MS}); we refer to \cite[Sect. 3.5]{FHS0} and \cite[Sect.
3.2]{MR2} for further discussions.

 We  focus on  causal relations between points of the
same line (Definition \ref{deflin}). To this aim, we introduce the
following definitions, which interpret the causal relations in our
case:
\begin{definition}\label{caurellin}
Consider a pair $(P,F)$ with $P=P(b_c^ +)$ and $F=F(b_{c'}^-)$ (or
$F=\emptyset$), and consider two arbitrary pairs $$(P(b_{c}^+
+k_1),F(b_{c'}^- +k_1)),(P(b_{c}^+ +k_2),F(b_{c'}^- +k_2))\in L(P,F).$$
Then, the line $L(P,F)$ is:
\begin{itemize}
\item {\em Timelike} if for $k_1<k_2$, $(P(b_{c}^+ +k_1),F(b_{c'}^-
+k_1))\ll(P(b_{c}^+ +k_2),F(b_{c'}^- +k_2))$. \item {\em Horismotic}
if for $k_1<k_2$, $(P(b_{c}^+ +k_1),F(b_{c'}^-
+k_1))\leq(P(b_{c}^+ +k_2),F(b_{c'}^- +k_2))$, but they are not
chronologically related. \item {\em Locally
horismotic
} if, for each $k_1$, there exists $\delta>0$ such that for any
$k_2\in \R$ with $|k_2-k_1|<\delta$, $(P(b_{c}^+ +k_1),F(b_{c'}^-
+k_1))$ and $(P(b_{c}^+ +k_2),F(b_{c'}^- +k_2))$ are causally but
not chronologically related.
\end{itemize}
\end{definition}
Note that, for a curve in a spacetime, previous definitions of the
causality of a line correspond with being timelike, lightlike with
no cut points and lightlike, resp. Moreover, from the stationary
character of the metric, $\delta$ in the definition of locally
horismotic can be chosen independently of $k_1$.


The causal relations for the c-boundary of a standard stationary
spacetime are summarized in the following proposition:
\begin{proposition}\label{pcaurellin}
Consider a point $(P,F)\in \partial V$  with $P\neq \emptyset$ and
 its associated line $L(P,F)$. Then:
\begin{itemize}
\item[(i)] If $F=\emptyset$, then $L(P,\emptyset)$ is horismotic.
\item[(ii)] If $F\neq \emptyset$, and so, $P=P(d_p^+)$,
$F=F(d_{p'}^-)$ (necessarily, with $p=(\Omega,x^+)\in \R\times
\partial_C^+ M$, $p'=(\Omega+\dpl(x^+,x^-),x^-)\in \R\times
\partial_C^- M$). Then:
\begin{enumerate}
\item  If $p=p'$, $L(P,F)$ is timelike. \item If $p\neq p'$,
$L(P,F)$ is locally horismotic.
\end{enumerate}
\end{itemize}
\end{proposition}
{\it Proof.}
(i) Consider two points in the line, $(P(b_c^+
+k_1),\emptyset),(P(b_c^+ +k_2),\emptyset)\in L(P,F)$. If
$k_1<k_2$ then $P(b_c^+ +k_1)\subset P(b_c^+ +k_2)$ and $\emptyset
\cap P(b_c^+ +k_2)=\emptyset$. So, $(P(b_c^+ +k_1),\emptyset)$ is
causally related with $(P(b_c^+ +k_2),\emptyset)$, but it is not
chronologically related, i.e. the line is horismotic.

(ii) If $F\neq \emptyset$, and so, $P=P(d_p^+), F=F(d_{p'}^-)$,
reasoning as in previous paragraph we can check that $P(d_p^+
+k_1)\subset P(d_p^+ +k_2)$ and $F(d_{p'}^- +k_2)\subset F(d_{p'}^-
+k_1)$ for $k_1<k_2$. Then, $(P(d_p^+ +k_1),F(d_{p'}^- +k_1))\leq (P(d_p^+
+k_2),F(d_{p'}^- +k_2))$. Now:
\begin{itemize}
\item[1.] For $p=p'=(\Omega,x^s)\in \R\times \partial_C^s M$,
observe that $(\overline{\Omega},\overline{x})\in P(d_p^+
+k_2)\cap F(d_{p}^- +k_1)$, with
$\overline{\Omega}=\Omega+(k_1+k_2)/2$ and $\overline{x}\in M$
such that $d^{\pm}(\overline{x},x^s)<(k_2-k_1)/2$ (recall that
$x^s\in\partial_{C}^{s}M$). In conclusion, both pairs are also
chronologically related. \item[2.]
For $p\neq p'$, the line $L(P,F)$ is locally horismotic. In fact,
if two points of the line were chronologically related, then we
would have $(t,x)\in P(d_p^+ +k_2)\cap F(d_{p'}^- +k_1)$, for some
$k_2>k_1$. Therefore, $t<\Omega+k_2$ and $t>\Omega+\dpl(x^+,x^-)+k_1$,
which is only possible if $k_2-k_1>\dpl(x^+,x^-)$. \cvd
\end{itemize}

\begin{remark}{\em
This last point 2 in the previous proof is sharper in the determination of
the locally horismotic character of the line $L(P,F)$ for $p\neq
p'$. In fact, when the two values $k_1, k_2$ of the parameter of
the line satisfies $|k_1-k_2|\leq \dpl(x^+,x^-)$, the
corresponding points in $L(P,F)$ are not chronologically related.
So, $\dpl(x^+,x^-)$ behaves as a lower bound for the appearance of
a  ``cut point'' along $L(P,F)$. }\end{remark}

\subsection{Topology of the partial boundaries and the c-boundary}\label{topcaubo}

In Proposition \ref{pepe}, we established that the point set
structures of the future completion and the future boundary can be regarded in a natural way as a cone. In fact, for any $x_0$ in $M$, the map $\pi_{x_0}^+$
defined in (\ref{ezzz}) is bijective.

In order to extend the cone structure of $\hat{V}$ and
$\hat{\partial} V$ to a topological level, some issues on the
chronological topology appear. On one hand, if the chronological
topology differs from the pointwise one, the map $\pi^+_{x_0}$ may
be non-continuous (this happens when a sequence $\{f_n\}$ which
converges chronologically to a function $f$,  does not converge
pointwise at $x_0$, i.e. $\{f_n(x_0)\}\not\rightarrow f(x_0)$). On
the other hand, the continuity of $\pi_{x_0}^+$ at some point may
depend on the choice
of $x_0\in M$: 
observe that, in Example \ref{ex43}, the chronological and
pointwise convergences of $\{x_n\}$ differ, but $\pi_{x_0}$,
$x_0=(a,b)$, preserves the limit if, and only if, $b=1/2$ (see
Figure \ref{peine} (B) for an intuitive idea). }

At any case, the continuity will be ensured if we consider the
pointwise topology. So, we are interested in the cases when both
topologies, the chronological and the pointwise one, coincide.
With this aim, note first (compare with Theorem \ref{below}):
\begin{lemma}\label{l625}
If $M_B^+$ is Hausdorff then the pointwise and the chronological
topologies coincide on $B^+(M)$.
\end{lemma}
{\it Proof.} Assume that $\{b^+_{c_{n}}\}\subset B^{+}(M)$
converges to some $b^+_{c}\in B^{+}(M)$ with the chronological
topology, and let us prove that it also converges with the
pointwise one (the converse follows from Prop. \ref{icont}). From Lemma \ref{lee}, $B^{+}(M)$ is ULL,  thus, $\hat{L}$ is of
first order (Proposition \ref{key}) and
$b^+_{c}\in\hat{L}(\{b^+_{c_{n}}\})$. Obviously,
$\{[b^+_{c_{n}}]\}\subset M_{B}^{+}$ converges to $[b^+_{c}]\in
M_{B}^{+}$ with its own chronological topology, and so, it also
converges with the pointwise one (as $M_B^+$ is Hausdorff, Theorem
\ref{below} applies). Therefore, there exist some constants
$\{k_{n}\}\subset\R$ such that $\{b^+_{c_{n}}-k_{n}\}\subset
B^{+}(M)$ converges pointwise to $b^+_{c}\in B^{+}(M)$. Moreover,
from the definition of chronological limit (see (\ref{ecu1})),
$k_n\rightarrow 0$.
In fact, otherwise,
there are two possibilities: either $k_{n_j}>\epsilon>0$, or
$k_{n_{j}}<-\epsilon<0$, for some subsequence
$\{n_{j}\}_{j}\subset \{n\}$. In the first case
$b^{+}_{c}+\epsilon/2<\liminf b^{+}_{c_{n_{j}}}$, and so,
$\{b^{+}_{c_{n}}\}$ does not converge chronologically to
$b_{c}^{+}$, since $b^{+}_{c}+\epsilon/2$ violates the maximal
character of $b_{c}^{+}$ (recall (b) in (\ref{ecu1})). In the
second case, $b_{c}^{+}>\liminf_{n}b^{+}_{c_{n_{j}}}$, which
violates condition (a) in (\ref{ecu1}).
\cvd
\begin{theorem}\label{pconocau}
Let $V=\R\times M$ be a standard stationary spacetime. If
$M_{B}^{+}$ is Hausdorff then the map $\pi^+_{x_0}$ in
(\ref{ezzz}) is an homeomorphism for all $x_0\in M$.

So, in this case, the future causal completion
$\hat{V}=B^{+}(M)\cup \{\infty\}$ is a (topological) cone with
basis $M_{B}^{+}$ and apex $i^{+}$, and the future causal boundary
$\hat{\partial}V=B^{+}(M)\setminus V$ is also a (topological) cone
with basis $\partial_{B}^{+}M$ and apex $i^{+}$.
\end{theorem}
{\it Proof.} By previous lemma, the subset $B^{+}(M)$ of $\hat{V}$
is endowed with the pointwise convergence topology,
 and then,  $\pi^+_{x_0}$ is obviously continuous. For
the continuity of the inverse, observe that if $\{[b^{+}_{c_n}]\}$
converges pointwise to $[b^{+}_{c}]$ and $\{b^{+}_{c_n}(x_0)\}$
converges to $b^{+}_{c}(x_0)$, necessarily $\{b^{+}_{c_n}\}$
converges pointwise to $b^{+}_c$. Finally, as $\{b^+_{c_{n}}\equiv
(t_n,x^+)\}\subset \R\times M_{B}^{+}$ converges to $i^{+}$ if and
only if $t_{n}\rightarrow\infty$, we deduce that the future causal
completion $\hat{V}$ (resp. the future causal boundary
$\hat{\partial}V$) is a topological cone with base $M_B^+$ (resp.
$\partial_B^+ M$) and apex $i^+$. \cvd

\smallskip

The definition of the chronological topology, combined with
Theorem \ref{tiden} and Proposition \ref{pcaurellin}, yields a
full description of the c-completion of a stationary spacetime.
Now, Theorem \ref{pconocau} shows that the topological structure
of the partial boundaries is particularly simple when
$M_B^\pm$ is Hausdorff. Even in this case, the
topological structure of the (total) c-boundary presents some
complications related with the ambiguity appeared at the point set
level (recall Theorem \ref{tiden}). So, in order to ensure a well
behaviour of the S-relation and obtain a nice description of the
topology of the (total) c-boundary, additional conditions are
required:
\begin{theorem}\label{oo}
Let $V=\R\times M$ be a standard stationary spacetime. If
$M_{C}^{s}$ is locally compact and $\dqm$ is a generalized
distance,
then $\overline{V}$ is simple, and so, $\overline{V}$ is homeomorphic to the quotient topological space
$(\hat{V}\cup_d \check{V})/\sim_{S}$.
\end{theorem}
{\it Proof.}
According to Definition \ref{simpletop} and Proposition \ref{propsimplepunt1}
(2), it suffices to show that a pair $(P,F)\in \overline{V}$ with
$P\neq \emptyset\neq F$ is limit of a sequence $\{(P_n,F_n)\}$
whenever $P\in \hat{L}(P_n)$ {\it or} $F\in\check{L}(F_n)$. In
fact, if say $P\in\hat{L}(P_n)$, let us see that
$F\in\check{L}(F_n)$. Since $\dqm$ is generalized and
$P\neq\emptyset\neq F$, Theorem \ref{tiden} (see also Corollary
\ref{corolarionew}) implies that $P=P(d^{+}_{p})$,
$F=F(d_{p}^{-})$ for some $p=(\Omega,x^s)\in\R\times M_{C}^{s}$.
As $M_C^s$ is locally compact, $P_{n}$ has to be of the form
$P(d^{+}_{p_{n}})$, $p_{n}=(\Omega_{n},x_{n}^s)\in \R\times
M_{C}^{s}$, for $n$ big enough, and such that $x_n^s\rightarrow
x^s$ (recall Prop. \ref{lemmabusemanM}). Moreover, from Theorem
\ref{tiden}, also
$F_n=F(d_{p_n}^-)$. Since $p_n\rightarrow p$, it is $F\in
\check{L}(F_{n})$ (recall Lemma \ref{lcf} and Proposition \ref{icont}.} \cvd

\smallskip

Finally, summarizing previous results:
\begin{corollary}\label{oo2} Let $V=\R\times M$ be a standard stationary
spacetime. If $\dqm$ is a generalized distance, $M_{C}^{s}$ is
locally compact and $M_B^\pm$ are Hausdorff, then $\overline{V}$
is simple and $\hat{\partial} V, \check{\partial} V$ have the
structure of topological cones with bases $\partial_B^+ M$,
$\partial_B^- M$ (being the apexes denoted $i^+,i^-$) resp.
Therefore, the completion $\overline{V}$ is homeomorphic to the
quotient topological space $\hat{V}\cup_d \check{V}/\sim_S$, where
$\hat{V},\check{V}$ has the structure of topological cones with
bases $M_B^+,M_B^-$ (and same apexes $i^{+}$, $i^{-}$), resp.
\end{corollary}

\begin{remark}\label{r629}{\em
Summarizing,  Figure \ref{figfinal} represents the point set and
causal structures of the {\cambios c-completion} of {\cambios a}
standard stationary spacetime. The c-boundary illustrated in the
figure is simple as a point set; if this condition did not hold,
only the addition of new locally horismotic lines (as the depicted
``helix'') through pairs type $(P_3,F'_3)$ or $(P'_3,F_3)$ would
be required. The first diagram in Figure \ref{figdiagram}
summarizes the implications at these pointset and chronological
levels.

In order to ensure also the topological structure suggested by
Figure \ref{figfinal}, two conditions have to be imposed: the
Hausdorffness of $M_B^\pm$ {\cambios (to obtain the topological
structure of double cone, recall Th. \ref{pconocau}), and the
(topological) simplicity
of $\partial V$ (to ensure the good behavior of the
identifications). The second diagram in Figure \ref{figdiagram}
summarizes when these properties are ensured. Recall that: (a)
when $M_B^\pm$ is Hausdorff, both the chronological topology on
$\overline{V}$ is also Hausdorff (recall Lemma \ref{l625} and
the identifications between terminal sets and $B^\pm (M)$),
and $\overline{V}$ is simple as a point set (Proposition
\ref{propsimplepunt1}(1)),  (b) the simplicity of $\partial V$ is
ensured under the conditions in Theorem \ref{oo}. }

Particularized to the static case (when $\dqm$ is directly a
generalized distance), Theorem \ref{pconocau} culminates the study
of $\hat\partial V$   in \cite{H, FH}, and Corollary \ref{oo2}
gives the correct topological description of $\partial V$, thus
completing \cite{AF}. Recall that the  pairings for the boundary
in the static case are {\cambios trivial}. So, Figure
\ref{figfinal2} gives a full intuitive picture of the point set
and causal structures of the boundary, which is also valid at the
topological level when $M_C$ is locally compact and $M_B$
Hausdorff.}\end{remark}

\begin{figure}
\centering
\ifpdf
  \setlength{\unitlength}{1bp}%
  \begin{picture}(380.92, 327.69)(10,0)
  \put(0,0){\includegraphics{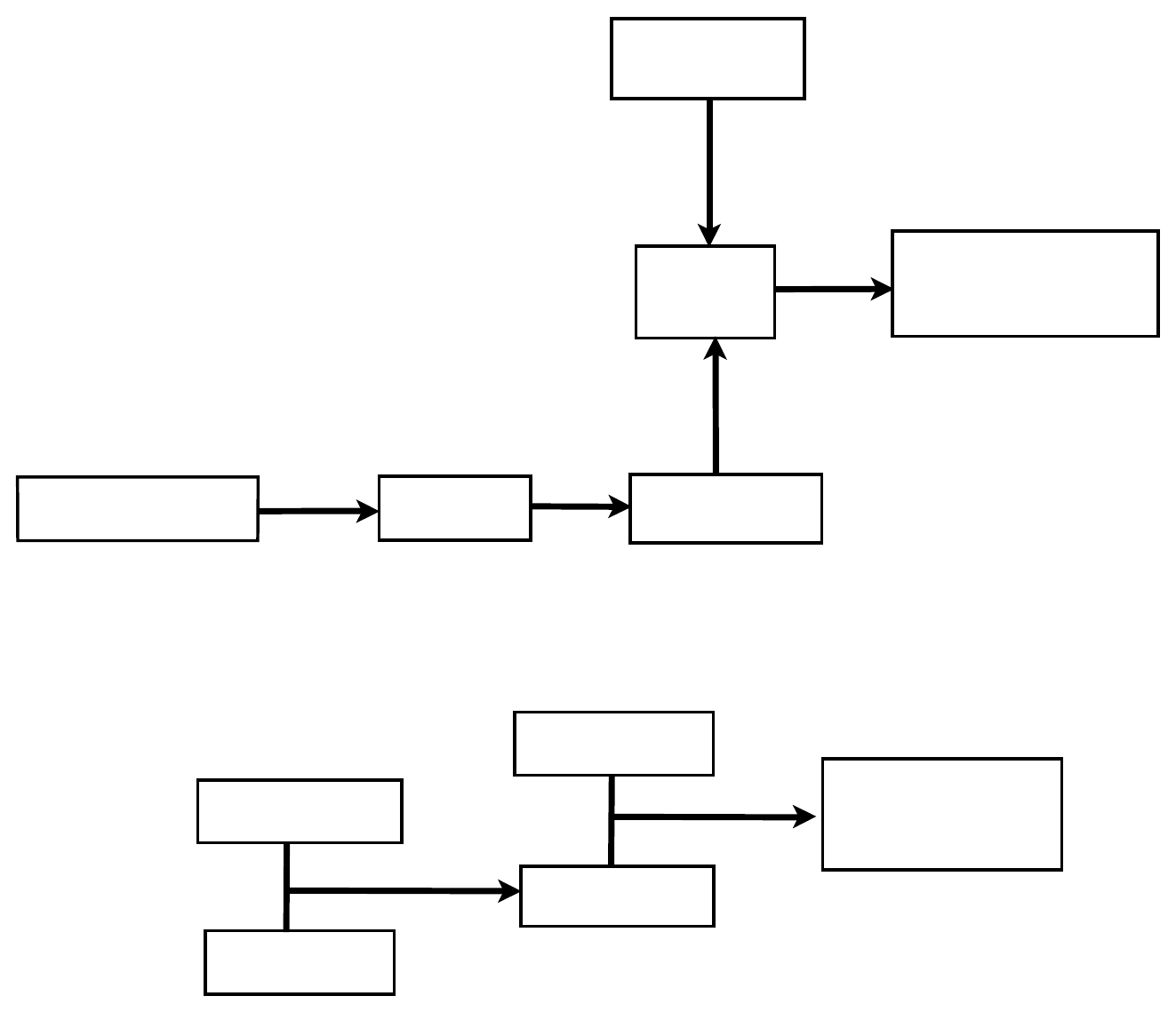}}
  \put(11.89,165.55){\fontsize{5.02}{6.02}\selectfont Riemannian distance}
  \put(125.18,163.13){\fontsize{5.02}{6.02}\selectfont $\dqm$ generalized}
  \put(214.80,165.46){\fontsize{5.02}{6.02}\selectfont Evenly paired }
  \put(205.82,308.72){\fontsize{5.02}{6.02}\selectfont $M_B^\pm$ Hausdorff}
  \put(212.39,235.68){\fontsize{5.02}{6.02}\selectfont $\overline{V}$ simple }
  \put(295.55,239.61){\fontsize{5.02}{6.02}\selectfont Point and chronological}
  \put(300.35,230.23){\fontsize{5.02}{6.02}\selectfont structure as in Fig. \ref{figfinal}}
  \put(66.98,62.74){\fontsize{5.02}{6.02}\selectfont $M_C^s$ locally compact}
  \put(74.69,14.72){\fontsize{5.02}{6.02}\selectfont $\dqm$ generalized}
  \put(173.73,85.17){\fontsize{5.02}{6.02}\selectfont $M_B^\pm$ Hausdorff}
  \put(182.88,37.88){\fontsize{5.02}{6.02}\selectfont $\overline{V}$ simple}
\put(171.88,31.88){\fontsize{5.02}{6.02}\selectfont (also topologically)}
  \put(274.91,65.68){\fontsize{5.02}{6.02}\selectfont Topology on $\overline{V}$, $\partial V$}
  \put(270.71,58.19){\fontsize{5.02}{6.02}\selectfont with two cones (Fig. \ref{figfinal})}
  \put(23.49,158.53){\fontsize{5.02}{6.02}\selectfont (Static case)}
  \put(216.46,158.38){\fontsize{5.02}{6.02}\selectfont bound. $\partial_C^\pm M$}
  \put(208.38,228.30){\fontsize{5.02}{6.02}\selectfont as a point set}
  \end{picture}%
\else
  \setlength{\unitlength}{1bp}%
  \begin{picture}(380.92, 327.69)(0,0)
  \put(0,0){\includegraphics{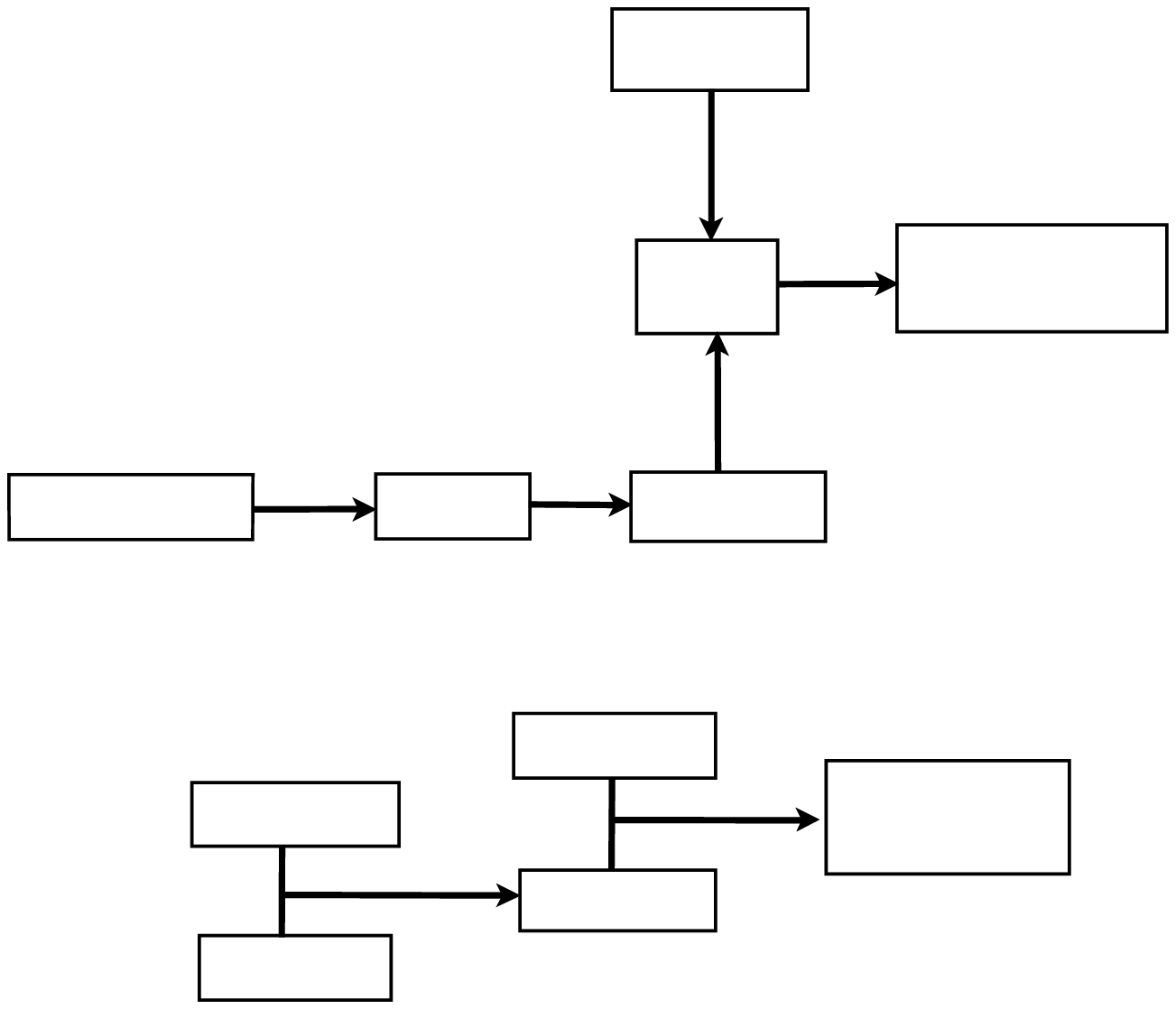}}
  \put(11.89,165.55){\fontsize{5.02}{6.02}\selectfont Riemannian distance}
  \put(125.18,163.13){\fontsize{5.02}{6.02}\selectfont $\dqm$ generalized}
  \put(214.80,165.46){\fontsize{5.02}{6.02}\selectfont Evenly paired }
  \put(205.82,308.72){\fontsize{5.02}{6.02}\selectfont $M_B^\pm$ Hausdorff}
  \put(212.39,235.68){\fontsize{5.02}{6.02}\selectfont $\overline{V}$ simple }
  \put(295.55,239.61){\fontsize{5.02}{6.02}\selectfont Point and chronological}
  \put(300.35,230.23){\fontsize{5.02}{6.02}\selectfont structure as in Fig. \ref{figfinal}}
  \put(66.98,62.74){\fontsize{5.02}{6.02}\selectfont $M_C^s$ locally compact}
  \put(74.69,14.72){\fontsize{5.02}{6.02}\selectfont $\dqm$ generalized}
  \put(173.73,85.17){\fontsize{5.02}{6.02}\selectfont $M_B^\pm$ Hausdorff}
  \put(182.88,37.88){\fontsize{5.02}{6.02}\selectfont $\overline{V}$ simple}
\put(171.88,31.88){\fontsize{5.02}{6.02}\selectfont (also topologically)}
  \put(274.91,65.68){\fontsize{5.02}{6.02}\selectfont Topology on $\overline{V}$, $\partial V$}
  \put(270.71,58.19){\fontsize{5.02}{6.02}\selectfont with two cones (Fig. \ref{figfinal})}
  \put(23.49,158.53){\fontsize{5.02}{6.02}\selectfont (Static case)}
  \put(216.46,158.38){\fontsize{5.02}{6.02}\selectfont bound. $\partial_C^\pm M$}
  \put(208.38,228.30){\fontsize{5.02}{6.02}\selectfont as a point set}
  \end{picture}%
\fi \caption{\label{figdiagram} Sufficient general conditions
which simplify the structure of the c-boundary both, as a point
set and topologically  (the arrows mean ``implies''), see Remark
\ref{r629}.\newline The character of {\em simple as a point set}
for $\overline{V}$ (which is satisfied either if $M^\pm_B$
is Hausdorff or if $V$ is static), is enough to ensure that
$\partial V$ includes a ``natural double cone'' as a point set.
However, to ensure that this double cone appears also at the
topological level, both conditions, $M^\pm_B$ is Hausdorff
and $\overline{V}$ is simple (also topologically), are imposed. In
terms of the geometry of $M$, the latter condition is ensured by
imposing both, $M_C^s$ is locally compact and $d^+_Q$ is a
generalized distance.}
\end{figure}
\begin{figure}
\centering
\ifpdf
  \setlength{\unitlength}{1bp}%
  \begin{picture}(460.52, 335.20)(0,0)
  \put(0,0){\includegraphics{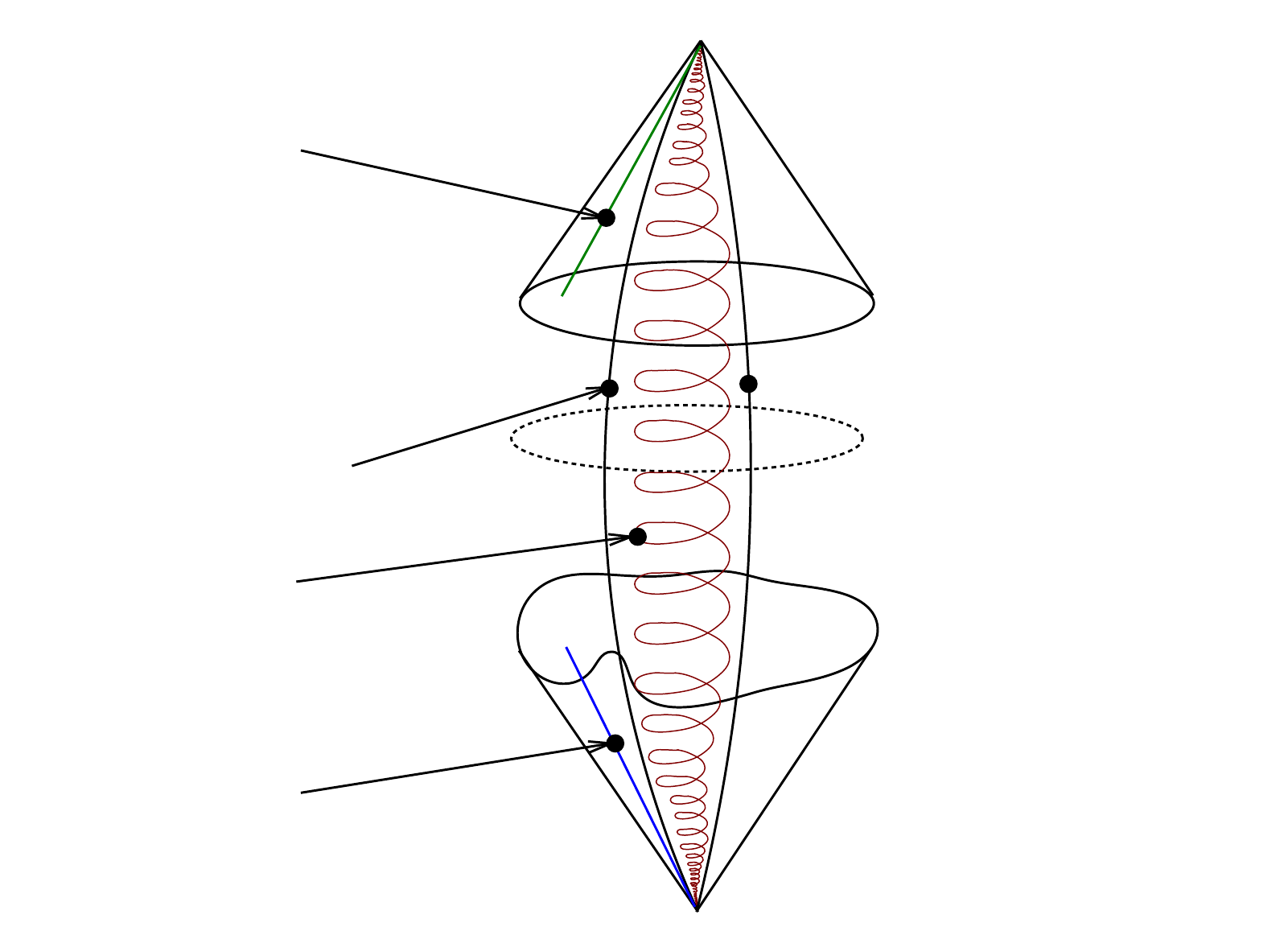}}
  \put(20.04,157.10){\fontsize{13.57}{16.29}\selectfont $(P_2,F_2): P_2,F_2$ from $\partial_{C}^s M$}
  \put(259.54,318.93){\fontsize{13.57}{16.29}\selectfont $i^{+}$}
  \put(256.23,8.59){\fontsize{13.57}{16.29}\selectfont $i^{-}$}
  \put(321.96,225.50){\fontsize{13.57}{16.29}\selectfont $\partial_{\cal B}^+ M$}
  \put(317.06,271.74){\fontsize{13.57}{16.29}\selectfont Future cone}
  \put(301.54,40.59){\fontsize{13.57}{16.29}\selectfont Past cone}
  \put(320.33,103.48){\fontsize{13.57}{16.29}\selectfont $\partial_{\cal B}^- M$}
  \put(326.30,172.75){\fontsize{13.57}{16.29}\selectfont $M$}
  \put(-10.24,36.15){\fontsize{13.57}{16.29}\selectfont $(\emptyset,F_1): F_1$ from $\partial_C^-M\setminus \partial_C^s M$}
  \put(-5.24,287.79){\fontsize{13.57}{16.29}\selectfont $(P_1,\emptyset): P_1$ from $\partial_C^+ M\setminus \partial_C^s M$ }
  \put(-60,114.94){\fontsize{13.57}{16.29}\selectfont $(P_3,F_3): P_3,F_3$ from $\partial_C^\pm M\setminus \partial_C^s M$ resp.}
  \end{picture}%
\else
  \setlength{\unitlength}{1bp}%
  \begin{picture}(460.52, 335.20)(0,0)
  \put(0,0){\includegraphics{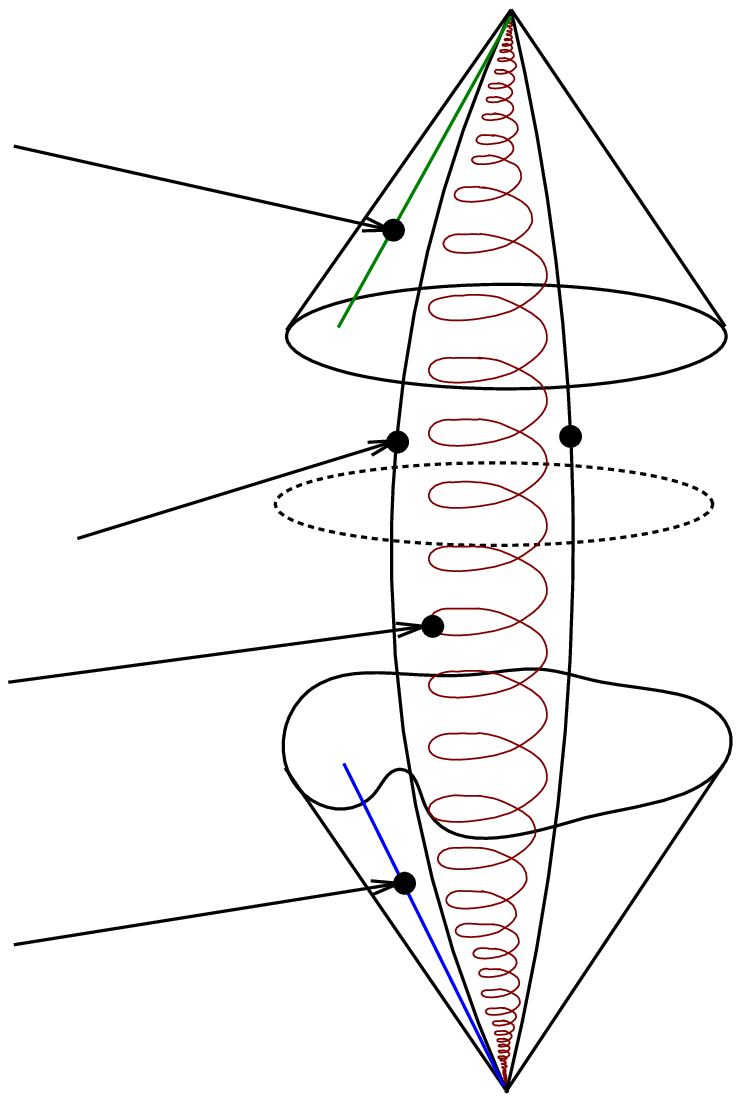}}
  \put(20.04,157.10){\fontsize{13.57}{16.29}\selectfont $(P_2,F_2): P_2,F_2$ from $\partial_{C}^s M$}
  \put(259.54,318.93){\fontsize{13.57}{16.29}\selectfont $i^{+}$}
  \put(256.23,8.59){\fontsize{13.57}{16.29}\selectfont $i^{-}$}
  \put(321.96,225.50){\fontsize{13.57}{16.29}\selectfont $\partial_{\cal B}^+ M$}
  \put(317.06,271.74){\fontsize{13.57}{16.29}\selectfont Future cone}
  \put(301.54,40.59){\fontsize{13.57}{16.29}\selectfont Past cone}
  \put(320.33,103.48){\fontsize{13.57}{16.29}\selectfont $\partial_{\cal B}^- M$}
  \put(326.30,172.75){\fontsize{13.57}{16.29}\selectfont $M$}
  \put(-10.24,36.15){\fontsize{13.57}{16.29}\selectfont $(\emptyset,F_1): F_1$ from $\partial_C^-M\setminus \partial_C^s M$}
  \put(-5.24,287.79){\fontsize{13.57}{16.29}\selectfont $(P_1,\emptyset): P_1$ from $\partial_C^+ M\setminus \partial_C^s M$ }
  \put(-60,114.94){\fontsize{13.57}{16.29}\selectfont $(P_3,F_3): P_3,F_3$ from $\partial_C^\pm M\setminus \partial_C^s M$ resp.}
  \end{picture}%
\fi \caption{\label{figfinal} Possibilities for the c-boundary of
a standard stationary spacetime:\newline The properly Busemann
boundaries $\partial_{\cal B}^\pm M$ determine the (forward,
backward) directions of $M$ at infinity. Each one yields a cone in
$\partial V$ with apex $i^+$ or $i^-$ and horismotic lines from
the apex.\newline Each point $x^+$ of the Cauchy boundary
$\partial_C^+M$ determines one (or more) lines with different
properties depending on if it generates $S$-relations (non-trivial
pairings) with one (or more) points $x^-\in
\partial_C^-M$ or not. If it does not, then $x^+$ generates a
single horismotic line [see $(P_1,\emptyset)$ in the figure;
$(\emptyset, F_1)$ represents the dual situation for
$\partial_C^-M$]. In the case that $x^+$ generates some
$S$-relation, there are two possibilities: (i) $x^+$ belongs to
$\partial_C^sM$, and so, it generates a single line, which is
timelike [see  $(P_2,F_2)$]; (ii) $x^+$ belongs to $\partial_C^+M
\backslash \partial_C^sM$ and generates one (or more) lightlike
(i.e. locally horismotic) lines [represented as the ``helix''
through $(P_3,F_3)$].}
\end{figure}

\begin{figure}
\centering
\ifpdf
  \setlength{\unitlength}{1bp}%
  \begin{picture}(250.00, 335.20)(0,0)
  \put(0,0){\includegraphics{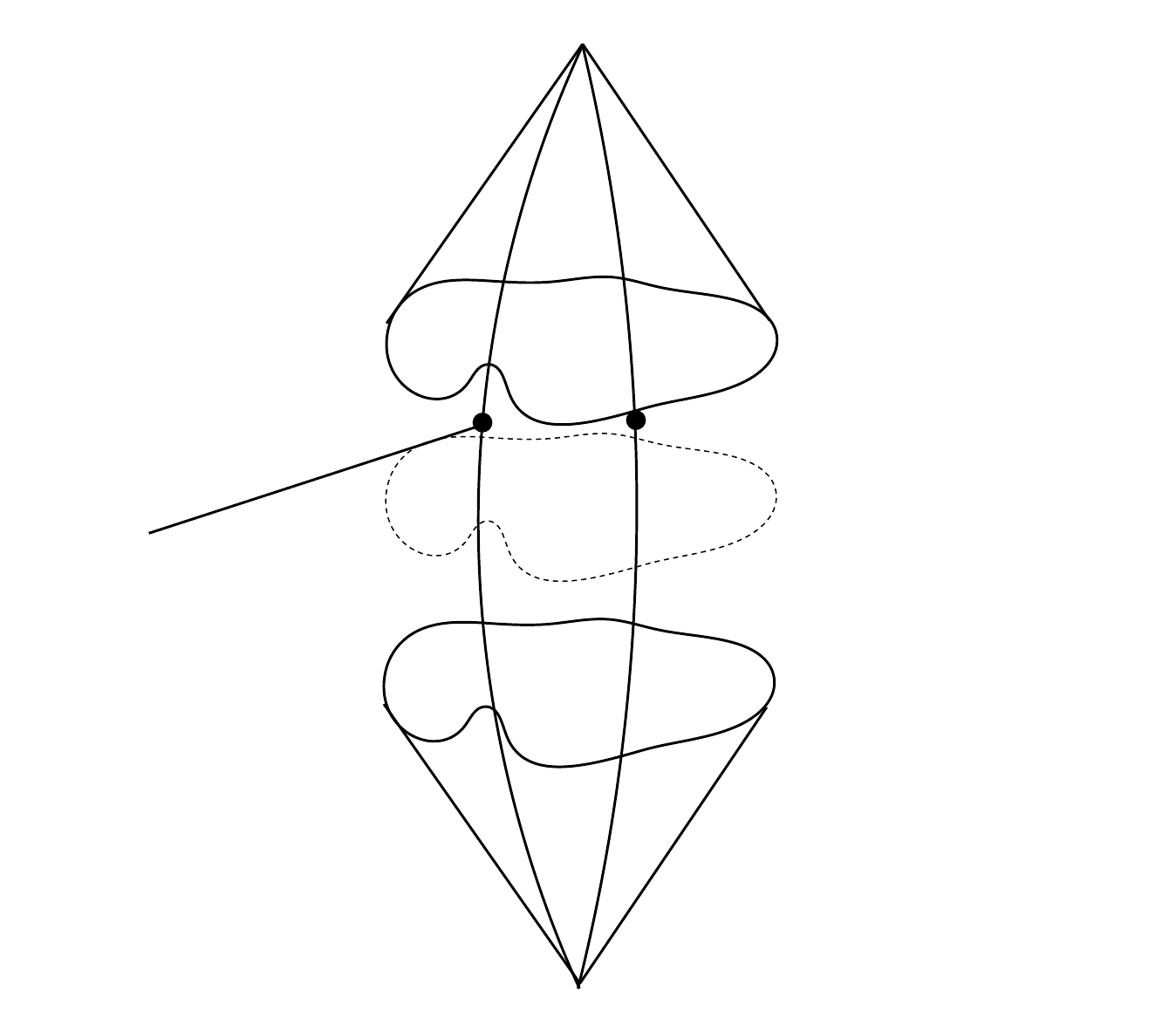}}
  \put(200.99,318.93){\fontsize{13.57}{16.29}\selectfont $i^{+}$}
  \put(197.68,8.59){\fontsize{13.57}{16.29}\selectfont $i^{-}$}
  \put(263.41,209.66){\fontsize{13.57}{16.29}\selectfont $\partial_{\cal B} M$}
  \put(258.51,271.74){\fontsize{13.57}{16.29}\selectfont Future cone}
  \put(242.99,40.59){\fontsize{13.57}{16.29}\selectfont Past cone}
  \put(261.78,103.48){\fontsize{13.57}{16.29}\selectfont $\partial_{\cal B} M$}
  \put(-30.67,145.03){\fontsize{13.57}{16.29}\selectfont $(P,F): P,F$ from $\partial_C M$}
  \end{picture}%
\else
  \setlength{\unitlength}{1bp}%
  \begin{picture}(250.00, 335.20)(0,0)
  \put(0,0){\includegraphics{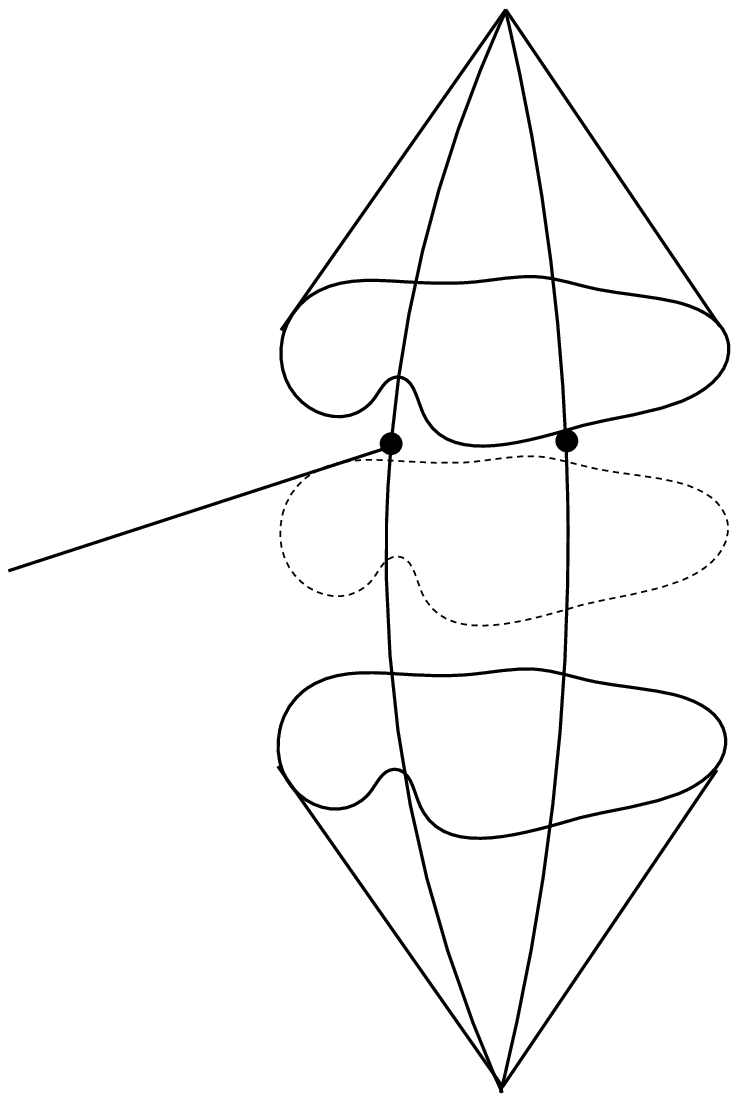}}
  \put(200.99,318.93){\fontsize{13.57}{16.29}\selectfont $i^{+}$}
  \put(197.68,8.59){\fontsize{13.57}{16.29}\selectfont $i^{-}$}
  \put(263.41,209.66){\fontsize{13.57}{16.29}\selectfont $\partial_{\cal B} M$}
  \put(258.51,271.74){\fontsize{13.57}{16.29}\selectfont Future cone}
  \put(242.99,40.59){\fontsize{13.57}{16.29}\selectfont Past cone}
  \put(261.78,103.48){\fontsize{13.57}{16.29}\selectfont $\partial_{\cal B} M$}
  \put(-30.67,145.03){\fontsize{13.57}{16.29}\selectfont $(P,F): P,F$ from $\partial_C M$}
  \end{picture}%
\fi \caption{\label{figfinal2} The static case (when the Finsler
metric on $M$ becomes Riemmanian) satisfies that $\partial_C^\pm
M=\partial_C^ s M=\partial_C M$ and $\partial_{\cal B}^\pm
M=\partial_{\cal B} M$; thus, the boundary is also symmetric under
time reversal. As a point set, the c-boundary of a standard static
spacetime is formed by two cones with bases $\partial_{\cal B} M$
and apexes $i^+,i^-$ resp., and timelike lines associated to
points of $\partial_C M$. However, to ensure this structure at the
topological level, the local compactness of $M_C$ plus the
Hausdorffness of $M_B$ must be imposed.}
\end{figure}

\subsection{Proof of Theorem \ref{th}}

Now, we are in conditions to summarize the proof of Theorem
\ref{th}.

\smallskip

\noindent {\it Proof of Theor. \ref{th}.} Assertion (1A) of the
point set structure  follows from Proposition \ref{pepe}. {\cambios
Assertion (1B) follows from Theorem \ref{tiden}.}

The assertion about the causal structure is summarized in Proposition \ref{pcaurellin}.

Finally, about the topological structure, assertions (3A) and (3B)
follow from Theorems \ref{pconocau} and \ref{oo} resp. The last paragraph is proved in Corollary \ref{oo2}.
\cvd

\section*{Acknowledgments}

The comments by one of the referees are acknowledged. The authors
are partially supported by the Spanish  grants MTM2010-18099
(MICINN)  and    P09-FQM-4496 (J. Andaluc\'{i}a), both with FEDER
funds. Also, the second-named author is supported by the Ph.D.
thesis Spanish MEC Grant AP2006-02237, \cite{He}. The third-named
author acknowledges IHES (Bures-sur-Yvette, France) a stay during
the month of August 2010, about the topics of this paper.

\end{document}